\documentclass{article}
\usepackage{graphicx}
\usepackage{amsmath,amssymb,amsthm}
\usepackage{epsfig}
\usepackage[english]{babel}
\usepackage[utf8]{inputenc}
\usepackage{mathtools}
\usepackage{subfigure}
\usepackage{enumerate}
\usepackage{hyperref}

\allowdisplaybreaks
\usepackage{geometry}
\theoremstyle{plain}
\newtheorem{thm}{Theorem}[section]
\newtheorem{lemm}[thm]{Lemma}
\theoremstyle{definition}
\newtheorem{exm}[thm]{Example}
\newtheorem{defi}[thm]{Definition}
\newtheorem{cor}[thm]{Corollary}
\newtheorem{rmk}[thm]{Remark}
\newtheorem{prop}[thm]{Proposition}
\newcommand{\defeq}{\stackrel{\text{def}}{=}}

\newcommand{\norm}[1]{\left\lVert#1\right\rVert}

\newcommand{\sqrtx}[2]{\sqrt{#1}\\ &\overline{\rule{0pt}{2.5ex}{#2 }}}
\newcommand\scalemath[2]{\scalebox{#1}{\mbox{\ensuremath{\displaystyle#2}}}}

\DeclareMathSymbol{\C}{\mathalpha}{AMSb}{"43}
\DeclareMathOperator{\supp}{supp}
\DeclareMathOperator{\diam}{diam}
\DeclareMathOperator{\Leb}{Leb}

\usepackage{mathtools}

\title{\bf On the Kunz-Souillard approach to localization for the discrete one dimensional generalized Anderson model}
\author{\textsc{ Valmir Bucaj}\footnote{The author was supported in part by NSF grant DMS--1361625.}}
\date{}
\begin{document}
\maketitle

\vskip 12pt
\begin{abstract}
\noindent We prove dynamical and spectral localization at all energies for the discrete {\it generalized Anderson model} via the Kunz-Souillard approach to localization.  This is an extension of the original Kunz-Souillard approach to localization for Schr\"odinger operators, to the case where a single random variable determines the potential on a block of an arbitrary, but fixed, size $\alpha$. For this model, we also give a description of the almost sure spectrum as a set and prove uniform positivity of the Lyapunov exponents. In fact, regarding positivity of the Lyapunov exponents, we prove a stronger statement where we also allow finitely supported distributions. We also show that for any size $\alpha$ {\it generalized Anderson model}, there exists some finitely supported distribution $\nu$ for which the Lyapunov exponent will vanish for at least one energy. Moreover, restricting to the special case $\alpha=1$, we describe a pleasant consequence of this modified technique to the original Kunz-Souillard approach to localization. In particular, we demonstrate that actually the single operator $T_1$ is a strict contraction in $L^2(\mathbb{R})$, whereas before it was only shown that the second iterate of $T_1$ is a strict contraction.   
\end{abstract}

\tableofcontents

\section{Introduction and setting}

The study of random Schr\"odinger operators is of particular importance, since such operators model disordered media (e.g. amorphous solids). In some instances, as it is the case for crystals, the structure of the solid is completely regular; that is, the atoms are distributed periodically on some lattice. Then, mathematically, in such regular crystals, the total potential that a single particle (e.g. electron) at some position in $\mathbb{R}^d$ feels is periodic with respect to the lattice at hand. Schr\"odinger operators with periodic potentials are well understood, see for example \cite{rs}, \cite{eastham}, \cite{techl1}.

 However, as it is often the case in nature, if the positions of the atoms in the solid deviate from, say, a lattice in some highly non-regular way, or if the solid is some kind of mixture of various materials, then it is natural to view the potential that, say, a single particle feels at some position, as some random quantity. Mathematically, this can be studied via Schr\"odinger operators with random potentials. So, understanding spectral properties of such operators is of great importance.

In this paper, we consider the case where the potentials of the Schr\"odinger operators are generated by independent and identically distributed random variables ($i.i.d.$). 
Specifically, the model which we study is as follows: Suppose $r:\mathbb{R}\to \mathbb{R}_{\geq 0}$  is bounded, compactly supported, and continuous away from a $\Leb-$zero measure set, with $\norm{r}_1=1$ and such that $\supp r$ contains a nontrivial interval. Define a measure $\nu$ on $\mathbb{R}$ via $d\nu(E)=r(E)dE.$ Let 

\begin{eqnarray*}
M&=&\sup\{|E|:E\in \supp(r)\}\\
I&=&[-M, M]\\
\Sigma_0&\defeq& [-M-2-\|V_0\|_\infty, M+2+\|V_0\|_\infty]\\
\Omega&=&I^{\mathbb{Z}}\\
d\mu(\omega)&=&\prod_{n\in \mathbb{Z}}r(\omega_n)d\omega_n.
\end{eqnarray*}
Above, we let $V_0=\{V_0(n)\}_{n\in\mathbb{Z}}\in\ell^\infty(\mathbb{Z})$, be some fixed bounded background potential.
For $n\in\mathbb{Z}$, and any fixed $\alpha\in \mathbb{Z}_+$, we define $V_{\omega}(\alpha n)=f_0\omega(n),\,V_\omega(\alpha n+1)=f_1\omega(n),\, V_\omega(\alpha n+2)=f_2\omega(n),\dots,V_\omega(\alpha n+\alpha-1)=f_{\alpha-1}\omega(n)$ for each $\omega\in \Omega,$ where $f_i>0\,$ for $i=0,1,\dots, \alpha-1$ are fixed real numbers. That is, the potential is of the form 
\begin{equation}\label{eq0001} V_\omega(n)=\sum_{k\in\mathbb{Z}}\omega(k)f(n-k\alpha),\end{equation}
 where the single site potential $f:\mathbb{Z}\to\mathbb{R}_{>0}$ is supported on $\{0,1,\dots,\alpha-1\},$ and we take $f(i)=f_i$ for $i=0,1,\dots, \alpha-1$. In other words, one $i.i.d$ random variable determines the value of the potential on a block of length $\alpha$. From now on, we will refer to this as an $\alpha-$block. 

 With this notation, we define a one parameter family of Schr\"odinger operators, $\{H_{\omega}\}_\omega$, on $\ell^2(\mathbb{Z})$ as follows

\begin{equation}\label{eq00}\left(H_{\omega}u\right)(n)=u(n+1)+u(n-1)+\left(V_\omega(n)+V_0(n)\right)u(n).\end{equation}

We wish to emphasize that we do not think of $H_\omega$ as a single operator, but rather as an operator valued function on the probability space $\Omega$ (i.e. $\omega\mapsto H_\omega$). As such, we are generally interested in statements about $H_\omega$ that hold almost surely; that is, with probability one on $\Omega$.  In the literature this model is known as the {\it generalized Anderson model} in the discrete setting, see for example \cite{dss}. 

In the special case when $\alpha=1$ and $f_0=1$, this model, initially introduced in 1958 by P. W. Anderson to explain various quantum mechanical effects of disordered media, is now known  in literature as the \textit{Anderson model}. The simplest non-trivial case, where the support of the distribution contains precisely two elements is known as the {\it Bernoulli-Anderson model}. It is well known that, the spectrum of the {\it Anderson model} has a simple description, namely we have $$\Sigma\defeq\sigma({H_\omega})=\big[-2,2\big]+ \supp( r)\defeq\big\{a+b:a\in\big[-2,2\big], b\in\supp(r)\big\},$$ for $\mu-$ almost every $\omega\in\Omega.$ That is, it is simply given by the union of the translates of the spectrum of the Laplacian by points in $\supp r$. Since, by assumption, $\supp r$ is compact, the above description shows that the spectrum of an Anderson model will always be given by a finite union of compact intervals. For a proof of this result, see \cite{gs}. Below, we also give a set description of the spectrum for the {\it generalized Anderson model}, though it is not as simple as for the original Anderson model. 

One interesting property to study for the Anderson model is the phenomenon of {\it localization}. There are typically two separate statements referring to localization: a spectral statement and a dynamical one. Given some interval $I\subset \mathbb{R}$, we say that the operators $H_\omega$ exhibit {\it spectral localization} in $I$ if, almost surely, they have pure point spectrum, with exponentially decaying eigenfunctions. Historically, the discovery that dense pure point spectrum can occur, came as a surprise to the mathematical community- this would be the case if $I$ is a nontrivial interval of $\sigma(H_\omega)$ above - as will be the situation in our case. On the other hand, different notions of {\it dynamical Anderson localization} have been used in literature. However, in essence, dynamical localization refers to an absence of transport in a random medium. This is typically quantified via (almost-sure) bounds on the moments of wave packets such as
$$\sup_{t}\sum_{n\in\mathbb{Z}}|n|^p\left|\langle{\delta_n,e^{-itH_\omega}\delta_0\rangle}\right|^2<\infty,$$ for all $p>0.$ In some instances, one can prove stronger statements, such as replacing the {\it almost sure} condition by an expectation $\mathbb{E}(\cdot),$ as is the case via the Kunz-Souillard approach to localization in dimension one. For a more elaborate discussion of this method in the case of Schr\"odinger operators, see \cite{cfks}.

  The first mathematically rigorous proof of strong dynamical localization for the actual Anderson model, for one dimensional discrete Schr\"odinger operators, was originally given by H. Kunz and B. Souillard in \cite{ks}. For a version of the continuum model, namely one with time-continuous randomness, Goldsheid, Molachov, and Pastur gave the first spectral localization proof.  Following Kunz-Souillard's work, there have been a few extensions of their method in different directions and settings. In \cite{bs}, Simon showed that it is possible to allow for the potential to decay at a specific rate and still obtain pure point spectrum. Though, in this situation one no longer necessarily gets exponential localization or exponentially decaying eigenfunctions of the operators at hand. Recently, in \cite{vb} we show that the original localization result of Kunz-Souillard and that of Simon in \cite{bs} also hold true for any Jacobi operators (i.e. the case where the non-diagonal entries are no longer 1, but rather any positive bounded sequence). 
In 2011, Damanik and Stolz in \cite{ds} developed a continuum analogue of the discrete Kunz-Souillard approach. It is their technique that we adopt and use extensively in this paper to develop another extension of the original Kunz-Souillard approach to localization for random Schr\"odinger operators in the discrete setting in one dimension. Recently, in \cite{dg}, Damanik and Gorodetski further generalize the Kunz-Souillard approach by allowing certain correlations among the random variables defining the potentials. In that paper, they also give some interesting applications to almost periodic Schr\"odinger operators.

In an appropriate formulation, it is known that dynamical localization implies spectral localization, while the converse is not true in general. For example, the so called random dimer model serves as a counterexample to this implication (see \cite{JS} and \cite {JSS} for a more elaborate description). One typically needs `` spectral localization +$\epsilon$" to imply dynamical localization in some suitable formulation. This relationship was studied by del Rio, Jitomirskaya, Last, and Simon in \cite{rjls}.

There are different approaches to localization: {\it Spectral averaging} can be used to study spectral localization; one can also study both spectral and dynamical localization via methods such as, {\it multi-scale analysis}, developed around 1983 by Fr\"ohlich and Spencer in \cite{fs}; {\it fractional moments method}, initially introduced by Aizenman and Molchanov in \cite{am};  and also, which is what we do in this paper, the {\it Kunz-Souillard} method. Each method has its advantages and disadvantages. While, for example, some results which have been proven using {\it multiscale analysis} are well out of reach of the {\it fractional moments method} or the {\it Kunz-Souillard} approach, it is extremely technical and thus harder to work with. On the other hand, though, for example, one can prove less results via the {\it fractional moments method}, first off, it is much more elementary in nature, and second, under appropriate stipulations, it can provide richer and stronger results than {\it multiscale analysis}. 

The basic idea behind the {\it fractional moments} method is very simple: one first tries to establish exponential decay for the fractional moments of the {\it Green's function} (i.e. the matrix elements of the resolvent of $H_\omega$) and conclude dynamical and spectral localization from these bounds. There are typically two fairly different approaches to establish dynamical localization from exponentially decay of the fractional moments of Green's function. The first one was initially developed by Graf in \cite{g}, and works directly in infinite volume; that is, it does not need to first restrict the Hamiltonian $H_\omega$ to a finite box. 
The second method begins by first considering restrictions of the Hamiltonian to some finite box and it relies on the so-called {\it finite volume eigenfunction correlators,} which arise from the eigenfunction expansion of the Hamiltonian. The second approach is actually quite similar in flavor to the Kunz-Souillard method, as far as the general ideas go. It has also proven to be very robust under generalizations. For example, it can be extended to prove localization for the continuum Anderson model, which actually took almost a decade for it to be achieved, but nonetheless, was finally settled in 2006 in \cite{aensg} and \cite{bnsg}. About four years later, in \cite{hss}, a fractional moment's method proof of localization at all energies was also provided for the one-dimensional continuum Anderson models, which finally filled a gap that existed in the literature up to that point.

 Similarly, while there are many limitations in the scope of the Kunz-Souillard approach to localization (maybe the biggest one being the lack of applicability to dimensions higher than one), its rather elementary nature and the rich results it produces make it very attractive and of considerable interest.

In the Kunz-Souillard approach one begins by first restricting the operator $H_\omega$ to some finite box, decomposing it in terms of its eigenspaces, and then via a change of variables rewriting the latter in terms of some integral operators. So, the real technical difficulty of the method lies in estimating the norms of these operators. Specifically, in the original work of Kunz and Souillard in \cite{ks}, the challenge was to show that the operator $T_1$ was a strict contraction. Though this was not achieved, it was  shown that the second iterate of $T_1$ is a strict contraction, which was enough to conclude dynamical localization. In this paper, as a consequence of our approach, in section $9$,  we actually manage to improve on this result, namely we deduce that the single operator $T_1$ is itself a strict contraction. We consider this to be an interesting result, since it is rarely the case that the generalization of an original result actually improves on it as well. Typically, whenever an extension of some previous work is done, the goal or the hope is that the extended work will continue to coincide with the original work, when restricted to a specific case, or something along those lines.


\section{Main results}

\subsection{Statements of the main results}

The main objective of this paper is to prove dynamical and spectral localization for the family of operators $\{H_\omega\}_{\omega\in\Omega}$ defined above. Specifically, we prove the following two theorems.

\begin{thm}\label{thm2.1} With $\Omega,\, \mu$, and $H_\omega$ as above, there exist constants $C,\gamma\in(0,\infty)$ such that 
$$\int_\Omega \left(\sup_{t\in\mathbb{R}}\left|\langle{\delta_m,e^{-itH_\omega}\delta_n\rangle}\right|\right)d\mu(\omega)\leq Ce^{-\gamma\left|\lfloor{\frac{m-n}{\alpha}\rfloor}\right|},$$ for all $m,n\in\mathbb{Z}$.
\end{thm}

For more pleasant exposition, we fix the following notation \begin{equation}a(m,n)=\int_\Omega \left(\sup_{t\in\mathbb{R}}\left|\langle{\delta_m,e^{-itH_\omega}\delta_n\rangle}\right|\right)d\mu(\omega).\end{equation}

\begin{thm}\label{thm2.2}
If there are constants $C,\gamma\in(0,\infty)$ such that $$\max_{n\in\{0,1\}} a(m,n)\leq Ce^{-\gamma|m|},$$ then for $\mu-$ almost every $\omega\in\Omega$, $H_\omega$ has pure point spectrum with exponentially decaying eigenfunctions. More precisely, these eigenfunctions obey estimates of the form $$|u(m)|\leq C_{\omega, \epsilon, u}e^{-(\gamma-\epsilon)|m|},$$ for small enough $\epsilon\in(0,\gamma).$
\end{thm}

\begin{proof} The proof of Theorem $\ref{thm2.2}$ is proven in almost identical way as in the special case when $\alpha=1,$ so we direct the reader to \cite{cfks} or \cite{dd}.
\end{proof}

We wish to note that the conclusion of Theorem $\ref{thm2.1}$ implies the hypothesis of Theorem $\ref{thm2.2}$. Thus, in the case of the Kunz-Souillard approach to localization, as stated before, one derives spectral localization from dynamical localization.

\begin{cor} Given any $V_0\in\ell^\infty(\mathbb{Z})$, then for $\mu-$ almost every $\omega\in\Omega$ the Schr\"odinger operator $H_\omega=\triangle+V_0+V_\omega,$ where $V_\omega$ is as in $(\ref{eq0001})$, has pure point spectrum. 
\end{cor}

This is an immediate consequence of the above Theorems. It is already known that any bounded potential can be perturbed by a random potential to obtain pure point spectrum, however, the above corollary establishes that this statement still remains true even when the perturbation is done with somewhat \emph{less random} potentials. 

In what follows, mainly for ease of notation, we will carry out all the computations for $V_0(n)\equiv 0$. For a discussion of the general case see the Appendix, where we will emphasize the key modifications one would have to do when running the below arguments with the fixed background potential present. 

\subsection{Discussion of the model}

Here we wish to give a brief heuristic discussion of our model, the importance of some of the specific assumptions we have made, which allow us to conclude {\it dynamical localization}, and some of the new challenges we faced.

First, it is of great importance, for our proofs to go through, that the {\it scaling} factors, $f_0,\dots, f_{\alpha-1}$, be all strictly positive. We use this assumption in a crucial way throughout some of the key proofs below, one of them being the proof that the change of variable is one-to-one.

The other crucial assumption is that the single-site distribution $\mu$ be purely absolutely continuous. For the discrete {\it generalized Anderson model}, in comparison to the continuum analogue in \cite{ds}, we actually  relax the condition on the density $r$, that is, we allow for a Lebesgue measure zero set of discontinuities. We remark, however, that following our arguments that justify this, one can relax this condition in the continuum case as well. 

We wish to point out that dynamical localization might no longer be present if one works with distributions supported on say a finite set. For example, if we start with the {\it Bernoulli-Anderson} model, say where $\supp \mu=\{\epsilon_a,\epsilon_b\}$, and each $\epsilon_a,\epsilon_b$ is assigned to each site with probability $q$ and $1-q$, and we continue to double up the sites, that is, we set $V_\omega(2n-1)=V_\omega(2n)=\omega_n$, for each $n\in\mathbb{Z}$, where $\omega=(\omega_n)_n\in\Omega\defeq\left(\supp\mu\right)^{\mathbb{Z}}$, then, in this case, one no longer obtains dynamical localization. This model, which in the literature now is known as the random-dimer model, was first introduced and studied in 1990 in \cite{dwp} by Dunlap, Wu, and Phillips. This is also the first example which demonstrates that, in general, {\it spectral localization} does not imply {\it dynamical localization}.

 Taking $q=1/2$, which also corresponds to the case of highest disorder, one may think of the random-dimer model as flipping a coin to determine whether the value of the potential at a pair of lattice sites should be $\epsilon_a$, or $\epsilon_b$. From an intuitive standpoint, the fact that one no longer gets localization in this scenario is a bit surprising, since it would seem reasonable to expect that this system is equivalent to the completely random system, with the only difference being, now having two sites per unit cell. 

In fact, it is still an interesting open problem, to see if one can extend the Kunz-Souillard approach to the case where the single-site distributions have a non-trivial singular part. 

One of the first challenges that we faced was finding the appropriate way to establish the injectivity of the change of variables map. Since we rely on Pr\"ufer variables, where the Pr\"ufer phase is initially only defined up to a multiple of $2\pi$, one needs to find the appropriate way to make such angles unique, while at the same time being able to establish injectivity of the change of variables. For this reason, the way one {\it forces } uniqueness, has to be, in some appropriate sense, natural. However, since we work in the discrete setting, in contrast to the continuum situation, one cannot exactly simply resolve this issue via imposing some continuity condition on the Pr\"ufer phase, since continuity for us doesn't make sense to begin with. We resolve this issue and develop some necessary results, of the continuum analogues, for the discrete Pr\"ufer phases in Section $\ref{dpv}$. Along these lines, part of the new challenge was developing some of the analogous background results for the discrete Pr\"ufer amplitudes and phases, which is something we do as need arises throughout the paper. The rest of the challenge was essentially finding the appropriate ways to discretize and adopt the continuum techniques developed by Damanik and Stolz in \cite{ds}.

We wish to point out that the discrete {\it generalized Anderson model} falls into the framework of the discrete one-dimensional random word models developed by Damanik, Sims, and Stolz in \cite{dss}. In their notation, the fundamental set of words for this model is $\mathcal{W}=\{\lambda w: \lambda \in \supp r\},$ where $w=(f_0,f_1,\dots,f_{\alpha -1})\in \mathbb{R}^\alpha.$ In other words, the {\it generalized Anderson model} that we consider here, is a special case of this more general framework. In that paper, the authors study discrete Schr\"odinger operators in one dimension whose potentials are obtained by randomly concatenating blocks form an underlying set of words. They use multiscale analysis to prove spectral localization at all energies, and dynamical localization away from a finite set of exceptional energies. 

The authors make the observation, which is an artifact of their techniques, that the number of these exceptional energies seems to decrease as the set of words increases (i.e. the support of the measure from which these random words are drawn). They go further as to conjecture, though they don't explicitly state it as such, that {\it for suitably rich word spaces there are no exceptional energies.} 

We wish to note that, in the special case of the {\it generalized Anderson model}, which is the central object of study in our paper, where we suppose that the distribution, $\mu$, has Lebesgue-a.e. continuous and compactly supported density, we give a positive answer to their \emph{conjecture}, for this particular class of distributions. Moreover, not only do we prove dynamical localization at all energies for the {\it generalized Anderson model}, but, via the Kunz-Souillard method, we also get exponential bounds, and thus dynamical localization in a stronger sense than in \cite{dss}. We do however wish to mention the important fact that their results are for a  larger class of operators, and as such, they hold true with minimal assumptions on the underlying word space. That is, they only assume a nontriviality condition on the word space, which for the {\it generalized Anderson model}, translates to $\mathcal{W}$ containing at least two words. 

We also wish to  note that $$H_{T_0\omega}=U_\alpha H_\omega U_\alpha^{-1}$$ where $T_0:\mathcal{W}^\mathbb{Z}\to \mathcal{W}^{\mathbb{Z}}$, is given by $(T_0\omega)_n=\omega_{n+1}$, and where $U_\alpha$ is the shift by $\alpha$  in $\ell^2(\mathbb{Z}).$

Finally, in the last section we prove uniform positivity of the Lyapunov exponents, where we even allow the distribution $\mu$ to be finitely supported. In fact, we provide a lower bound on $\#(\supp \mu)$, for which we have uniform positivity of the Lyapunov exponents for the {\it generalized Anderson model}. 

In \cite{hss}, Stolz. et. al, among other things, prove F\"ursteberg at all energies, for the continuum Anderson model in one dimension, under the assumption that $\supp \mu$ is not discrete. Their proof is based on some of the results and techniques developed in \cite{dss2} and a slightly more general result than the one proved in \cite{ks}. We believe that once one has proved a discrete analogue of the result in \cite{ks}, then the same technique could be adopted to proving positivity of the Lyapunov exponent at all energies even for the discrete {\it generalized Anderson model}. We however, do not take this route, but instead use a more elementary and direct approach which in turn allows us to relax the assumptions on the support of the distribution $\mu$; that is, we allow discrete supports for our distributions. We point out that, by itself, this could count as a stronger result than the result one could potentially obtain by the aforementioned approach, however, since for the {\it generalized Anderson model}, which we study in this paper, we begin by assuming that $\mu$ is absolutely continuous, the outcome would have been the same. 

Moreover, we demonstrate that the dependence on the size of the block for the lower bound on the cardinality of the support of $\mu$ is necessary. Specifically, we show that for any given $\alpha$, there exists some finitely supported distribution $\mu$, for which the Lyapunov exponent for the {\it generalized Anderson model} vanishes for at least one energy. Thus, it is not possible to find an $\alpha-$ independent lower bound for $\#(\supp \mu)$ for which one would get uniformly positive Lyapunov exponents for the {\it generalized Anderson model}. 


\subsection{Almost sure spectrum of $H_\omega$}
In this section we give a description of the almost sure spectrum, $\Sigma$, of $H_\omega$ as a set. Specifically, we prove the following theorem

\begin{thm}\label{thm2.3-2.4} Let $\Sigma$ denote the almost sure spectrum of $H_\omega$. Then, we have $$\Sigma=\overline{\bigcup_{\substack{\omega=\{\omega_n\}\in\supp\mu\\ \omega=\text{periodic}}}\sigma(H_\omega)}.$$ 
\end{thm}

The proof follows the same general guidelines as in the original Anderson model (eg. see \cite{wk}). We begin by proving the following technical lemma,

\begin{lemm} \label{lemmm02.5}There is a full measure set $\Omega_0$ such that for any $\omega\in\Omega_0$, any finite set $I\subset\mathbb{Z}$, any periodic $\bar \omega=\{\bar\omega_n\}\in\supp\mu$, and any $\epsilon>0$, there exists a sequence $\{j_n\}$ in $\mathbb{Z}$ with $|j_n|\to \infty$ such that $$\sup_{i\in I}|V_\omega(i+j_n)-V_{\bar \omega}(i)|<\epsilon.$$ 

\end{lemm}

\begin{proof} Let $I$, $\bar\omega=\{\bar\omega_n\}$ and $\epsilon>0$, as in the statement, be fixed, and set $$A\defeq\{\omega\in\Omega:\,\sup_{i\in I}|V_\omega(i)-V_{\bar\omega}(i)|<\epsilon\}.$$  Using independence of the random variables and the definition of the topological support of $\nu$ we first show that $\mu(A)>0.$ Let $\Lambda_k\defeq[\alpha k,\alpha k+\alpha-1]$. Then, for some $r$ we have 
$$I\subset K_r\defeq \bigcup_{n=1}^{r}\Lambda_{k_n},$$ where $k_{n_i}\neq k_{n_j}$ for $i\neq j.$
With no loss of generality we can assume that $I\cap \Lambda_{k_n}\neq\emptyset$ for all $n=1,2,\dots, J$, else we can just disregard the $\Lambda_k-$s that result in empty intersections with $I$.  
 Observe that
\begin{align*}
A&=\{\omega\in\Omega:\, \sup_{i\in I}|V_{\omega}(i)-V_{\bar\omega}(i)|<\epsilon\}\\
&=\bigcap_{n=1}^r\{\omega\in\Omega:\, b_n |\omega_{k_n}-\bar\omega_{p_n}|<\epsilon\, \}\\
&=\bigcap_{n=1}^r\{\omega\in\Omega:\, \omega_{k_n}\in(\bar\omega_{p_n}-\epsilon\, b_n^{-1}, \bar\omega_{p_n}+\epsilon\, b_n^{-1})\},\\
\end{align*} 
where $b_n=\max_{j\in I\cap \Lambda_{k_n}}\{f_j\}.$
So, using the fact that $\{\omega_k\}$ are random $i.i.d$ variables and the fact that $\bar\omega_{p_n}\in\supp \nu,$ we get
\begin{align*}
\mu(A)&=\mu\left(\bigcap_{n=1}^r\{\omega\in\Omega:\, \omega_{k_n}\in(\bar\omega_{p_n}-\epsilon\, b_n^{-1}, \bar\omega_{p_n}+\epsilon\, b_n^{-1})\}\right)\\
&=\prod_{n=1}^r\nu (\bar\omega_{p_n}-\epsilon\, b_n^{-1}, \bar\omega_{p_n}+\epsilon\, b_n^{-1})\\
&>0.
\end{align*}
 
\noindent
Next, pick a sequence $\{l_n\}\in\mathbb{Z}$ such that for any $n\neq m$ we have $$|l_n-l_m|>2\alpha\diam K_r\,,$$ where $\alpha$ is the block size, and such that $V_\omega(i)$ and $V_\omega(i+l_n)$ have the same coefficient $f_j$ that multiplies the corresponding random $i.i.d$ at the respective sites. Then, the events $$A_n=A_n(I, q, \epsilon) \defeq\{\omega\in\Omega:\, \sup_{i\in I}|V_\omega(i+l_n)-V_{\bar\omega}(i)|<\epsilon\}$$ are independent, and it is easy to see that $\mu(A_n)=\mu(A)>0.$ So, since $\{A_n\}$ are independent and $\displaystyle \sum_n\mu(A_n)=\infty$, from Borel-Canelli lemma we know that $$\Omega_{I, q, \epsilon}\defeq\{\omega\in\Omega:\, \omega\in A_n \text{ for infinitely many } n\}$$ has probability one; that is, $\mu\left(\Omega_{I,q,\epsilon}\right)=1.$ We wish too point out that the desired sequence $\{j_n\}$ will be a subsequence of $\{l_n\}$.  Let $ C_0$ be a dense countable subset of $\supp \nu$. Then, since the collection of all finite subsets $I$ of $\mathbb{Z}$ is countable, we get that the set 

\begin{equation}\label{eq:04}\Omega_0\defeq \bigcap_{I\subset \mathbb{Z}, q\in C_0, n\in\mathbb{Z}_+}\Omega_{I, q, \frac{1}{n}}\end{equation}
\noindent
as a countable intersection of sets of measure one, has measure one; that is, $\mu\left(\Omega_0\right)=1.$
It is easy to see that, by definition, $\Omega_0$ satisfies the requirements of the claim, and thus concluding the proof! 
\end{proof}

\begin{proof}[Proof of Theorem $\ref{thm2.3-2.4}$] The inclusion $``\subset"$ is a standard result; it essentially follows from strong resolvent convergence, see for example \cite[p.290, Thm.VIII.24]{rs'}.
To prove the reverse inclusion $`` \supset"$, we will use the \emph{Weyl's criterion}, $$ \lambda\in\sigma(H_\omega) \text{ iff }\, \exists \, \psi_n\in\ell^2(\mathbb{Z}),\,\text{ with }\, \|\psi_n\|=1 \text{ such that } \|(H_\omega-\lambda)\psi_n\|\to 0 \text{ as } n\to \infty.$$ 
\noindent
Let $\displaystyle\lambda\in \bigcup_{\substack{\omega=\{\omega_n\}\in\supp\mu\\ \omega=periodic}}\sigma(H_\omega)$. Then $\lambda\in\sigma(H_{\bar\omega})$ for some periodic $\bar\omega\in\supp\mu.$ 
So, by the {\it Weyl's criterion}, there exist some finitely supported sequence $\varphi_n\in\ell^2(\mathbb{Z})$ with $\|\varphi_n\|=1$ and $\|(H_{\bar\omega}-\lambda)\varphi_n\|\to 0$ as $n\to\infty.$
Set $\displaystyle\varphi^{(j)}(i)=\varphi(i-j)$. 
Our goal is to construct a Weyl sequence, $\psi_n$, for the operator $H_\omega$, where $\omega\in\Omega_0.$ By Lemma $\ref{lemmm02.5}$, for every $\omega\in\Omega_0$, where $\Omega_0$ is defined in equation $(\ref{eq:04})$, there is a sequence $\{j_n\}$ with $|j_n|\to \infty$ such that

\begin{equation} \sup_{i\in\supp \varphi_n}|V_\omega(i+j_n)-V_{\bar\omega}(i)|<\frac{1}{n}.\end{equation}

We claim that $\psi_n\defeq \varphi_n^{(j_n)}$ is a Weyl sequence for $H_\omega$ and $\lambda.$ Let $H_{\omega, j_n}$ be the operator with potential $V_{\omega,j_n}(i)=V_\omega(i+j_n).$ It is a straightforward computation to show that $$\left(H_\omega\varphi_n^{(j_n)}\right)(i+j_n)=\left(H_{\omega, j_n}\varphi_n\right)(i).$$ So, as a result one also has \begin{equation}\|(H_\omega-\lambda)\varphi_n^{(j_n)}\|=\|(H_{\omega,j_n}-\lambda)\varphi_n\|.\end{equation} On the other hand,

\begin{align*}
\|(H_{\omega,j_n}-H_{\bar\omega})\varphi_n\|^2&=\|(V_{\omega,j_n}-V_{\bar\omega})\varphi_n\|^2\\
&=\sum_{i\in\mathbb{Z}}\left|\left(V_{\omega,j_n}(i)-V_{\bar\omega}(i)\right)\varphi_n(i)\right|^2\\
&=\sum_{i\in\mathbb{Z}}|V_{\omega}(i+j_n)-V_{\bar\omega}(i)|^2|\varphi_n(i)|^2\\
&\leq \sup_{i\in\supp\varphi_n}|V_{\omega}(i+j_n)-V_{\bar\omega}(i)|^2\|\varphi_n\|^2\\
&<\frac{1}{n^2}.
\end{align*}
Then,
\begin{align*}
\|(H_\omega-\lambda)\psi_n\|&=\|(H_\omega-\lambda)\varphi_n^{(j_n)}\|\\
&=|(H_{\omega,j_n}-\lambda)\varphi_n\|\\
&\leq\|(H_{\omega,j_n}-H_{\bar\omega})\varphi_n\|+\|(H_{\bar\omega}-\lambda)\varphi_n\|\\
&\leq \frac{1}{n}+\|(H_{\bar\omega}-\lambda)\varphi_n\|\\
&\to 0 \text{ as } n\to \infty.
\end{align*}
\noindent This shows that $\lambda\in\sigma(H_\omega)$ as desired. 
\end{proof}

\begin{rmk}
We wish to point out that in the special case where $f_i\equiv 1$ one gets that the almost sure spectrum for the {\it generalized Anderson model} is the same as for the original Anderson model; that is, $$\displaystyle \Sigma=[-2,2]+\supp\nu.$$
\end{rmk}

\subsection{The strategy to prove Theorem $\ref{thm2.1}$}\label{sec2.3.}

To prove Theorem $\ref{thm2.1}$ we begin by considering the restriction of our operator $H_\omega$ to some finite box. Specifically, for some fixed $L\in\mathbb{Z}_+$, denote by $H_\omega^L$ the restriction of $H_\omega$ to $\ell^2\left([-\alpha L,\alpha L-1]\cap\mathbb{Z}\right).$ That is,

\begin{equation*}
  H_{\omega}^L=\left(
  \begin{array}{cccccc} 
    V_{\omega}(-\alpha L)& 1& &&&  0 \\ 
    1& V_{\omega}(-\alpha L+1)&&& &  0\\ 
 0&1&&&&\\
    \vdots&  &&\ddots&1  & \vdots \\ 
    &&&&V_{\omega}(\alpha L-2)&1\\
    0& &&\ldots&1& V_{\omega}(\alpha L-1) \\ 
  \end{array}
  \right).
\end{equation*}

Let $\{E_l(\omega)\}_l$ and $\{v_\omega^{L,l}\}_l$, for $l=1,2,\dots, 2\alpha L$, be the eigenvalues and the corresponding normalized eigenfunctions of $H_\omega^L$, respectively. Define 

\begin{eqnarray*}
a_L(m,n)&=&\int_\Omega \left(\sup_{t\in\mathbb{R}}\left|\langle{\delta_m,e^{-itH_\omega^L}\delta_n\rangle}\right|\right)d\mu(\omega)\\
\rho_L(m,n)&=&\int_\Omega\left(\sum_l\left|\langle{\delta_m,v_\omega^{L,l}\rangle}\right|\left|\langle{\delta_n,v_\omega^{L,l}\rangle}\right|\right)d\mu(\omega),
\end{eqnarray*}
and notice that the above integrals are simply $2L-$fold iterated integrals. This follows from the fact that $H_\omega^L$ depends only on the entries $\omega_{-L},\dots, \omega_{L-1}$, and the fact that $\omega_n$'s are $i.i.d$, so the measure $d\mu(\omega)$ is simply the product measure.\\

In laying down the groundwork for the proof of Theorem $\ref{thm2.1}$ we begin by stating the following two lemmas. They are easy to prove, see for example \cite[pp.192-193]{cfks}, however, for completeness and reader's convenience, we include the brief arguments here.

\begin{lemm}\label{lemma02.3} For $m,n\in \mathbb{Z}$ we have $$a(m,n)\leq\liminf_{L\to\infty}a_L(m,n).$$
\end{lemm}

\begin{proof}
First, regarding $H_{\omega}^{L}$ as an operator in $\ell^2(\mathbb{Z})$, in the natural way, we observe that $H_{\omega}^{L}$ converges strongly to $H_{\omega}.$ As a consequence, $e^{-itH_{\omega}^{L}}$ converges strongly to $e^{-itH_{\omega}}$, for each $t\in \mathbb{R},$ and every $\omega.$ As a result, we also have
$$\lim_{L\to \infty}\left|\langle{\delta_m, e^{-itH_{\omega}^{L}}\delta_n}\rangle\right|=\left|\langle{\delta_m, e^{-itH_{\omega}}\delta_n}\rangle\right|.$$ Next, for each $t\in \mathbb{R}$, we have 
$$\left|\langle{\delta_m, e^{-itH_{\omega}^{L}}\delta_n}\rangle\right|\leq\sup_{t'\in\mathbb{R}}\left|\langle{\delta_m, e^{-it'H_{\omega}^{L}}\delta_n}\rangle\right|.$$ Taking lim inf of both sides we obtain:

$$\left|\langle{\delta_m, e^{-itH_{\omega}}\delta_n}\rangle\right|=\lim_{L\to \infty}\left|\langle{\delta_m, e^{-itH_{\omega}^{L}}\delta_n}\rangle\right|\leq\liminf_{L\to\infty}\sup_{t'\in\mathbb{R}}\left|\langle{\delta_m, e^{-it'H_{\omega}^{L}}\delta_n}\rangle\right|.$$ Hence, 
$$\sup_{t\in\mathbb{R}}\left|\langle{\delta_m, e^{-itH_{\omega}}\delta_n}\rangle\right|\leq\liminf_{L\to\infty}\sup_{t\in\mathbb{R}}\left|\langle{\delta_m, e^{-itH_{\omega}^{L}}\delta_n}\rangle\right|.$$
The result follows by an application of Fatou's lemma.
\end{proof}

\begin{lemm}\label{lemma02.4} For $L\in\mathbb{Z}_+$, and $m,n\in\mathbb{Z}$ we have $$a_L(m,n)\leq\rho_L(m,n).$$
\end{lemm}

\begin{proof} Using the eigenfunction expansion of $H_\omega^L$ we have

\begin{align*}
a_L(m,n)&=\int_{\Omega}\left(\sup_{t\in\mathbb{R}}\left|\langle{\delta_m, e^{-itH_{\omega}^{L}}\delta_n}\rangle\right|\right)d\mu(\omega)\\
&=\int_{\Omega}\left(\sup_{t\in\mathbb{R}}\left|\langle{\delta_m, e^{-itH_{\omega}^{L}}\sum_l\langle{\delta_n,v_{\omega}^{L,l}}\rangle}v_{\omega}^{L,l}\rangle\right|\right)d\mu(\omega)\\
&\leq\int_{\Omega}\left(\sup_{t\in\mathbb{R}}\sum_l\left|\langle{\delta_m, e^{-itE_{\omega}^{L,l}}\langle{\delta_n,v_{\omega}^{L,l}}\rangle}v_{\omega}^{L,l}\rangle\right|\right)d\mu(\omega)\\
&=\int_{\Omega}\left(\sup_{t\in\mathbb{R}}\sum_l \left|\langle{\delta_m,v_{\omega}^{L,l}}\rangle\right|\left|\langle{\delta_n,v_{\omega}^{L,l}}\rangle\right|\right)d\mu(\omega)\\
&=\rho_L(m,n).
\end{align*}
\end{proof}

So, the strategy of the proof is to show that $\rho_L(m,n)\leq Ce^{-\gamma\left|\lfloor{\frac{m-n}{\alpha}\rfloor}\right|}$ where $C$ is some $L-$independent constant. For simplicity, we will only consider the case for $n=0$, since the general case is proved analogously. The first step in establishing this bound will be to rewrite $\rho_L(m,0)$, via a change of variables, in a suitable way that will make it easier to estimate the norm of the resulting expression. This will be achieved in Lemma $\ref{lemma4.9}$. Before we begin describing the change of variables, we will take a short detour to introduce and develop some of the relevant results regarding discrete Pr\"ufer variables, which we will need in the later sections.


\section{Discrete Pr\"{u}fer variables}\label{dpv}

Let $u_{-L}(\circ, \omega, E)$ be the solution of \begin{equation}\label{eqd1}u(n+1)+u(n-1)+V_{\omega}(n)u(n)=Eu(n),\end{equation} with potential $V_{\omega}$ as above, satisfying the initial conditions $u(-\alpha L-1)=0,\, u(-\alpha L)=1$ where $\alpha$ is as above and $L$ is some fixed positive integer. Here the subscript $-L$ is there to indicate that we start solving the difference equation recurrently from left to right.

We define the corresponding Pr\"ufer phase $\phi_{-L}(\circ, \omega, E)$, and amplitude $R_{-L}(\circ, \omega, E)$ to be 
\begin{eqnarray}\label{eq2}u_{-L}(n,\omega, E)&=&R_{-L}(n, \omega, E)\sin \phi_{-L}(n, \omega, E)\\
u_{-L}(n+1, \omega, E)&=&R_{-L}(n,\omega, E)\cos \phi_{-L}(n, \omega, E),\nonumber\end{eqnarray}

\noindent normalized at $-\alpha L-1$ (i.e. $R_{-L}(-\alpha L-1)=1$).

Since from the boundary condition at $-\alpha L-1$ it follows that the \emph{initial} angle $\phi_{-L}(-\alpha L-1,\omega, E)$ is fixed to be $0$ (in fact, apriori it is only fixed up to $\mod \pi$, but we choose it to be precisely zero) we make the Pr\"{u}fer angle $\phi_{-L}(\cdot,\omega, E)$ unique in a way which was motivated by a discussion in \cite{rj}.


First, observe that we have:

\begin{equation}
\left(
\begin{array}{c}
u(\alpha k+\alpha)\\
u(\alpha k+\alpha-1)
\end{array}
\right)=T_{E,\omega}(\alpha k+\alpha-1)\dots T_{E,\omega}(\alpha k)\left(\begin{array}{c}u(\alpha k)\\ u(\alpha k-1)\end{array}\right),
\end{equation}

where $T_E(m)$ is the one-step transfer matrix; that is,

\begin{equation*}
T_{E,\omega}(m)=\left(\begin{array}{cc}E-V_\omega(m)&-1\\1&0\end{array}\right)
\end{equation*}

Next, for each  $k=-L,\dots, L-1,$ we fix a homotopy \begin{equation}F_{E,\alpha}^i(k,s)=\left\{\begin{array}{lr}\left(\begin{array}{cc}\cos f(s)&-\sin f(s)\\ \sin f(s)&\cos f(s)\end{array}\right)& 0\leq s\leq 1/2\\ 
\left(\begin{array}{cc} g(s)\big(E-V_\omega(\alpha k+i)\big)&-1\\1&0\end{array}\right)&  1/2 \leq s\leq 1\\ \end{array}\right.
\end{equation}
where $i=0,\dots,\alpha-1$ and $f:[0,\frac{1}{2}]\to[0,\pi/2]$ and $g:[\frac{1}{2},1]\to [0,1]$ are $C^{\infty}$ functions such that 

\begin{enumerate}[(i)]
\item $f(s)=0$ for $s$ close to $0$ and $f(s)=\pi/2$ for $s$ near $\frac{1}{2}$;
\item $g$ is strictly increasing with $g(1/2)=0$, and $g(1)=1$. 
\end{enumerate}

So, we have $F_{E,\alpha}^i(k,0)=I_2$ and $F_{E,\alpha}^i(k,1)=T_{E,\omega}(\alpha k+i),$ where $I_2$ is the identity matrix. 
Next, set \begin{equation} F_{E,\alpha}(k,s)=\prod_{i=\alpha-1}^0 F_{E,\alpha}^i(k,s)\end{equation}

and observe that we also get $F_{E,\alpha}(k,0)=I_2$ and $F_{E,\alpha}(k,1)=T_{E,\omega}(\alpha k+\alpha-1)\cdot\dots\cdot T_{E,\omega}(\alpha k)$.\\

So, given that the phase $\phi_{-L}(\alpha(k-1)+\alpha-1)$ has been fixed, then the {\it next} phase $\phi_{-L}(\alpha k+\alpha -1)$ is uniquely determined via the homotopy $F_{E,\alpha}$ above (we consider only these phases for reasons which will become clear later). This is true since the angle now changes continuously from site to site via the above homotopy. \\


\begin{lemm}\label{lemma03.1} For each $-L\leq k\leq L-1$, we have
\begin{equation}\label{eq7'}
\left|\phi_{-L}(\alpha k+\alpha-1,\omega_{-L},\dots,\omega_k,E)-\phi_{-L}(\alpha(k-1)+\alpha -1,\omega_{-L},\dots,\omega_{k-1},E)\right|< \pi B(\alpha)<\infty,\end{equation} for some constant $B(\alpha)>0$, uniform in $k$, $E\in\Sigma_0$ and $\omega_j'$s.
\end{lemm}
\begin{proof} For each $E,\alpha, k$ let us define a function $G_{E,\alpha}^k: [0,1]\to\mathbb{R}^2$ by $$G_{E,\alpha}^k(s)=F_{E,\alpha}(k,s)\left(u_{-L}(\alpha k), u_{-L}(\alpha k-1)\right)^t,$$ where $F_{E,\alpha}(\cdot,\cdot)$ is as above. Since $G_{E,\alpha}^k(0)$ determines the angle $\phi_{-L}(\alpha k-1)$ and $G_{E,\alpha}^k(1)$ the angle $\phi_{-L}(\alpha k+\alpha-1)$, it suffices to prove that the number of times that the curve  $\Gamma\defeq G_{E,\alpha}^k([0,1])$ winds around the origin is bounded above by a finite number, and the bound is uniform in $k\in\{-L,\dots,L-1\}$, $E\in \Sigma_0$, and $\omega_k\in\supp\nu$. To this end, first we note that as $s$ ranges from $0$ to $1/2$, $F_{E,\alpha}(k,s)$ is simply $R_{f(s)}^\alpha$ where 

$$R_{f(s)}=\left(\begin{array}{cc}\cos f(s)&-\sin f(s)\\ \sin f(s)&\cos f(s)\end{array}\right)$$ is the rotation matrix. Since $f$ goes from $0$ to $\pi/2$ as $s$ ranges from $0$ to $1/2$, then the curve $\Gamma_0\defeq G_{E,\alpha}^k([0,1/2])$ would have traversed an angle of length at most $\alpha \left(\frac{\pi}{2}\right),$ for all $k=-L,\dots, L-1$ all $E\in \Sigma_0$ and all $\omega_k\in\supp\nu$. On the other hand, for $s\in[1/2,1]$, from the definition of $F_{E,\alpha}(k,s)$ and $G_{E,\alpha}^k$, we will have 
\begin{equation}
G_{E,\alpha}^k(s)=\left(\begin{array}{c}P_{\alpha, E,\omega_k}(g(s))\\ Q_{\alpha-1, E,\omega_k}(g(s))\end{array}\right)
\end{equation}
where $P_{m, E,\omega_k}$ and $ Q_{m, E,\omega_k}$ are polynomials of degree $m$. Similarly as before, we wish to show that the number of times $G_{E,\alpha}^k(s)$ winds around the origin as $s$ ranges from $1/2$ to $1$ is bounded above by a finite number, uniform in $k\in\{-L,\dots, L-1\}$, $E\in\Sigma_0$ and $\omega_k\in\supp\nu$. Next, since the curve $\Gamma_1\defeq G_{E,\alpha}^k([1/2,1])$ is smooth, to count the number of times it winds around the origin is equivalent to counting the number of $s\in[1/2,1]$ for which  $P_{\alpha, E,\omega_k}(g(s))=0.$ Since $P_{\alpha , E,\omega_k}$ is a polynomial of degree $\alpha$, then it has at most $\alpha$ real zeros, for all $E\in\Sigma_0$ and all $\omega_k\in\supp\nu$. In other words, since as $E$ ranges over $\Sigma_0$ and $\omega_k$ rangers over $\supp\nu$ the only thing that changes are coefficients of $P_{\alpha, E,\omega_k}(\cdot)$, then the number of real roots of $P_{\alpha, E,\omega_k}(\cdot)$ will be uniformly bounded above by $\alpha$.  So, computing the winding number of $\Gamma_1$ is now equivalent to counting the number of $s\in[1/2,1]$ for which $g(s)=x_0$, where $x_0$ is any of the real roots of $P_{\alpha, E,\omega_k}$.  Since, by hypothesis, $g$ is injective in $[1/2, 1]$ and since there are at most $\alpha$ real roots of $P_{\alpha, E,\omega_k}$, we conclude that there are at most $\alpha$ possible values of $s\in[1/2,1]$ for which $g(s)=x_0$, for all $E\in\Sigma_0$ and all $\omega_k\in\supp\nu.$ Finally, since $\Gamma=\Gamma_0\cup\Gamma_1$, putting together the two arguments above, we conclude that there is some integer $B(\alpha)$, uniform in $E$ and $\omega_k$, so that the curve $\Gamma$ traverses an angle of length at most $B(\alpha)$.  Moreover, from the definition of the homotopy $F_{E,\alpha}(k,s)$ it follows immediately that such a bound $B(\alpha)$ is independent of $k$; that is, it is the same on each $\alpha-$block. This concludes the proof!
\end{proof}



Rewriting expression $(\ref{eqd1})$, we get
\begin{equation}\label{eq3}\frac{u(n+1)}{u(n)}+\frac{u(n-1)}{u(n)}=E-V_{\omega}(n).\end{equation}
Using the expressions in $(\ref{eq2})$ the expression in $(\ref{eq3})$ becomes: 
\begin{equation}\label{eq4} \cot\phi_{-L}(n,\omega, E)+\tan \phi_{-L}(n-1, \omega, E)=E-V_\omega(n).\end{equation}
Similarly, we also get \begin{equation}\label{eq4'}R_{-L}(n,\omega, E)\sin\phi_{-L}(n,\omega, E)=R_{-L}(n-1,\omega,E)\cos\phi_{-L}(n-1,\omega,E).\end{equation}
We rewrite $(\ref{eq4})$ using the explicit definition of the potential $V_\omega$:
\begin{equation}\label{eq7}
\cot\phi_{-L}(n,\omega, E)+\tan \phi_{-L}(n-1, \omega, E)=E-f_i\omega_k,\end{equation}
where $k=\lfloor{\frac{n}{\alpha}\rfloor},$ and $i\in\{0,1,\dots, \alpha-1\}.$ 

\begin{lemm}\label{lemm1} For $j<\lfloor{\frac{n}{\alpha}\rfloor}$, we have
$$R_{-L}^2(n,\omega,E)\frac{\partial}{\partial \omega_j}\phi_{-L}(n,\omega,E)=\sum_{i=0}^{\alpha-1}f_iu_{-L}^2(\alpha j+i,\omega, E).$$
\end{lemm}

\begin{proof}
From equation $(\ref{eq7})$ we get: 

\begin{equation}\label{eq8}
\cot\phi_{-L}(n,\omega, E)=-\tan \phi_{-L}(n-1, \omega, E)+E-f_i\omega_k.
\end{equation}
Differentiating $(\ref{eq8})$ with respect to $\omega_j$ and since by assumption $j<k\defeq \lfloor{\frac{n}{\alpha}\rfloor}$, we get

\begin{equation}\nonumber
-\frac{1}{\sin^2\phi_{-L}(n,\omega,E)}\frac{\partial}{\partial \omega_j}\phi_{-L}(n,\omega,E)=-\frac{1}{\cos^2\phi_{-L}(n-1,\omega,E)}\frac{\partial}{\partial\omega_j}\phi_{-L}(n-1,\omega,E).
\end{equation}
Equivalently,
\begin{equation}\label{eq9}
\frac{\partial}{\partial \omega_j}\phi_{-L}(n,\omega,E)=\frac{\sin^2\phi_{-L}(n,\omega,E)}{\cos^2\phi_{-L}(n-1,\omega,E)}\frac{\partial}{\partial\omega_j}\phi_{-L}(n-1,\omega,E).
\end{equation}
Multiplying both sides of $(\ref{eq9})$ by $R^2_{-L}(n,\omega,E)$, and using the relation in $(\ref{eq4'})$, we get

\begin{equation}\label{eq10}
R^2_{-L}(n,\omega,E)\frac{\partial}{\partial \omega_j}\phi_{-L}(n,\omega,E)=R^2_{-L}(n-1,\omega,E)\frac{\partial}{\partial\omega_j}\phi_{-L}(n-1,\omega,E).
\end{equation}
Iterating $(\ref{eq10})$ and using the relation in $(\ref{eq2})$ we get

\begin{align}\label{eq11}\nonumber
R^2_{-L}(n,\omega,E)\frac{\partial}{\partial \omega_j}\phi_{-L}(n,\omega,E)&=R^2_{-L}(\alpha j+\alpha-1,\omega,E)\frac{\partial}{\partial\omega_j}\phi_{-L}(\alpha j+\alpha-1,\omega,E)\\ \nonumber
&=R^2_{-L}(\alpha j+\alpha-2,\omega,E)\frac{\partial}{\partial\omega_j}\phi_{-L}(\alpha j+\alpha-2,\omega,E)\\
&+f_{\alpha-1}R^2_{-L}(\alpha j+\alpha-1,\omega,E)\sin^2(\alpha j+\alpha-1,\omega,E)\\ \nonumber
&\vdots\\ 
&=\sum_{i=0}^{\alpha-1}f_iR_{-L}^2(\alpha j+i,\omega, E)\sin^2(\alpha j+i,\omega, E)\\
&=\sum_{i=0}^{\alpha-1}f_iu_{-L}^2(\alpha j+i,\omega, E),
\end{align}
as desired.
\end{proof}

Though obvious, for the record, we mention  special cases that arise if $j=\lfloor{\frac{n}{\alpha}\rfloor}$, namely, for $ 0\leq N\leq \alpha-1$

\begin{equation}\label{eq12}R_{-L}(\alpha j+N,\omega,E)\frac{\partial}{\partial \omega_j}\phi_{-L}(\alpha j+N,\omega,E)=\sum_{i=0}^{N}f_iu_{-L}^2(\alpha j+i,\omega, E).\end{equation}

One can also obtain a formula for the partial derivative of the phase with respect to the energy $E$. One can find a proof of this in \cite[Lem. 2]{JSS}. Since our situation is slightly different, we present here a proof of this result with the corresponding changes. 

\begin{lemm}\label{lemm2}  $$R^2_{-L}(n,\omega,E)\frac{\partial}{\partial E}\phi_{-L}(n,\omega,E)=-\sum_{j=n}^{-\alpha L}u_{-L}^2(j,\omega,E).$$
\end{lemm}

\begin{proof}
From equation $(\ref{eq7})$ we get: 

\begin{equation}\label{eq20}
\cot\phi_{-L}(n,\omega, E)=-\tan \phi_{-L}(n-1, \omega, E)+E-f_i\omega_k.
\end{equation}
Differentiating $(\ref{eq20})$ with respect to $E$ we get

\begin{equation}\nonumber
-\frac{1}{\sin^2\phi_{-L}(n,\omega,E)}\frac{\partial}{\partial E}\phi_{-L}(n,\omega,E)=-\frac{1}{\cos^2\phi_{-L}(n-1,\omega,E)}\frac{\partial}{\partial E}\phi_{-L}(n-1,\omega,E)+1
\end{equation}
Equivalently,
\begin{equation}\label{eq211}
\frac{\partial}{\partial E}\phi_{-L}(n,\omega,E)=\frac{\sin^2\phi_{-L}(n,\omega,E)}{\cos^2\phi_{-L}(n-1,\omega,E)}\frac{\partial}{\partial E}\phi_{-L}(n-1,\omega,E)-\sin^2\phi_{-L}(n,\omega, E)
\end{equation}
Multiplying both sides of $(\ref{eq211})$ by $R^2_{-L}(n,\omega,E)$, and using the relation in $(\ref{eq4'})$, and $(\ref{eq2})$ we get

\begin{equation}\label{eq222}
R^2_{-L}(n,\omega,E)\frac{\partial}{\partial E}\phi_{-L}(n,\omega,E)=R^2_{-L}(n-1,\omega,E)\frac{\partial}{\partial E}\phi_{-L}(n-1,\omega,E)-u^2_{-L}(n,\omega,E)
\end{equation}
Iterating $(\ref{eq222})$ and using the boundary condition at site $-\alpha L-1$, we get:

\begin{align}\nonumber
R^2_{-L}(n,\omega,E)\frac{\partial}{\partial E}\phi_{-L}(n,\omega,E)&=-\sum_{j=n}^{-\alpha L}u_{-L}^2(j,\omega,E).
\end{align}
This concludes the proof.
\end{proof}

\section{Change of variables}

\subsection{Introducing the change of variables}
We introduce the following change of variables 

\begin{eqnarray*}J_L:[-M,M]^{2L}\times \{0,1,\dots,2\alpha L-1\}&\to& \Sigma_0\times\mathbb{T}_\alpha^{2L-1}\times\{0,1,\dots,2B(\alpha)-1\}\\
(\bar \omega, l)&\mapsto& (E,\theta_{-L},\dots, \theta_{L-2}, N)\\
\end{eqnarray*} 

where $\mathbb{T}_\alpha\defeq \mathbb{R}/(2\pi B(\alpha) \mathbb{Z}),$  and $\Sigma_0=[-2-M,2+M],$ as follows: 

\begin{itemize}
\item for $k=-L,\dots,L-2$, we pick $\theta_k\in\mathbb{T}_\alpha$ such that $$\phi_{-L}(\alpha k+\alpha-1,\omega_{-L},\dots, \omega_k,E_l(\omega))\equiv \theta_k\mod 2\pi B(\alpha).$$ 
\item $E\in\Sigma_0$ is given by $$E=E_l(\omega).$$
\item $N\in\{0,1,\dots, 2B(\alpha)-1\}$ is defined such that $$l\equiv N\mod 2B(\alpha).$$
\end{itemize}

Since we are taking the phase angles modulo $2\pi B(\alpha)$, we need to make sure that this process does not produce any ambiguities. In other words, we need to ensure that the change of variables $J_L$ is one-to-one. For completeness, we prove this in the lemma below, by following the argument for the continuum case almost verbatim. For a proof of the original statement in the continuum setting, see \cite[Lemm. 2.3]{ds}. 
\begin{lemm}\label{chov}
The change of variables $J_L$ is one-to-one.
\end{lemm}

\begin{proof}
We first note that, for the $l^{th}$ eigenfunction, the Pr\"ufer phase $\phi_{-L}$ runs from $0$ at $-\alpha L-1$ to $\displaystyle \frac{\pi}{2}+N\pi$, for some $N\in \mathbb{Z}$, at $\alpha L-1$. More precisely, since the Pr\"ufer angle $\phi_{-L}$ is strictly monotone decreasing in $E$, it follows that for the $l^{th}$ eigenvalue, $E_l$ of $H_\omega^L$, we actually have $\phi_{-L}(\alpha L-1, \omega_{-L},\dots,\omega_{L-1}, E_l(\omega))=\pi/2-l\pi,$ where $l=0,1,\dots, 2\alpha L.$ Suppose that $(E,\theta_{-L},\dots, \theta_{L-2}, N)$ belongs to the range of $J_L$ and that there is 
some $(\omega_{-L},\dots, \omega_{L-1},l)\in [-M,M]^{2L}\times\{0,1,\dots, 2\alpha L\}$, such that

$$J_L(\omega_{-L},\dots, \omega_{L-1}, l)=(E,\theta_{-L}\dots, \theta_{L-2},N).$$ By the definition of $J_L$, for $k=-L,\dots, L-2$, we have 
\begin{equation}\label{eq019}
\phi_{-L}(\alpha k+\alpha-1,\omega_{-L},\dots, \omega_k,E_l(\omega))\equiv \theta_k\mod 2\pi B(\alpha).
\end{equation}
We will show, iteratively in $k$, that this in fact determines $\omega_k$ uniquely. We begin by fixing $\theta_{-L-1}\defeq 0$ and $\theta_{L-1}\defeq \pi/2+N\pi.$ In $(\ref{eq7'})$ we showed that  
\begin{equation}\label{eq020}
\left|\phi_{-L}(\alpha k+\alpha-1,\omega_{-L},\dots, \omega_k, E)-\phi_{-L}(\alpha(k-1)+\alpha-1,\omega_{-L},\dots,\omega_{k-1},E)\right|<B(\alpha)\pi.
\end{equation}
 Assume to the contrary that there is some other $\omega_k'\neq\omega_k$ such that $(\ref{eq019})$ holds as well. Then, there exist $n_1, n_2, n_3\in \mathbb{Z}$ such that 
\begin{eqnarray*}
\phi_{-L}(\alpha k+\alpha-1,\omega_{-L},\dots, \omega_k,E_l(\omega))&=&\theta_k+n_1(2\pi B(\alpha))\\
\phi_{-L}(\alpha k+\alpha-1,\omega_{-L},\dots, \omega_k',E_l(\omega))&=&\theta_k+n_2(2\pi B(\alpha))\\
\phi_{-L}(\alpha (k-1)+\alpha-1,\omega_{-L},\dots, \omega_{k-1},E_l(\omega))&=&\theta_{k-1}+n_3(2\pi B(\alpha))
\end{eqnarray*}
From Lemma $\ref{lemm1}$ we see that $\phi_{-L}(\alpha k+\alpha-1,\omega_{-L},\dots, \omega_k, E)$ is strictly increasing in $\omega_k$, so with no loss of generality we may assume that $n_1>n_2.$
Then,
\begin{align*}
\left|\phi_{-L}(\alpha k+\alpha-1)-\phi_{-L}(\alpha(k-1)+\alpha-1)\right|&=\left|(\theta_k-\theta_{k-1})+2\pi B(\alpha)(n_1-n_3)\right|\\
&=\left|(\theta_k-\theta_{k-1})+2\pi B(\alpha)(n_1-n_2)+2\pi B(\alpha)(n_2-n_3)\right|\\
&\geq\left|2\pi B(\alpha)(n_1-n_2)\right|-\left|(\theta_k-\theta_{k-1})+2\pi B(\alpha)(n_2-n_3)\right|\\
&\geq2\pi B(\alpha)(n_1-n_2)-B(\alpha)\pi\\
&\geq B(\alpha) \pi,
\end{align*}
where the last inequality follows from the fact that $n_1-n_2\geq 1.$ But, this contradicts $(\ref{eq020})$, so there must be a unique $\omega_k\in[-M,M]$ that satisfies $(\ref{eq019})$. Having reconstructed the unique $\bar \omega=(\omega_{-L},\dots, \omega_{L-1})$ for which $J_L(\bar \omega, l)=(E,\theta_{-L},\dots, \theta_{L-2},N)$, then from the fact that $\phi_{-L}(\alpha L-1,\bar\omega, E_l(\omega))=\pi/2-l\pi=\theta_{L-1}\defeq \pi/2+N\pi$, it follows that $l$ is also determined uniquely, thus concluding the proof.
\end{proof}
Our next goal is to show that, locally, any two given {\it consecutive} Pr\"ufer angles, $\theta_{k-1},\theta_k$ determine a unique {\it coupling constant} $\omega_k$. We do this, via adopting the analogous argument in \cite{ds}.\\

For $E,\lambda, \theta_{k-1},\theta_k\in\mathbb{R}$ let $u_{-}(\cdot, \theta_{k-1},\lambda, E)$ be the unique solution of the difference equation \begin{equation}\label{eq024} u(n+1)+u(n-1)+\lambda f(n-\alpha k) u(n)=Eu(n)\end{equation} with the initial conditions $u(\alpha k-1)=\sin \theta_{k-1}$ and $u(\alpha k)=\cos\theta_{k-1}$, where the subscript $"-"$ means that we are solving starting from left to right. Also, let $u_{+}(\cdot, \theta_k,\lambda, E)$, be the unique solution of $(\ref{eq024})$ with initial conditions $u(\alpha(k+1))=\cos \theta_k$ and $u(\alpha k+\alpha -1)=\sin \theta_k$, where the subscript $"+"$ means that we solve recursively from right to left. We wish to point out that $\lambda f(n-\alpha k)$ above is precisely the restriction of the potential $V_\omega$ to the $\alpha-$block $[\alpha k, \alpha k+\alpha-1].$ Let $\phi_{-}(\cdot, \theta_{k-1},\lambda, E)$, $R_{-}(\cdot, \theta_{k-1}, \lambda, E)$ and $\phi_{+}(\cdot, \theta_k,\lambda, E), \, R_{+}(\cdot, \theta_k,\lambda, E)$ be the Pr\"{u}fer phase and amplitude for $u_{-}(\cdot, \theta_{k-1}, \lambda, E)$ and $u_{+}(\cdot, \theta_k,\lambda, E)$, respectively, normalized at $\alpha k-1$ and $\alpha k+\alpha -1$ respectively; that is $\phi_{-}(\alpha k-1,\theta_{k-1},\lambda, E)=\theta_{k-1}$, $R_{-}(\alpha k-1,\theta_{k-1},\lambda, E)=1$, and $\phi_{+}(\alpha k+\alpha-1,\theta_k,\lambda, E)=\theta_k$, $R_{+}(\alpha k+\alpha-1,\theta_k,\lambda, E)=1.$ We also make the Pr\"ufer phases $\phi_{-}(\cdot, \theta_{k-1},\lambda, E)$ and $\phi_{+}(\cdot, \theta_k,\lambda, E)$  unique via the same argument as before. Moreover, similar relations as in $(\ref{eq4})$ and $(\ref{eq4'})$ hold in this set up as well. 

\begin{lemm}\label{lemm3.2'}
With the same notation as above we have
$$\frac{\partial}{\partial \theta_{k-1}}\phi_{-}(\alpha k+\alpha-1,\theta_{k-1},\omega_k,E)=\frac{1}{R_{-}^2(\alpha k+\alpha-1,\theta_{k-1},\omega_k,E)}.$$
\end{lemm}

\begin{proof}
Differentiating $$\cot\phi_{-}(\alpha k+\alpha-1,\theta_{k-1},\omega_k,E)+\tan \phi_{-}(\alpha k+\alpha-2,\theta_{k-1},\omega_k,E)=E-V_\omega(\alpha k+\alpha-1)$$ with respect to the initial angle $\theta_{k-1}$ we get \begin{equation}\label{eq025} \frac{\partial}{\partial \theta_{k-1}}\phi_{-}(\alpha k+\alpha-1,\theta_{k-1},\omega_k,E)=\frac{\sin^2\phi_{-}(\alpha k+\alpha-1,\theta_{k-1},\omega_k,E)}{\cos^2\phi_{-}(\alpha k+\alpha-2,\theta_{k-1},\omega_k,E)}\frac{\partial}{\partial \theta_{k-1}}\phi_{-}(\alpha k+\alpha-2,\theta_{k-1},\omega_k,E).\end{equation} Then, the result follows by iterating $(\ref{eq025})$, using the analogue of expression $(\ref{eq4'})$  for $R_{-}$, and the fact that $R_{-}(\alpha k-1)=1.$
\end{proof}
Next, we show that when $\phi_{-}(\alpha k+\alpha-1, \cdot,\lambda, E)$ is considered as a map from $\mathbb{T}_{\alpha}$ to $\mathbb{T}_\alpha$  it is well-defined. Let $ [\theta]\in \mathbb{T}_\alpha$, and let $x,y$ be two different representatives from this class; that is, $x=\theta+n_12\pi B(\alpha)$, and $y=\theta+n_22\pi B(\alpha)$ for some $n_1, n_2\in \mathbb{Z}.$ Then, since in general, by linearity, we have $u_{-}(\cdot, \theta_{k-1}+\pi, \lambda, E)=-u_{-}(\cdot, \theta_{k-1},\lambda, E)$, it follows that $\phi_{-}(\cdot, \theta_{k-1}+\pi,\lambda, E)=\phi_{-}(\cdot, \theta_{k-1},\lambda, E)+(2l(\cdot,\theta_{k-1})-1)\pi$, for some $l(\cdot,\theta_{k-1})\in \mathbb{Z}$. 
Since, apriori, $l$ depends on the space variable and the initial angle $\theta_{k-1}$, it may not be constant, however, we argue that this cannot be the case. By the way we have extended the Pr\"ufer phase continuously via the homotopy $F_{E,\alpha}$ it is not difficult to convince yourself that $\phi_{-}$ depends continuously on the space variable and the initial angle $\theta_{k-1}$. In particular, as a result, it is possible to perturb either the space variable or the initial angle $\theta_{k-1}$ sufficiently small such that the change in the difference $\phi_{-}(\cdot, \theta_{k-1}+\pi, \lambda, E)-\phi_{-}(\cdot, \theta_{k-1},\lambda, E)$ is strictly less than $\pi$, which would imply that the change in $(2l(\theta_{k-1})-1)\pi$ must be strictly less than $\pi$ as well. But, since $l(\cdot, \theta_{k-1})$ is an integer, this is impossible, unless it is constant. In particular, $l(x,\theta_{k-1})$ is constant for all $x$. But, since from the initial condition we have $\phi_{-}(\alpha k-1,\theta_{k-1},E)=\theta_{k-1}$, it follows that $l(x,\theta_{k-1})\equiv 1.$ So, we have \begin{equation}\label{eqeq1} \phi_-(\cdot,\theta_{k-1}+\pi,\lambda, E)=\phi_-(\cdot, \theta_{k-1},\lambda, E)+\pi.\end{equation}
Then,
\begin{align*}
\phi_{-}(\alpha k+\alpha-1,x,\lambda, E)\mod 2\pi B(\alpha)&=\phi_{-}(\alpha k+\alpha-1, \theta+n_12\pi B(\alpha),\lambda, E)\mod 2\pi B(\alpha)\\
&=\phi_{-}(\alpha k+\alpha-1, \theta, \lambda, E)+n_12\pi B(\alpha)\mod 2\pi B(\alpha)\\
&=\phi_{-}(\alpha k+\alpha-1, \theta, \lambda, E)\mod 2\pi B(\alpha)\\
&=\phi_{-}(\alpha k+\alpha-1, \theta, \lambda, E)+n_22\pi B(\alpha)\mod 2\pi B(\alpha)\\
&=\phi_{-}(\alpha k+\alpha-1, \theta+n_22\pi B(\alpha),\lambda, E)\mod 2\pi B(\alpha)\\
&=\phi_{-}(\alpha k+\alpha-1,y,\lambda, E)\mod 2\pi B(\alpha),
\end{align*}
as claimed. Similarly, one shows that $\phi_{+}(\alpha k-1,\cdot,\lambda, E)$ is well-defined, when considered as a map from $\mathbb{T}_\alpha$ to $\mathbb{T}_\alpha.$

Given any $x,y\in \mathbb{T}_\alpha$, if there is some coupling constant $\lambda \in [-M, M]$ such that $\phi_{-}(\alpha k+\alpha-1, y,\lambda, E)=x$ (or $\phi_{+}(\alpha k-1,x,\lambda, E)=y$) we set $\lambda(y,x,E)=\lambda.$ Next, we show that if such a coupling constant exists, then it must be unique. Given $x,y\in \mathbb{T}_\alpha$ suppose there are $\lambda_1,\lambda_2\in[-M,M]$ such that $\phi_{-}(\alpha k+\alpha-1,x,\lambda_1,E)=y=\phi_{-}(\alpha k+\alpha-1,x,\lambda_2,E)$; that is, there are some $n_1,n_2\in\mathbb{Z}$ such that $\phi_{-}(\alpha k+\alpha-1,x,\lambda_1,E)=y+n_12\pi B(\alpha)$ and $\phi_{-}(\alpha k+\alpha-1,x,\lambda_2,E)=y+n_22\pi B(\alpha).$ With no loss of generality, suppose that $\lambda_1<\lambda_2. $ But, since $f_i>0$ for all $i=0,1,\dots,\alpha-1$, from Lemma $\ref{lemm1}$ it follows that $\phi_{-}$ is strictly increasing when viewed as a function of the coupling constant $\lambda$, thus, it follows that $n_2>n_1.$  Let $\phi_{-}(\alpha k-1,x,\lambda_1,E)=x+n_42\pi B(\alpha)$ and $\phi_{-}(\alpha k-1,x,\lambda_2, E)=x+n_32\pi B(\alpha)$, so, for the same reason as before $n_3>n_4$. Then

\begin{align*}\medmuskip=1mu
|\phi_{-}(\alpha k+\alpha-1,x,\lambda_2,E)-\phi_{-}(\alpha k-1,x,\lambda_2, E)|=&|(y-x)+2\pi B(\alpha)(n_2-n_3)|\\
=&|(y-x)+2\pi B(\alpha)(n_2-n_1)+2\pi B(\alpha)(n_1-n_3)|\\
\geq &|2\pi B(\alpha)(n_2-n_1)|-|(y-x)+2\pi B(\alpha)(n_1-n_3)|\\
=&2\pi B(\alpha)(n_2-n_1)\\
&-|(y-x)+2\pi B(\alpha)(n_1-n_4)+2\pi B(\alpha)(n_4-n_3)|\\
\geq& 2\pi B(\alpha)(n_2-n_1)-|(y-x)+2\pi B(\alpha)(n_1-n_4)|\\
&-|2\pi B(\alpha)(n_4-n_3)|\\
=&2\pi B(\alpha)(n_2-n_1+n_3-n_4)\\
&-|\phi_{-}(\alpha k+\alpha-1,x,\lambda_1,E)-\phi_{-}(\alpha k-1, x,\lambda_1,E)|\\
\geq& 2\pi B(\alpha)(n_2-n_1+n_3-n_4)-\pi B(\alpha)\\
>& \pi B(\alpha),
\end{align*}
contradicting $(\ref{eq7'})$, and thus showing that there exists at most one such coupling constant $\lambda$. In other words, this shows that for each $E$, the function $\lambda_E: \mathbb{T}_\alpha^2\to [-M,M]$ given by $\lambda_E(x,y)\defeq \lambda(x,y,E)$ is well defined. Below, we prove an important property of this function which we will use later.

\begin{prop}\label{Leb}
Let $A\subset [-M,M]$. If $Leb(A)=0$, then $Leb\left(\lambda_E^{-1}(A)\right)=0$, where $Leb$ denotes Lebesgue measure. 
\end{prop}

\begin{proof}
It suffices to show that $\lambda_E\in C^1(\mathcal{D})$ and that $Leb\big(\{(x,y)\in\mathcal{D}: \nabla \lambda_E(x,y)=0\}\big)=0$ where $\mathcal{D}\defeq\{(x,y)\in\mathbb{T}_\alpha^2: \lambda_E(x,y)\, \textit{exists}\}.$ A quick computation shows that 
\noindent
\begin{align}\label{Leb01}\nonumber
\frac{\partial}{\partial y}\lambda_E(x,y)&=\frac{1}{\frac{\partial x}{\partial \lambda}}\\\nonumber
&=\frac{1}{\frac{\partial \phi_{+}(\alpha k-1,x,\lambda, E)}{\partial \lambda}}\\
&=-\frac{R_{+}^{2}(\alpha k -1,x,\lambda,E)}{\sum_{i=0}^{\alpha-1}f_iu_{+}^2(\alpha k+i,x,\lambda,E)},
\end{align}
and
\begin{align}\label{Leb02}\nonumber
\frac{\partial}{\partial x}\lambda_E(x,y)&=\frac{1}{\frac{\partial x}{\partial \lambda}}\\ \nonumber
&=\frac{1}{\frac{\partial \phi_{-}(\alpha k+\alpha-1,y,\lambda, E)}{\partial \lambda}}\\
&=\frac{R_{-}^{2}(\alpha k+\alpha -1,y,\lambda,E)}{\sum_{i=0}^{\alpha-1}f_iu_{-}^2(\alpha k+i,y,\lambda,E)}.
\end{align}
We claim that $\{(x,y)\in \mathcal{D}: \nabla \lambda_E(x,y)=0\}=\emptyset.$ If not, then this would imply that $R_-^2(\alpha k+\alpha-1,y,\lambda, E)=0$ and $R_+^2(\alpha k-1,x,\lambda,E)=0.$ Then, from the definition of $R_{\pm }$ we would have $u_-^2(\alpha k+\alpha-1,y,\lambda, E)+u_-^2(\alpha k+\alpha,y,\lambda, E)=0$ and $u_+^2(\alpha k-1,x,\lambda, E)+u_+^2(\alpha k,x,\lambda, E)=0.$ But, this in turn would imply that $u_{\pm}\equiv 0$, which is a contradiction. 

Finally, since $u_{\pm}(\alpha k+i,\cdot,\lambda, E)$ are polynomials in $\sin(\cdot), \cos (\cdot)$ and $R^2_{\pm}$ is a sum of two such functions, from $(\ref{Leb01})$ and $(\ref{Leb02})$, it follows that $\frac{\partial}{\partial x}\lambda_E$ and $\frac{\partial}{\partial y}\lambda_E$ are continuous in $y,x$, and thus $\lambda_E\in C^1(\mathcal{D})$. This concludes the proof!
\end{proof}

\begin{lemm}\label{lemmm3.2}
With the same notation as above, for $k=-L,\dots, L-2$, we have 

\begin{enumerate}[$(i)$]\item $R_{-L}(\alpha k+\alpha-1,\theta_{-L},\dots, \theta_k, E)=R_{-L}(\alpha k-1,\theta_{-L},\dots, \theta_{k-1}, E)R_{-}(\alpha k+\alpha -1,\theta_{k-1},\lambda(\theta_{k-1},\theta_k,E), E)$
\item $R_{-L}(\alpha k-1,\theta_{-L},\dots, \theta_{k-1}, E)=R_{-L}(\alpha k+\alpha-1,\theta_{-L},\dots, \theta_{k}, E)R_{+}(\alpha k-1,\theta_{k},\lambda(\theta_{k-1},\theta_k,E), E)$
\end{enumerate}
\end{lemm}
\begin{proof}
From $(\ref{eq4'})$ and the analogous expression for $R_{-}$ we get \begin{equation} R_{-L}(\alpha k+\alpha-1)=R_{-L}(\alpha k-1)\prod_{i=\alpha-1}^0\frac{\cos\phi_{-L}(\alpha k+i-1)}{\sin\phi_{-L}(\alpha k+i)}\end{equation} and \begin{equation}\label{eq030} R_{-}(\alpha k+\alpha-1)=R_{-}(\alpha k-1)\prod_{i=\alpha-1}^0\frac{\cos\phi_{-}(\alpha k+i-1)}{\sin\phi_{-}(\alpha k+i)},\end{equation} respectively. Above, we have suppressed dependence on the energy and the angle. Now, since the way we extend the initial angles in both {\it global} and {\it local} settings is via the exact same homotopy, and since the {\it initial} angle in the {\it local} setting is precisely the angle one gets at that site when solving the {\it global} difference equation, and finally, since by assumption we have $R_{-}(\alpha k-1)=1$, we conclude that $(i)$ follows. Part $(ii)$ is proven similarly; that is,

\begin{equation} R_{-L}(\alpha k-1)=R_{-L}(\alpha k+\alpha-1)\prod_{i=0}^{\alpha -1}\frac{\sin \phi_{-L}(\alpha k+i)}{\cos \phi_{-L}(\alpha k+i-1)}\end{equation}
and \begin{equation} R_{+}(\alpha k-1)=R_{+}(\alpha k+\alpha-1)\prod_{i=0}^{\alpha -1}\frac{\sin \phi_{+}(\alpha k+i)}{\cos \phi_{+}(\alpha k+i-1)}.\end{equation} Now, the result follows via a similar argument as in the previous paragraph and the fact that $R_{+}(\alpha k+\alpha -1)=1.$
\end{proof}

\begin{lemm}\label{lemmm3.3}
With the same notation as above, for $k=-L,\dots, L-1$, we have \begin{enumerate}[$(i)$]\item$\sum_{i=0}^{\alpha-1}f_iu_{-L}^2(\alpha k+i, \theta_{-L},\dots, \theta_{k},E)=R_{-L}(\alpha k-1,\theta_{-L},\dots, \theta_{k-1},E)\\
     \cdot \left(\sum_{i=0}^{\alpha-1}f_iu_{-}^2(\alpha k+i,\theta_{k-1},\lambda(\theta_{k-1},\theta_k,E),E)\right)$
\item $\sum_{i=0}^{\alpha-1}f_iu_{-L}^2(\alpha k+i, \theta_{-L},\dots, \theta_{k},E)=R_{-L}(\alpha k+\alpha-1,\theta_{-L},\dots, \theta_{k},E)\\
\cdot\left(\sum_{i=0}^{\alpha-1}f_iu_{+}^2(\alpha k+i,\theta_{k},\lambda(\theta_{k-1},\theta_k,E),E)\right)$\end{enumerate}
\end{lemm}
\begin{proof}

First note that 
\begin{align*}u_{-L}(\alpha k)&=R_{-L}(\alpha k)\sin \phi_{-L}(\alpha k)\\
&=R_{-L}(\alpha k-1)\cos\phi_{-L}(\alpha k-1)\\
&=R_{-L}(\alpha k-1)\cos \theta_{k-1}\\
&=R_{-L}(\alpha k-1)u_{-}(\alpha k).
\end{align*} 
Similarly,
\begin{align*}
u_{-L}(\alpha k+\alpha -1)&=R_{-L}(\alpha k+\alpha -1)\sin \phi_{-L}(\alpha k+\alpha-1)\\
&=R_{-L}(\alpha k+\alpha -1)\sin \theta_k\\
&=R_{-L}(\alpha k+\alpha-1)u_{+}(\alpha k+\alpha -1).
\end{align*}
Now, since given an initial angle, we extend it uniquely in exactly the same way in both the {\it global} and {\it local } setting, we see that the corresponding Pr\"ufer phases associated to $u_{-L}$, $u_{-}$ and $u_{+}$ are exactly the same, starting with $\theta_{k-1}$ or $\theta_k$ respectively. Finally, one concludes the proof via a similar argument as in Lemma $\ref{lemmm3.2}$.

\end{proof}
\vskip 8pt


\subsection{Computing the Jacobian}\label{CJ}
Before we begin to carry out the process of changing variables, note that if $v_l$ is a normalized eigenvector of $H^{(L)}_\omega\defeq H_\omega\big|_{\ell^2(-\alpha L,\dots, \alpha L-1)}$, corresponding to the energy $E_l$, then in terms of the new variables we can express it as \begin{equation}v_l(\circ)=\frac{u_{-L}(\circ,\omega(\theta_{-L},\dots,\theta_{L-2},E_l),E_l)}{\norm{u_{-L}(\circ,\omega(\theta_{-L},\dots,\theta_{L-2},E_l),E_l)}}.
\end{equation}

Using Lemmas \ref{lemm1}, \ref{lemm2} and the Feynman-Hellman formula, we get
\begin{enumerate}
\item For $j>k$:
\begin{align*}
\frac{\partial \theta_k}{\partial \omega_j}&=\frac{\partial\phi_{-L}(\alpha k+\alpha-1,\omega_{-L},\dots,\omega_k,E_l(\omega))}{\partial E}\frac{\partial E_l}{\partial \omega_j}\\
&=-\frac{1}{R_{-L}^2(\alpha k+\alpha-1,\omega,E)}\left(\sum_{i=\alpha k+\alpha-1}^{-\alpha L}u_{-L}^2(i,\omega,E)\right)\sum_{i=0}^{\alpha-1}f_iv_l^2(\alpha j+i,\omega, E)\\
&=-\frac{1}{R_{-L}^2(\alpha k+\alpha-1,\omega,E)}\frac{\sum_{i=\alpha k+\alpha-1}^{-\alpha L}u_{-L}^2(i,\omega,E)}{\sum_{i=-\alpha L}^{\alpha L-1}u_{-L}^2(i,\omega,E)}\sum_{i=0}^{\alpha-1}f_iu_{-L}^2(\alpha j+i,\omega, E)\\
\end{align*}
\item For $j\leq k$:
\begin{align*}
\frac{\partial \theta_k}{\partial \omega_j}&=\frac{\partial\phi_{-L}(\alpha k+\alpha -1,\omega_{-L},\dots,\omega_k,E_l(\omega))}{\partial \omega_j}+\frac{\partial\phi_{-L}(\alpha k+\alpha-1,\omega_{-L},\dots,\omega_k,E_l(\omega))}{\partial E}\frac{\partial E_l}{\partial \omega_j}\\
&=\frac{\sum_{i=0}^{\alpha-1}f_iu_{-L}^2(\alpha j+i,\omega, E)}{R_{-L}^2(\alpha k+\alpha-1,\omega, E)}-\frac{1}{R_{-L}^2(\alpha k+\alpha-1,\omega,E)}\frac{\sum_{i=\alpha k+\alpha-1}^{-\alpha L}u_{-L}^2(i,\omega,E)}{\sum_{i=-\alpha L}^{\alpha L-1}u_{-L}^2(i,\omega,E)}\\
&\cdot\sum_{i=0}^{\alpha-1}f_iu_{-L}^2(\alpha j+i,\omega, E)\\
&=\frac{1}{R_{-L}^2(\alpha k+\alpha-1,\omega,E)}\left(\frac{\sum_{i=\alpha k+\alpha}^{\alpha L-1}u_{-L}^2(i,\omega,E)}{\sum_{i=-\alpha L}^{\alpha L-1}u_{-L}^2(i,\omega,E)}\right)\sum_{i=0}^{\alpha-1}f_iu_{-L}^2(\alpha j+i,\omega, E).
\end{align*}
\end{enumerate}

We fix the following notation
\begin{eqnarray*}S&\defeq&\sum_{i=-\alpha L}^{\alpha L-1}u_{-L}^2(i,\omega,E);\\
S_a^b&\defeq&\sum_{i=a}^bu^2_{-L}(i,\omega,E);\\
M_j&\defeq&\sum_{i=0}^{\alpha-1}f_iu_{-L}^2(\alpha j+i,\omega, E).
\end{eqnarray*}


Now we are in a position to compute the Jacobian of the change of variables map $J_L^l$. Here $J_L^l$ represents the change of variables for a fixed $l$, as defined in $(\ref{eq:eq38})$.

\begin{equation*}
\det \partial J_L^l/\partial \omega=\det\left(\begin{array}{cccccccccc}
\frac{\partial E}{\partial \omega_{-L}}&\frac{\partial E}{\partial \omega_{-L+1}}&\frac{\partial E}{\partial \omega_{-L+2}}&\dots&\dots&\dots&\frac{\partial E}{\partial \omega_{L-3}}&\frac{\partial E}{\partial \omega_{L-2}}&\frac{\partial E}{\partial \omega_{L-1}}\\
\frac{\partial \theta_{-L}}{\partial \omega_{-L}}&\frac{\partial \theta_{-L}}{\partial \omega_{-L+1}}&\frac{\partial \theta_{-L}}{\partial \omega_{-L+2}}&\dots&\dots&\dots&\frac{\partial \theta_{-L}}{\partial \omega_{L-3}}&\frac{\partial \theta_{-L}}{\partial \omega_{L-2}}&\frac{\partial \theta_{-L}}{\partial \omega_{L-1}}\\
\frac{\partial \theta_{-L+1}}{\partial \omega_{-L}}&\frac{\partial \theta_{-L+1}}{\partial \omega_{-L+1}}&\frac{\partial \theta_{-L+1}}{\partial \omega_{-L+2}}&\dots&\dots&\dots&\frac{\partial \theta_{-L+1}}{\partial \omega_{L-3}}&\frac{\partial \theta_{-L+1}}{\partial \omega_{L-2}}&\frac{\partial \theta_{-L+1}}{\partial \omega_{L-1}}\\
\vdots&\vdots&\vdots&\vdots&\vdots&\vdots&\vdots&\vdots&\vdots\\
\frac{\partial \theta_{-1}}{\partial \omega_{-L}}&\frac{\partial \theta_{-1}}{\partial \omega_{-L+1}}&\frac{\partial \theta_{-1}}{\partial \omega_{-L+2}}&\dots&\dots&\dots&\frac{\partial \theta_{-1}}{\partial \omega_{L-3}}&\frac{\partial \theta_{-1}}{\partial \omega_{L-2}}&\frac{\partial \theta_{-1}}{\partial \omega_{L-1}}\\
\frac{\partial \theta_{0}}{\partial \omega_{-L}}&\frac{\partial \theta_{0}}{\partial \omega_{-L+1}}&\frac{\partial \theta_{0}}{\partial \omega_{-L+2}}&\dots&\dots&\dots&\frac{\partial \theta_{0}}{\partial \omega_{L-3}}&\frac{\partial \theta_{0}}{\partial \omega_{L-2}}&\frac{\partial \theta_{0}}{\partial \omega_{L-1}}\\
\vdots&\vdots&\vdots&\vdots&\vdots&\vdots&\vdots&\vdots&\vdots\\
\frac{\partial \theta_{L-3}}{\partial \omega_{-L}}&\frac{\partial \theta_{L-3}}{\partial \omega_{-L+1}}&\frac{\partial \theta_{L-3}}{\partial \omega_{-L+2}}&\dots&\dots&\dots&\frac{\partial \theta_{L-3}}{\partial \omega_{L-3}}&\frac{\partial \theta_{L-3}}{\partial \omega_{L-2}}&\frac{\partial \theta_{L-3}}{\partial \omega_{L-1}}\\
\frac{\partial \theta_{L-2}}{\partial \omega_{-L}}&\frac{\partial \theta_{L-2}}{\partial \omega_{-L+1}}&\frac{\partial \theta_{L-2}}{\partial \omega_{-L+2}}&\dots&\dots&\dots&\frac{\partial \theta_{L-2}}{\partial \omega_{L-3}}&\frac{\partial \theta_{L-2}}{\partial \omega_{L-2}}&\frac{\partial \theta_{L-2}}{\partial \omega_{L-1}}\\
\end{array}
\right)
\end{equation*}

\begin{equation*}
=\det\left(\scalemath{0.7}{\begin{array}{cccccccccc}
\frac{M_{-L}}{S}&\frac{M_{-L+1}}{S}&\frac{M_{-L+2}}{S}&\dots&&\frac{M_{L-2}}{S}&\frac{M_{L-1}}{S}\\
\frac{M_{-L}S_{-\alpha L+\alpha}^{ \alpha L-1}}{SR_{-L}^2(-\alpha L+\alpha-1)}&\frac{-M_{-L+1}S_{-\alpha L+\alpha-1}^{-\alpha L}}{SR_{-L}^2(-\alpha L+\alpha-1)}&\frac{-M_{-L+2}S_{-\alpha L+\alpha-1}^{-\alpha L}}{SR_{-L}^2(-\alpha L+\alpha-1)}&\dots&&\frac{-M_{L-2}S_{-\alpha L+\alpha-1}^{-\alpha L}}{SR_{-L}^2(-\alpha L+\alpha-1)}&\frac{-M_{L-1}S_{-\alpha L+\alpha-1}^{-\alpha L}}{SR_{-L}^2(-\alpha L+\alpha-1)}\\
\frac{M_{-L}S_{-\alpha L+2\alpha}^{ \alpha L-1}}{SR_{-L}^2(-\alpha L+2\alpha-1)}&\frac{M_{-L+1}S_{-\alpha L+2\alpha}^{\alpha L-1}}{SR_{-L}^2(-\alpha L+2\alpha-1)}&\frac{-M_{-L+2}S_{-\alpha L+2\alpha-1}^{-\alpha L}}{SR_{-L}^2(-\alpha L+2\alpha-1)}&\dots&&\frac{-M_{L-2}S_{-\alpha L+2\alpha-1}^{-\alpha L}}{SR_{-L}^2(-\alpha L+2\alpha-1)}&\frac{-M_{L-1}S_{-\alpha L+2\alpha-1}^{-\alpha L}}{SR_{-L}^2(-\alpha L+2\alpha-1)}\\
\vdots&\vdots&\vdots&\vdots&\vdots&\vdots&\\
\frac{M_{-L}S_{0}^{ \alpha L-1}}{SR_{-L}^2(-1)}&\dots&\frac{M_{-1}S_{0}^{ \alpha L-1}}{SR_{-L}^2(-1)}&-\frac{M_{0}S_{-1}^{ -\alpha L}}{SR_{-L}^2(-1)}&-\frac{M_{1}S_{-1}^{-\alpha L}}{SR_{-L}^2(-1)}&\dots&-\frac{M_{L-1}S_{-1}^{-\alpha L}}{SR_{-L}^2(-1)}\\
\frac{M_{-L}S_{\alpha}^{ \alpha L-1}}{SR_{-L}^2(\alpha-1)}&\dots&\frac{M_{-1}S_{\alpha}^{ \alpha L-1}}{SR_{-L}^2(\alpha-1)}&\frac{M_{0}S_{\alpha}^{ \alpha L-1}}{SR_{-L}^2(\alpha-1)}&-\frac{M_{1}S_{\alpha-1}^{-\alpha L}}{SR_{-L}^2(\alpha-1)}&\dots&-\frac{M_{L-1}S_{\alpha-1}^{-\alpha L}}{SR_{-L}^2(\alpha-1)}\\
\vdots&\vdots&\vdots&\vdots&\vdots&\vdots&\vdots\\
\frac{M_{-L}S_{\alpha L-2\alpha}^{ \alpha L-1}}{SR_{-L}^2(\alpha L-2\alpha-1)}&\frac{M_{-L+1}S_{\alpha L-2\alpha}^{\alpha L-1}}{SR_{-L}^2(\alpha L-2\alpha-1)}&\frac{M_{-L+2}S_{\alpha L-2\alpha}^{\alpha L-1}}{SR_{-L}^2(\alpha L-2\alpha-1)}&\dots&\frac{M_{L-3}S_{\alpha L-2\alpha}^{ \alpha L-1}}{SR_{-L}^2(\alpha L-2\alpha-1)}&\frac{-M_{L-2}S_{\alpha L-2\alpha-1}^{-\alpha L}}{SR_{-L}^2(\alpha L-2\alpha-1)}&\frac{-M_{L-1}S_{\alpha L-2\alpha-1}^{-\alpha L}}{SR_{-L}^2(\alpha L-2\alpha-1)}\\
\frac{M_{-L}S_{\alpha L-\alpha}^{ \alpha L-1}}{SR_{-L}^2(\alpha L-\alpha-1)}&\frac{M_{-L+1}S_{\alpha L-\alpha}^{\alpha L-1}}{SR_{-L}^2(\alpha L-\alpha-1)}&\frac{M_{-L+2}S_{\alpha L-\alpha}^{\alpha L-1}}{SR_{-L}^2(\alpha L-\alpha-1)}&\dots&\frac{M_{L-3}S_{\alpha L-\alpha}^{ \alpha L-1}}{SR_{-L}^2(\alpha L-\alpha-1)}&\frac{M_{L-2}S_{\alpha L-\alpha}^{\alpha L-1}}{SR_{-L}^2(\alpha L-\alpha-1)}&\frac{-M_{L-1}S_{\alpha L-\alpha-1}^{-\alpha L}}{SR_{-L}^2(\alpha L-\alpha-1)}\\
\end{array}}
\right)
\end{equation*}

Factoring out common factors in rows and columns, we get

\begin{align}\nonumber
\det \partial J_L^l/\partial \omega 
&=\frac{\prod_{j=-L}^{L-1}\sum_{i=0}^{\alpha-1}f_iu_{-L}^2(\alpha j+i,\omega,E)}{\prod_{k=-L}^{L-2}R^2_{-L}(\alpha k+\alpha -1,\omega, E)} \left(\sum_{i=-\alpha L}^{\alpha L-1}u_{-L}^2(i,\omega,E)\right)^{-2L}\times \det A,\\
\end{align}
where

\begin{equation*}\footnotesize
\arraycolsep=3pt
\medmuskip = 1mu
A=\left(\begin{array}{cccccccccc}
1&1&1&\dots&&&\dots&1&1\\
S_{-\alpha L+\alpha}^{\alpha L-1}&-S_{-\alpha L+\alpha -1}^{-\alpha L}&-S_{-\alpha L+\alpha -1}^{-\alpha L}&\dots&&&\dots&-S_{-\alpha L+\alpha -1}^{-\alpha L}&-S_{-\alpha L+\alpha -1}^{-\alpha L}\\
S_{-\alpha L+2\alpha}^{\alpha L-1}&S_{-\alpha L+2\alpha}^{\alpha L-1}&-S_{-\alpha L+2\alpha -1}^{-\alpha L}&\dots&&&\dots&-S_{-\alpha L+2\alpha -1}^{-\alpha L}&-S_{-\alpha L+2\alpha -1}^{-\alpha L}\\
\vdots&\vdots&\vdots&\vdots&\vdots&\vdots&\vdots&\vdots&\vdots\\
S_{0}^{\alpha L-1}&\dots&\dots&S_{0}^{\alpha L-1}&-S_{-1}^{-\alpha L}&-S_{-1}^{-\alpha L}&\dots&\dots&-S_{-1}^{-\alpha L}\\
S_{\alpha}^{\alpha L-1}&\dots&\dots&S_{\alpha}^{\alpha L-1}&S_{\alpha}^{\alpha L-1}&-S_{\alpha -1}^{-\alpha L}&\dots&\dots&-S_{\alpha -1}^{-\alpha L}\\
\vdots&\vdots&\vdots&\vdots&\vdots&\vdots&\vdots&\vdots&\vdots\\
S_{\alpha L-\alpha}^{\alpha L-1}&S_{\alpha L-2\alpha}^{\alpha L-1}&S_{\alpha L-2\alpha}^{\alpha L-1}&\dots&\dots&\dots&S_{\alpha L-2\alpha}^{\alpha L-1}&-S_{\alpha L-2\alpha-1}^{-\alpha L}&-S_{\alpha L-2\alpha -1}^{-\alpha L}\\
S_{\alpha L-\alpha}^{\alpha L-1}&S_{\alpha L-\alpha}^{\alpha L-1}&S_{\alpha L-\alpha}^{\alpha L-1}&\dots&\dots&\dots&S_{\alpha L-\alpha}^{\alpha L-1}&S_{\alpha L-\alpha}^{\alpha L-1}&-S_{\alpha L-\alpha -1}^{-\alpha L}\\
\end{array}
\right)
\end{equation*}

To compute the determinant of $A$ we use the following lemma. For a proof see \cite[Lemma 2.5]{ds}

\begin{lemm}\label{lem2.1}
We have
\begin{equation*}
\det\left(\begin{array}{cccccc}
a_1&a_1&a_1&\dots&a_1&a_1\\
b_2&a_2&a_2&\dots&a_2&a_2\\
b_3&b_3&a_3&\dots&a_3&a_3\\
\vdots&\vdots&\vdots&&\vdots&\vdots\\
b_n&b_n&b_n&\dots&b_n&a_n\\
\end{array}
\right)=a_1(a_2-b_2)(a_3-b_3)\dots(a_n-b_n).
\end{equation*}
\end{lemm}

Noting that $A$ has the same structure as the matrix in Lemma $\ref{lem2.1}$, we find that 

\begin{align*}
\det A=&(-S_{-\alpha L+\alpha-1}^{-\alpha L}-S_{-\alpha L+\alpha}^{\alpha L-1})(-S_{-\alpha L+2\alpha -1}^{-\alpha L}-S_{-\alpha L+2\alpha}^{\alpha L-1})\dots(-S_{-1}^{-\alpha L}-S_0^{\alpha L-1})(-S_{\alpha-1}^{-\alpha L}-S_{\alpha}^{\alpha L-1})\\
&\dots(-S_{\alpha L-2\alpha-1}^{-\alpha L}-S_{\alpha L-2\alpha}^{\alpha L-1})(-S_{\alpha L-\alpha-1}^{-\alpha L-\alpha}-S_{\alpha L-\alpha}^{\alpha L-1})\\
&=-\left(\sum_{i=-\alpha L}^{\alpha L-1}u_{-L}^2(i)\right)^{2L-1}.
\end{align*}

Thus we have proved the following lemma.

\begin{lemm}\label{lemm2.2} With the same notation and set up as above we have

$$\det \partial J_L^l/\partial \omega=-\frac{\prod_{j=-L}^{L-1}\sum_{i=0}^{\alpha-1}f_iu_{-L}^2(\alpha j+i,\omega,E)}{\prod_{k=-L}^{L-2}R^2_{-L}(\alpha k+\alpha -1,\omega, E)} \left(\sum_{i=-\alpha L}^{\alpha L-1}u_{-L}^2(i,\omega, E_l)\right)^{-1}.$$
\end{lemm}

 Let $\{v_{\bar \omega}^{L,l}\}_l$ be the normalized eigenvalues of $H_\omega^{(L)}$ corresponding to energies $\{E_l(\omega)\}_l.$
Now, we are in a position to rewrite $\rho_L(m,0)$ in terms of the new variables; that is, we have the following lemma,

\begin{lemm}\label{lem04.7}
We have $$\rho_L(m,0)\leq \int_{\Sigma_0}\rho_L(m,0,E)dE,$$ where for $E\in\Sigma_0$ we write

\begin{align*}
\rho_L(m,0,E)&\defeq \sum_{N=0}^{2B(\alpha)-1}\int_{\mathbb{T}_\alpha^{2L-1}}r(\lambda(0,\theta_{-L},E_l))r(\lambda(\theta_{L-2},N\pi+\pi/2,E_l))
\left(\prod_{n=-L+1}^{L-2}r(\lambda(\theta_{n-1},\theta_n, E_l))\right)\\ \nonumber
&\cdot \frac{R_{-L}(m,\omega, E)R_{-L}(0,\omega,E)}{\sum_{i=0}^{\alpha-1}f_iu_{-L}^2(\alpha (L-1)+i,\theta_{-L},\dots,\theta_{L-2},\pi/2+N\pi,E)}\\
&\cdot \frac{\prod_{k=-L}^{L-2}R^2_{-L}(\alpha k+\alpha -1,\theta_{-L},\dots,\theta_k, E)}{\prod_{j=-L}^{L-2}\sum_{i=0}^{\alpha-1}f_iu_{-L}^2(\alpha j+i,\theta_{-L},\dots, \theta_j,E)}d\theta_{-L}\dots d\theta_{L-2}.
\end{align*}
We interpret $r(\lambda(\cdot))$ as zero if $\lambda(\cdot,\cdot,\cdot)$ does not exist.
\end{lemm}

\begin{proof} Let $\Omega_L\defeq [-M,M]^{2L},$ and set

$$A_l\defeq \int_{R^{2L}}\chi_{\Omega_L}\left|v_{\bar{\omega}}^{L,l}(m)\right|\left|v_{\bar \omega}^{L,l}(0)\right|\prod_{n=-L}^{L-1}r(\omega_n)d\omega_{-L}\dots d\omega_{L-1}.$$ On the set $\Omega_L$ we define the change of variables as follows

\begin{eqnarray}\label{eq:eq38}
J_L^l: \Omega_L&\to& \Sigma_0\times \mathbb{T}_\alpha^{2L-1}\\
\omega&\mapsto&(E_l(\omega), \theta_{-L},\dots, \theta_{L-2}).
 \nonumber\end{eqnarray}

Pick $N\in\{0,1,\dots, 2B(\alpha)-1\}$ so that $N\equiv l\mod 2B(\alpha)$. Next, we perform the change of variables

\begin{align}\label{eq038}\nonumber
A_k&=\int_{R^{2L}}\chi_{\Omega_L}\left|v_{\bar{\omega}}^{L,l}(m)\right|\left|v_{\bar \omega}^{L,l}(0)\right|r(\omega_{-L})r(\omega_{L-1})\prod_{n=-L+1}^{L-2}r(\omega_n)d\omega_{-L}\dots d\omega_{L-1}\\ \nonumber
&=\int_{\Sigma_0}\int_{\mathbb{T}_\alpha^{2L-1}}\chi_{J_L^l(\Omega_L)}\big| \det \partial J_L^l/\partial \omega \big|^{-1}r(\lambda(0,\theta_{-L},E))r(\lambda(\theta_{L-2},N\pi+\pi/2,E))\\ \nonumber
&\prod_{n=-L+1}^{L-2}r(\omega(\theta_{n-1},\theta_n, E))
\frac{|u_{-L}(m,\omega,E)||u_{-L}(0,\omega,E)|}{\sum_{i=-\alpha L}^{\alpha L-1}u_{-L}^2(i,\omega, E)}d\theta_{-L}\dots d\theta_{L-2}d E\\ \nonumber
&=\int_{\Sigma_0}\int_{\mathbb{T}_\alpha^{2L-1}}\chi_{J_L^l(\Omega_L)}\frac{\prod_{k=-L}^{L-2}R^2_{-L}(\alpha k+\alpha -1,\theta_{-L},\dots,\theta_k, E)}{\prod_{j=-L}^{L-1}\sum_{i=0}^{\alpha-1}f_iu_{-L}^2(\alpha j+i,\omega,E)}\left(\sum_{i=-\alpha L}^{\alpha L-1}u_{-L}^2(i,\omega, E)\right) r(\lambda(0,\theta_{-L},E))\\ \nonumber
&\cdot r(\lambda(\theta_{L-2},N\pi+\pi/2,E))\prod_{n=-L+1}^{L-2}r(\lambda(\theta_{n-1},\theta_n, E))\frac{|u_{-L}(m,\omega, E)||u_{-L}(0,\omega,E)|}{\sum_{i=-\alpha L}^{\alpha L-1}u_{-L}^2(i,\omega, E)}d\theta_{-L}\dots d\theta_{L-2}d E\\ \nonumber
&=\int_{\Sigma_0}\int_{\mathbb{T}_\alpha^{2L-1}}\chi_{J_L^l(\Omega_L)}r(\lambda(0,\theta_{-L},E))r(\lambda(\theta_{L-2},N\pi+\pi/2,E))
\left(\prod_{n=-L+1}^{L-2}r(\lambda(\theta_{n-1},\theta_n, E))\right)\\ \nonumber 
&\cdot \frac{|u_{-L}(m,\omega,E)||u_{-L}(0,\omega,E)|}{\sum_{i=0}^{\alpha-1}f_iu_{-L}^2(\alpha (L-1)+i,\theta_{-L},\dots,\theta_{L-2},N\pi+\pi/2,E)}\\ \nonumber
&\cdot \frac{\prod_{k=-L}^{L-2}R^2_{-L}(\alpha k+\alpha -1,\theta_{-L},\dots,\theta_k, E)}{\prod_{j=-L}^{L-2}\sum_{i=0}^{\alpha-1}f_iu_{-L}^2(\alpha j+i,\theta_{-L},\dots,\theta_j,E)}d\theta_{-L}\dots d\theta_{L-2}d E\\ \nonumber
&\leq\int_{\Sigma_0}\int_{\mathbb{T}_\alpha^{2L-1}}\chi_{J_L^l(\Omega_L)}r(\lambda(0,\theta_{-L},E))r(\lambda(\theta_{L-2},N\pi+\pi/2,E))
\left(\prod_{n=-L+1}^{L-2}r(\lambda(\theta_{n-1},\theta_n, E))\right)\\ \nonumber
&\cdot \frac{R_{-L}(m,\omega, E)R_{-L}(0,\omega,E)}{\sum_{i=0}^{\alpha-1}f_iu_{-L}^2(\alpha (L-1)+i,\theta_{-L},\dots,\theta_{L-2},N\pi+\pi/2,E)}\\
&\cdot \frac{\prod_{k=-L}^{L-2}R^2_{-L}(\alpha k+\alpha -1,\theta_{-L},\dots,\theta_k, E)}{\prod_{j=-L}^{L-2}\sum_{i=0}^{\alpha-1}f_iu_{-L}^2(\alpha j+i,\theta_{-L},\dots, \theta_j,E)}d\theta_{-L}\dots d\theta_{L-2}d E.\\ \nonumber
\end{align} 

Now, from Lemma \ref{chov} it follows that the sets $\left\{J_L^l(\Omega_L)\right\}_{l=0}^{2\alpha L}$ are pairwise disjoint. Thus, we get

\begin{align*}
\rho_L(m,0)&=\sum_l A_l\\
&\leq\int_{\Sigma_0}\sum_l\int_{\mathbb{T}_\alpha^{2L-1}}\chi_{J_L^l(\Omega_L)}r(\lambda(0,\theta_{-L},E))r(\lambda(\theta_{L-2},N\pi+\pi/2,E))
\left(\prod_{n=-L+1}^{L-2}r(\lambda(\theta_{n-1},\theta_n, E))\right)\\ \nonumber
&\cdot \frac{R_{-L}(m,\omega, E)R_{-L}(0,\omega,E)}{\sum_{i=0}^{\alpha-1}f_iu_{-L}^2(\alpha (L-1)+i,\theta_{-L},\dots,\theta_{L-2},N\pi+\pi/2,E)}\\
&\cdot \frac{\prod_{k=-L}^{L-2}R^2_{-L}(\alpha k+\alpha -1,\theta_{-L},\dots,\theta_k, E)}{\prod_{j=-L}^{L-2}\sum_{i=0}^{\alpha-1}f_iu_{-L}^2(\alpha j+i,\theta_{-L},\dots, \theta_j,E)}d\theta_{-L}\dots d\theta_{L-2}d E.\\ 
& \leq\int_{\Sigma_0}\int_{\mathbb{T}_\alpha^{2L-1}}r(\lambda(0,\theta_{-L},E_l))r(\lambda(\theta_{L-2},N\pi+\pi/2,E))
\left(\prod_{n=-L+1}^{L-2}r(\lambda(\theta_{n-1},\theta_n, E))\right)\\ \nonumber
&\cdot \frac{R_{-L}(m,\omega, E)R_{-L}(0,\omega,E)}{\sum_{i=0}^{\alpha-1}f_iu_{-L}^2(\alpha (L-1)+i,\theta_{-L},\dots,\theta_{L-2},N\pi+\pi/2,E)}\\
&\cdot \frac{\prod_{k=-L}^{L-2}R^2_{-L}(\alpha k+\alpha -1,\theta_{-L},\dots,\theta_k, E)}{\prod_{j=-L}^{L-2}\sum_{i=0}^{\alpha-1}f_iu_{-L}^2(\alpha j+i,\theta_{-L},\dots, \theta_j,E)}d\theta_{-L}\dots d\theta_{L-2}d E.\\ 
&\leq \int_{\Sigma_0}\rho_L(m,0,E)dE,
\end{align*}
concluding the proof.
\end{proof}

\vskip 8pt



\subsection{Integral operator formula for $\rho_L(m,0, E)$}\label{IOF}

Below, we define a family of integral operators, and express $\rho_L(m,0)$ in terms of these operators. In doing so, we reduce the problem of bounding $\rho_L(m,0)$ to integral operator bounds.
\begin{defi}\label{defi3.7}
For $k=-L+1,\dots, L-2$, and $E\in \mathbb{R}$, we define a family of operators on $L^p(\mathbb{T}_\alpha)$:
\begin{displaymath}
   \left(T_{E,\alpha}^{k}F\right)(x)= \left\{
     \begin{array}{lr}
       \displaystyle \int_{\mathbb{T}_\alpha}\frac{R_{+}^2(\alpha k -1,y,\lambda(x,y,E),E)}{\sum_{i=0}^{\alpha-1}f_iu_{+}^2(\alpha k+i, y,\lambda(x,y,E),E)}r(\lambda(x,y,E))F(y)dy,&  k> 0\\ 
\displaystyle \int_{\mathbb{T}_\alpha} \frac{R_{-}^2(\alpha k+\alpha -1,y,\lambda(y,x,E),E)}{\sum_{i=0}^{\alpha-1}f_iu_{-}^2(\alpha k+i,y,\lambda(y,x,E),E)}r(\lambda(y,x,E))F(y)dy.&k\leq 0
     \end{array}
   \right.
\end{displaymath} 

and for $k>0$ we also define:

$$\left(\tilde{T}_{E,\alpha}^{k}F\right)(x)= \int_{\mathbb{T}_\alpha} \frac{R_{+}(\alpha k -1,y,\lambda(x,y,E),E)}{\sum_{i=0}^{\alpha-1}f_iu_{+}^2(\alpha k+i, y,\lambda(y,x,E),E)}r(\lambda(x,y,E))F(y)dy$$
\end{defi}

\begin{defi}\label{defi3.8.}For $E\in \mathbb{R}$, we also define the following two special functions:
\begin{eqnarray*}\psi_{E,\alpha}^{-L}(x)&=&\frac{R^{2}_{-}(-\alpha L+\alpha-1,0,\lambda(0,x,E),E)}{\sum_{i=0}^{\alpha-1}f_iu_{-}^2(-\alpha L+i,0,\lambda(0,x,E),E)}r(\lambda(0,x,E))\\
\psi_{E,\alpha,N}^L(x)&=&\frac{1}{\sum_{i=0}^{\alpha-1}f_iu_{+}^2(\alpha(L-1)+i,x,\lambda(x,N\pi+\pi/2,E),E)}r(\lambda(x,N\pi+\pi/2,E)).
\end{eqnarray*}
\end{defi}
A quick computation shows that
\begin{align}\label{eq22}\nonumber
\displaystyle \left(T_{E,\alpha}^{0}T_{E,\alpha}^{-1}\dots T_{E,\alpha}^{-L+1}\psi_{E,\alpha}^{-L}\right)(\theta_0)&=\int_{\mathbb{T}_\alpha^{L}}\frac{\prod_{k=-L+1}^{0}R_{-}^2(\alpha k+\alpha-1, \theta_{k-1},\lambda(\theta_{k-1},\theta_k,E),E)}{\prod_{j=-L+1}^{0}\left(\sum_{i=0}^{\alpha-1}f_iu_{-}^2(\alpha j+i,\theta_{j-1},\lambda(\theta_{j-1},\theta_j,E),E)\right)}\\ \nonumber
&\cdot \frac{R_{-}^2(-\alpha L+\alpha-1,0,\lambda(0,\theta_{-L},E),E)}{\sum_{i=0}^{\alpha-1}f_iu_{-}^2(-\alpha L+i,0,\lambda(0,\theta_{-L},E),E)}\\
&\cdot r(\lambda(\theta_{-L},0,E))\prod_{n=-L+1}^{0}r(\lambda(\theta_{n-1},\theta_n,E))d\,\theta_{-1}\dots d\,\theta_{-L}.
\end{align}
Similarly, 
\begin{align}\label{eq23}\nonumber
\left(\tilde T_{E,\alpha}^{1}\dots \tilde T_{E,\alpha}^{k_0},T_{E,\alpha}^{k_0+1}\dots T_{E,\alpha}^{L-2}\psi_{E,\alpha,N}^L\right)(\theta_0)&=\int_{\mathbb{T}_\alpha^{L-2}}\left(\prod_{k=1}^{k_0}\frac{R_{+}(\alpha k-1,\theta_{k},\lambda(\theta_{k-1},\theta_k,E),E)}{\sum_{i=0}^{\alpha-1}f_iu_{+}^2(\alpha k+i,\theta_{k},\lambda(\theta_{k-1},\theta_k,E),E)}\right)\\ \nonumber
\times &\left(\prod_{k=k_0+1}^{L-2}\frac{R_{+}^2(\alpha k-1,\theta_{k},\lambda(\theta_{k-1},\theta_k,E),E)}{\sum_{i=0}^{\alpha-1}f_iu_{+}^2(\alpha k+i,\theta_{k},\lambda(\theta_{k-1},\theta_k,E),E)}\right)\\
\times &\frac{r(\lambda(\theta_{L-1},N\pi,E)\prod_{n=1}^{L-2}r(\lambda(\theta_{n-1},\theta_{n},E))}{\sum_{i=0}^{\alpha-1}u_{+}^2(\alpha( L-1),\theta_{L-2},\lambda(\theta_{L-2},N\pi+\pi/2,E),E)}d\,\theta_1\dots d\,\theta_{L-2}.
\end{align}

\begin{lemm} \label{lemma4.9}With the same notation and setup as above, for some constant $\tilde C>0$, we have:
$$\rho_L(m,0,E)\leq \tilde C\sum_{N=0}^{2B(\alpha)-1}\Big\langle{\tilde T_{E,\alpha}^{1}\dots \tilde T_{E,\alpha}^{k_0}T_{E,\alpha}^{k_o+1}\dots T_{E,\alpha}^{L-2}\psi_{E,\alpha,N}^L,T_{E,\alpha}^{0}T_{E,\alpha}^{-1}\dots T_{E,\alpha}^{-L+1}\psi_{E,\alpha}^{-L}\Big\rangle}_{L^2(\mathbb{T}_\alpha,d\theta_0)},$$ where by $\big\langle{\cdot ,\cdot \big\rangle}$ we have denoted the inner product on $L^2(\mathbb{T}_\alpha),$ and where $m=\alpha k_0+j$, for some $j\in\{0,\dots,\alpha-1\}.$
\end{lemm}

\begin{proof}
The proof is essentially a combination of Lemmas $\ref{lemmm3.2}, \,\ref{lemmm3.3}$, the results in $(\ref{eq22})$, $(\ref{eq23})$, and the discussion below. Specifically, from Lemmas $\ref{lemmm3.2}$ and $\ref{lemmm3.3}$ we get:

\begin{eqnarray*}
\frac{R_{-L}^2(\alpha k+\alpha-1, \theta_{-L},\dots, \theta_k,E)}{\sum_{i=0}^{\alpha-1}f_iu_{-L}^2(\alpha k+i, \theta_{-L},\dots,\theta_{k},E)}&=&\frac{R_{-}^2(\alpha k+\alpha-1,\theta_{k-1},\lambda(\theta_{k-1},\theta_k,E),E)}{\sum_{i=0}^{\alpha-1}f_iu_{-}^2(\alpha k+i,\theta_{k-1},\lambda(\theta_{k-1},\theta_k,E),E)};\\
\frac{R_{-L}^2(\alpha k-1, \theta_{-L},\dots, \theta_{k-1},E)}{\sum_{i=0}^{\alpha-1}f_iu_{-L}^2(\alpha k+i, \theta_{-L},\dots,\theta_{k},E)}&=&\frac{R_{+}^2(\alpha k-1,\theta_{k},\lambda(\theta_{k-1},\theta_k,E),E)}{\sum_{i=0}^{\alpha-1}f_iu_{+}^2(\alpha k+i,\theta_{k},\lambda(\theta_{k-1},\theta_k,E),E)};\\
\end{eqnarray*}
and 
\[\frac{R_{-L}(\alpha k-1, \theta_{-L},\dots, \theta_{k-1},E)R_{-L}(\alpha k+\alpha-1, \theta_{-L},\dots, \theta_k,E)}{\sum_{i=0}^{\alpha-1}f_iu_{-L}^2(\alpha k+i, \theta_{-L},\dots,\theta_{k},E)}\]
\begin{align*}
&=\frac{R_{-L}^2(\alpha k-1, \theta_{-L},\dots, \theta_{k-1},E)}{\sum_{i=0}^{\alpha-1}f_iu_{-L}^2(\alpha k+i, \theta_{-L},\dots,\theta_{k},E)}\cdot\frac{R_{-L}(\alpha k+\alpha-1, \theta_{-L},\dots, \theta_k,E)}{R_{-L}(\alpha k-1, \theta_{-L},\dots, \theta_{k-1},E)}\\
&=\frac{R_{+}^2(\alpha k-1,\theta_{k},\lambda(\theta_{k-1},\theta_k,E),E)}{\sum_{i=0}^{\alpha-1}f_iu_{+}^2(\alpha k+i,\theta_{k},\lambda(\theta_{k-1},\theta_k,E),E)}\cdot \frac{1}{R_{+}(\alpha k-1,\theta_k,\lambda(\theta_{k-1},\theta_k,E),E)}\\
&=\frac{R_{+}(\alpha k-1,\theta_{k},\lambda(\theta_{k-1},\theta_k,E),E)}{\sum_{i=0}^{\alpha-1}f_iu_{+}^2(\alpha k+i,\theta_{k},\lambda(\theta_{k-1},\theta_k,E),E)}.
\end{align*}

Next, suppose that $\alpha (k_0-1)+\alpha-1<m\leq \alpha k_0+\alpha-1$, for some $ -L\leq k_0\leq L-1;$ that is $m=\alpha k_0+j$, for some $j\in\{0,\dots,\alpha-1\}.$ We rewrite the integrand in $(\ref{eq038})$, that is,

\[
\frac{R_{-L}(m,\omega, E)R_{-L}(0,\omega,E)\prod_{k=-L}^{L-2}R^2_{-L}(\alpha k+\alpha -1,\theta_{-L},\dots,\theta_k, E)}{\sum_{i=0}^{\alpha-1}f_iu_{-L}^2(\alpha (L-1)+i,\theta_{-L},\dots,\theta_{L-2},N\pi+\pi/2,E)\prod_{k=-L}^{L-2}\sum_{i=0}^{\alpha-1}f_iu_{-L}^2(\alpha k+i,\theta_{-L},\dots, \theta_k,E)}\]
\[=R_{+}(0,\theta_0,\lambda(\theta_{-1},\theta_0,E),E)R_{+}(m,\theta_{k_0},\lambda(\theta_{k_0-1},\theta_{k_0},E),E)\]
\[\cdot\left(\frac{R_{-L}(\alpha-1,\theta_{-L},\dots,\theta_0,E)R_{-L}(\alpha k_0+\alpha-1,\theta_{-L},\dots,\theta_{k_0},E)\prod_{k=-L}^{L-2}R^2_{-L}(\alpha k+\alpha -1,\theta_{-L},\dots,\theta_k, E)}{\sum_{i=0}^{\alpha-1}f_iu_{-L}^2(\alpha (L-1)+i,\theta_{-L},\dots,\theta_{L-2},N\pi+\pi/2,E)\prod_{k=-L}^{L-2}\sum_{i=0}^{\alpha-1}f_iu_{-L}^2(\alpha k+i,\theta_{-L},\dots, \theta_k,E)}\right)\]
\[=R_{+}(0,\theta_0,\lambda(\theta_{-1},\theta_0,E),E)R_{+}(m,\theta_{k_0},\lambda(\theta_{k_0-1},\theta_{k_0},E),E)
\cdot\left(\prod_{k=-L}^0\frac{R_{-L}^2(\alpha k+\alpha-1,\theta_{-L},\dots,\theta_k,E)}{\sum_{i=0}^{\alpha-1}f_iu_{-L}^2(\alpha k+i,\theta_{-L},\dots, \theta_k,E)}\right)\] \[\cdot \left(\prod_{k=1}^{k_0}\frac{R_{-L}(\alpha k-1, \theta_{-L},\dots, \theta_{k-1},E)R_{-L}(\alpha k+\alpha-1, \theta_{-L},\dots, \theta_k,E)}{\sum_{i=0}^{\alpha-1}f_iu_{-L}^2(\alpha k+i,\theta_{-L},\dots, \theta_k,E)}\right)\]
\[\cdot \left(\prod_{k=k_0+1}^{L-1}\frac{R_{-L}^2(\alpha k-1,\theta_{-L},\dots,\theta_{k-1},E)}{\sum_{i=0}^{\alpha-1}f_iu_{-L}^2(\alpha k+i,\theta_{-L},\dots, \theta_k,E)}\right)\]
\[=R_{+}(0,\theta_0,\lambda(\theta_{-1},\theta_0,E),E)R_{+}(m,\theta_{k_0},\lambda(\theta_{k_0-1},\theta_{k_0},E),E)\]
\[\cdot\left(\prod_{k=-L}^0\frac{R_{-}^2(\alpha k+\alpha-1,\theta_{k-1},\lambda(\theta_{k-1},\theta_k,E),E)}{\sum_{i=0}^{\alpha-1}f_iu_{-}^2(\alpha k+i,\theta_{k-1},\lambda(\theta_{k-1},\theta_k,E),E)}\right)\left(\prod_{k=1}^{k_0}\frac{R_{+}(\alpha k-1,\theta_{k},\lambda(\theta_{k-1},\theta_k,E),E)}{\sum_{i=0}^{\alpha-1}f_iu_{+}^2(\alpha k+i,\theta_{k},\lambda(\theta_{k-1},\theta_k,E),E)}\right)\]
\[\cdot \left(\prod_{k=k_0+1}^{L-1}\frac{R_{+}^2(\alpha k-1,\theta_{k},\lambda(\theta_{k-1},\theta_k,E),E)}{\sum_{i=0}^{\alpha-1}f_iu_{+}^2(\alpha k+i,\theta_{k},\lambda(\theta_{k-1},\theta_k,E),E)}\right).\\
\]
The only remaining task is to produce uniform bounds for $R_{+}(0,\theta_0,\lambda(\theta_{-1},\theta_0,E),E)$ and 

\noindent $R_{+}(m,\theta_{k_0},\lambda(\theta_{k_0-1},\theta_{k_0},E),E)$, which we do in the following lemma. This concludes the proof!

\end{proof}

\begin{lemm}\label{lemm3.10}
For $k=-L,\dots, L-1,$ and $ j=0,1,\dots, \alpha-1$, we have \begin{enumerate}[$(i)$]
\item $R_{-}^2(\alpha k+j, y,\lambda(y,x,E),E)\leq C_1<\infty,$
\item $R_{+}^2(\alpha k+j, x,\lambda(x,y,E),E)\leq C_2<\infty,$
\item $\sum_{i=0}^{\alpha-1}f_iu_{+}^2(\alpha k+i, x,\lambda(y,x,E),E)\geq C_3>0,$
\item $\sum_{i=0}^{\alpha-1}f_iu_{-}^2(\alpha k+i,y,\lambda(y,x,E),E)\geq C_4>0$.
\end{enumerate}
 uniformly in $E\in\Sigma_0$ and $x,y$ such that $\lambda(y,x)\in [-M,M],$ and $j,k$.
\end{lemm}
\begin{proof}

Here, we only prove $(i)$ and $(iv)$; since the proofs of $(ii)$ and $(iii)$ are completely analogous, we omit them.
Recall that \begin{equation}\label{eq37}
\left(
\begin{array}{c}
u_{-}(\alpha k+j+1, y,\lambda, E)\\
u_{-}(\alpha k+j,y,\lambda, E)
\end{array}
\right)=T_{E,\omega}(\alpha k+j)\dots T_{E,\omega}(\alpha k)\left(\begin{array}{c}u_{-}(\alpha k,y,\lambda,E)\\ u_{-}(\alpha k-1,y,\lambda,E)\end{array}\right),
\end{equation}
where \begin{equation*}
T_{E,\omega}(\alpha k+i)=\left(\begin{array}{cc}E-f_i\omega_k&-1\\1&0\end{array}\right)
\end{equation*} for $i=0,\dots,j.$
Let $A_{E,\omega}(\alpha k+j)=T_{E,\omega}(\alpha k+j)\dots T_{E,\omega}(\alpha k)$, then $(\ref{eq37})$ becomes

\begin{equation}\label{eq38}
\left(
\begin{array}{c}
u_{-}(\alpha k+j+1, y,\lambda, E)\\
u_{-}(\alpha k+j,y,\lambda, E)
\end{array}
\right)=A_{E,\omega}(\alpha k+j)\left(\begin{array}{c}\cos y\\ \sin y\end{array}\right),
\end{equation}where we have used the fact that our solution $u_{-}$ is normalized at the left end-point. Let $\norm{\cdot}_{\mathrm{HS}}$ denote the Hilbert-Schmidt norm of an operator, then we have 

\begin{align}\label{eq40} \nonumber
\norm{T_{E,\omega}(\alpha k+i)}^2&\leq \norm{T_{E,\omega}(\alpha k+i)}_{\mathrm{HS}}^2\\ \nonumber
&=|E-f_i\omega_k|^2+2\\ \nonumber
&\leq \left(2+\left(1+\max_{0\leq i\leq j}f_i\right)M\right)^2+2\\
&\leq \left(2+\left(1+\max_{0\leq i\leq \alpha-1}f_i\right)M\right)^2+2\\
&\defeq \bar C_1<\infty.
\end{align}
where we have used the fact that $\supp r\in[-M,M]$, and $E\in \Sigma_0=[-2+M,2+M].$
\noindent
Now, from $(\ref{eq38})$ and $(\ref{eq40})$ it follows that 
\begin{equation}\label{eq41}
|u_{-}(\alpha k+j+1,y,\lambda, E)|^2+|u_{-}(\alpha k+j, y,\lambda, E)|^2\leq \bar C_1^{j+1}\leq \bar C_1^\alpha\defeq C_1. 
\end{equation}
Then, by the definition of the Pr\"{u}fer amplitude $R_{-}$, from $(\ref{eq41})$ it follows that 

$$R_{-}^2(\alpha k+j,y,\lambda, E)\leq C_1,$$ as desired!

To prove $(iv)$, note that from the analogous expressions of $(\ref{eq38})$, for $i=0,\dots, \alpha-1$, we also get
\begin{align}\label{eq42}\nonumber
|u_{-}(\alpha k+i,y,\lambda, E)|^2+|u_{-}(\alpha k+i+1, y,\lambda, E)|^2&\geq \frac{1}{\norm{A_{E,\omega}(\alpha k+i)^{-1}}}\\
&=\frac{1}{\norm{A_{E,\omega}(\alpha k+i)}}\\
&\geq \frac{1}{\bar C_1^\alpha}\\
&\defeq \bar C_4>0.
\end{align}
As a result, since $f_i>0$ for all $i=0,\dots,\alpha-1$, there is some $C_4>0$ such that

$$\sum_{i=0}^{\alpha-1}f_iu_{-}^2(\alpha k+i,y,\lambda(y,x,E),E)\geq C_4,$$ thus, proving $(iv)$.
\end{proof}

\section{Norm estimates for the integral operators}

In this section we will provide estimates on the norms of the integral operators defined in Definition $\ref{defi3.7}$. 

\begin{defi}We denote the norm of a linear operator $T$ from $L^p(\mathbb{T}_\alpha)$ to $L^q(\mathbb{T}_\alpha)$ by $\norm{T}_{p,q}.$
\end{defi}

Following the line of argument for the continuum case in \cite{ds}, we state and prove the following series of lemmas.

\begin{lemm}\label{lemma05.2}
We have

\begin{enumerate}[$(i)$]

\item $\norm{T_{E,\alpha}^{k}}_{1,2}\leq C_0<\infty$ uniformly in $E\in \Sigma_0$ and $k.$
\item $\norm{\psi_{E,\alpha}^{-L}}_1\leq C_1<\infty$ and $\norm{\psi_{E,\alpha}^{L}}_1\leq C_2<\infty$, uniformly in $E\in\Sigma_0$, where $\norm{\cdot}_1$ denotes the $L^1-$norm on $\mathbb{T}_\alpha$. 
\end{enumerate}
\end{lemm}

\begin{proof}
We prove the case when $k\leq 0$. The case when $k>0$ is proved in an identical fashion. 

Let $F\in L^1(\mathbb{T}_\alpha)$, then 

\begin{align*}
\norm{T_{E,\alpha}^{k}F}_2^2&=\int_{\mathbb{T}_\alpha}\left|\left(T_{E,\alpha}^{k}F\right)(x)\right|^2dx\\
&=\int_{\mathbb{T}_\alpha}\left|\int_{\mathbb{T}_\alpha} \frac{R_{-}^2(\alpha k+\alpha -1,y,\lambda(y,x,E),E)}{\sum_{i=0}^{\alpha-1}f_iu_{-}^2(\alpha j+i,y,\lambda(y,x,E),E)}r(\lambda(y,x,E))F(y)dy\right|^2dx\\
&\leq\int_{\mathbb{T}_\alpha}\left(\frac{C_1}{C_4}\right)^2\norm{r}_\infty^2\int_{\mathbb{T}_\alpha}\left|F(y)\right|dy dx\\
&=2\pi B(\alpha)\left(\frac{C_1}{C_4}\right)^2\norm{r}_\infty^2\norm{F}_1^2\\
&\defeq C_0,
\end{align*}
where the inequality in the third line follows from Lemma $\ref{lemm3.10}$ and elementary estimates. This concludes the proof of $(i)$. Next, making use of Lemma $\ref{lemm3.10}$ once again, we get
\begin{align*}
\norm{\psi_{E,\alpha}^{-L}}_1&=\int_{\mathbb{T}_\alpha}\left|\psi_{E,\alpha}^{-L}(x)\right|dx\\
&=\int_{\mathbb{T}_\alpha}\left|\frac{R^{2}_{-}(-\alpha L+\alpha-1,0,\lambda(0,x,E),E)}{\sum_{i=0}^{\alpha-1}f_iu_{-}^2(-\alpha L+i,0,\lambda(0,x,E),E)}r(\lambda(0,x,E))\right| dx\\
&\leq2\pi B(\alpha)\frac{C_1}{C_4}\norm{r}_\infty\\
&\defeq C_1.
\end{align*}
Similarly for $\psi_{E,\alpha}^L.$
\end{proof}

\begin{lemm}\label{lemm4.3.}For all $-L+1\leq k\leq L-2$ and $E\in \Sigma_0$ we have $$\norm{T_{E,\alpha}^{k}}_{1,1}=1$$

\end{lemm}

\begin{proof}We consider the case $k\leq 0.$ Let $F\in L^1(\mathbb{T}_\alpha)$, then

\begin{align}\label{eq43}\nonumber
\norm{T_{E,\alpha}^{k}F}_1&=\int_{\mathbb{T}_\alpha}\left|\left(T_{E,\alpha}^{k}F)(x)\right)\right|dx\\\nonumber
&=\int_{\mathbb{T}_\alpha}\left|\int_{\mathbb{T}_\alpha}\frac{R_{-}^2(\alpha k+\alpha -1,y,\lambda(y,x,E),E)}{\sum_{i=0}^{\alpha-1}f_iu_{-}^2(\alpha k+i,y,\lambda(y,x,E),E)}r(\lambda(y,x,E))F(y)dy\right|dx\\\nonumber
&\leq \int_{\mathbb{T}_\alpha}\int_{\mathbb{T}_\alpha}\frac{R_{-}^2(\alpha k+\alpha -1,y,\lambda(y,x,E),E)}{\sum_{i=0}^{\alpha-1}f_iu_{-}^2(\alpha k+i,y,\lambda(y,x,E),E)}r(\lambda(y,x,E))\left|F(y)\right|dy\,dx\\
&=\int_{\mathbb{T}_\alpha}\left(\int_{\mathbb{T}_\alpha}\frac{R_{-}^2(\alpha k+\alpha -1,y,\lambda(y,x,E),E)}{\sum_{i=0}^{\alpha-1}f_iu_{-}^2(\alpha k+i,y,\lambda(y,x,E),E)}r(\lambda(y,x,E))dx\right)\left|F(y)\right|dy.\\\nonumber
\end{align}

For fixed $y$, let us introduce the following change of variables $J_y: \lambda\to x$, where $\lambda$ comes from the set of coupling constants such that $\lambda(y,x,E)=\lambda$ (i.e. we are really thinking of $x$ as the inverse function of $\lambda(y,x,E)$ for a fixed $y$; it is not difficult to see that such an inverse exists, since $\lambda(y,x, E)$ is strictly monotonic in $x$ for each fixed $y$). From $\phi_{-}(\alpha k+\alpha-1, y,\lambda, E)=x$ and Lemma $\ref{lemm1}$ we get

\begin{align}\label{eq45}\nonumber
\frac{\partial x}{\partial \lambda}&=\frac{\partial \phi_{-}(\alpha k+\alpha-1,y,\lambda, E)}{\partial \lambda}\\
&=R_{-}^{-2}(\alpha k+\alpha -1,y,\lambda,E)\sum_{i=0}^{\alpha-1}f_iu_{-}^2(\alpha k+i,y,\lambda,E).
\end{align}
Then, after carrying out the change of variables $(\ref{eq43})$ becomes

\begin{equation*}
\int_{\mathbb{T}_\alpha}\left(\int_{\mathbb{R}}r(\lambda)d\lambda\right)\left|F(y)\right|=\norm{r}_1\norm{F}_1=\norm{F}_1.
\end{equation*}
This shows that $\norm{T_{E,\alpha}^{k}}_{1,1}\leq 1.$

Now, since for $F\geq 0$ in $(\ref{eq43})$ we have all equal signs, it follows that, indeed, $\norm{T_{E,\alpha}^{k}}_{1,1}=1,$ as claimed. In a completely similar way, one proves the case for $k>0.$
\end{proof}

\begin{defi}\label{defi5.4} For $x,y\in\mathbb{T}_\alpha$, and $k\geq 0$, let 
\begin{displaymath}
   K_{E,\alpha}^{k}(x,y)= \left\{
     \begin{array}{lr}
       \displaystyle \frac{R_{+}(\alpha k -1,y,\lambda(x,y,E),E)}{\sum_{i=0}^{\alpha-1}f_iu_{+}^2(\alpha k+i, y,\lambda(x,y,E),E)}r(\lambda(x,y,E)),&\text{if } \lambda(x,y,E)\, \text{exists} \\ 
0& \text{otherwise}
     \end{array}
   \right.
\end{displaymath} be the integral kernel for the operator $\tilde T_{E,\alpha}^{k}$. 
\end{defi}

\begin{lemm}\label{lemma4.5}\hspace{1cm}

\begin{enumerate}[$(a)$]
  \item For all $x,y\in\mathbb{T}_{\alpha}$ and all $k\geq 0$, we have $$K_{E,\alpha}^{k}(x+\pi,y+\pi)=K_{E,\alpha}^{k}(x,y).$$
  \item $K_{E,\alpha}^k(\cdot,\cdot)$ is almost everywhere continuous on $\mathbb{T}_\alpha^2.$
  \end{enumerate}
\end{lemm}
\begin{proof}
We begin by proving the first claim. In general, by linearity one has $u_{+}(\cdot, y+\pi,\lambda, E)=-u_+(\cdot, y,\lambda, E).$ From this it follows that $R_+(\cdot, y+\pi,\lambda, E)=R_+(\cdot, y,\lambda, E)$, and $\phi_+(\cdot, y+\pi,\lambda, E)=\phi_+(\cdot, y,\lambda, E)+\pi$, where the latter follows by an identical argument as the one in the paragraph preceding  $(\ref{eqeq1})$. Now, from this observation, and the definition of $\lambda(x,y,E)$ it immediately follows that $\lambda(x+\pi,y+\pi,E)=\lambda(x,y,E).$ That is, if for any $x,y\in \mathbb{T}_\alpha$ there is a coupling constant $\lambda$ such that $\phi_+(\alpha k-1,y,\lambda, E)=x$, we set $\lambda(x,y,E)=\lambda.$ Now, from the argument above, since given $x+\pi,y+\pi\in\mathbb{T}_\alpha$, for the same coupling constant $\lambda$, we have $\phi_{+}(\alpha k-1, y+\pi,\lambda, E)=\phi_+(\alpha k-1, y,\lambda, E)+\pi=x+\pi$, hence $\lambda(x+\pi,y+\pi,E)=\lambda,$ as desired. Finally, the claim is merely a combination of these facts!

We will prove continuity of $K_{E,\alpha}^k(\cdot,\cdot)$ on $\mathbb{T}_\alpha^2\setminus \lambda_E^{-1}(B)$, where $B$ denotes the set of discontinuities of the density function $r$. To that end, we need to first show that \begin{equation}\label{eq47'}\mathcal{D}\defeq\{(y,x)\in\mathbb{T}_\alpha^2:\lambda(y,x) \, \text{exists}\} \,\text{is open and }\, \lambda(\cdot,\cdot) \,\text{ is continuous on}\, \mathcal{D}.\end{equation} 

Let $(y,x)\in \mathcal{D}$ be some fixed point and $\epsilon>0$, be given. We first fix $y$ and increase $x$ a little; that is, for some $\delta_1>0$, we consider $x+\delta_1$. Since, $(y,x)\in\mathcal{D}$ , we have $\bar \lambda\defeq\lambda(y,x)$, for some $\bar\lambda$, that is $\phi_+(\alpha k-1,y,\bar \lambda, E)=x$, and since $\phi_+(\alpha k-1,y,\lambda,E)$ is continuous in $\lambda$ it means that for small enough $\delta_1$ there will be some $\lambda_1$ such that  $\phi_+(\alpha k-1,y, \lambda_1, E)=x+\delta_1,$ which is the same as $\lambda(y,x+\delta_1)=\lambda_1.$ On the other hand, since $\phi_+$ is strictly increasing in $\lambda$, for every $0<\delta< \delta_1$ we will have $\phi_+(\alpha k-1,y,\bar\lambda_1,E)=x+\delta$, for some $\bar\lambda<\bar\lambda_1<\lambda_1.$ We can actually pick $\delta_1$ small enough so that $\lambda(y,x+\delta_1)\leq \lambda(y,x)+\epsilon/2.$ So, from above, it follows that also \begin{equation}\label{eq47}\bar\lambda_1\defeq\lambda(y,x+\delta)\leq\lambda(y,x)+\epsilon/2, \,\text{for all}\, 0<\delta\leq\delta_1\end{equation} Via a similar argument, it follows that if we keep $x+\delta_1$ fixed and decrease $y$ a little, we can find a small enough $\delta_2>0$ so that $\lambda(y-\delta_2,x+\delta_1)\leq\lambda(y,x+\delta_1)+\epsilon/2$. Putting these two together, we get $\lambda(y-\delta_2,x+\delta_1)\leq\lambda(y,x)+\epsilon.$ Again, via an almost identical argument, we can find small enough $\delta_3,\delta_4>0$ so that $\lambda(y+\delta_4,x-\delta_3)\geq \lambda(y,x)-\epsilon.$ We wish to note that, using the fact that $\phi_+$ is increasing in $\lambda$, a similar argument as in $(\ref{eq47})$ can be made in the other three cases above, as well. Let, $B(y,x)\defeq\{(y+\delta,x+\tilde \delta): -\delta_2\leq \delta\leq \delta_4, -\delta_1\leq \tilde\delta\leq \delta_1\}$. The fact that $B(y,x)\subset\mathcal{D},$ follows immediately from $(\ref{eq47})$ and the last comment! From above it also follows that $\lambda(B(y,x))\subset [\lambda(y,x)-\epsilon,\lambda(y,x)+\epsilon]$, which proves continuity of $\lambda(\cdot,\cdot).$ To prove $Leb-$almost everywhere continuity of $K_{E,\alpha}^{k}(y,x)$ we first rewrite it using the closed subset $A\defeq\lambda^{-1}(\supp r)$ of $\mathcal{D},$

\begin{displaymath}
   K_{E,\alpha}^{k}(y,x)= \left\{
     \begin{array}{lr}
       \displaystyle \frac{R_{+}(\alpha k -1,x,\lambda(y,x,E),E)}{\sum_{i=0}^{\alpha-1}f_iu_{+}^2(\alpha k+i, x,\lambda(y,x,E),E)}r(\lambda(y,x,E)),&\text{if } (y,x)\in\mathcal{D} \\ 
0& \text{if }\,  (y,x)\in\mathbb{T}^2_\alpha\setminus A
     \end{array}
   \right.
\end{displaymath} First, we note that this is well-defined, since if $(y,x)\in \mathcal{D}\setminus A$, then $r(\lambda(y,x))=0$, so there is no ambiguity. Since $\{\mathcal{D}, \mathbb{T}_\alpha^2\setminus A\}$ forms an open cover of $\mathbb{T}_\alpha^2$ it suffices to prove continuity on each of the open sets separately. Continuity on $\mathbb{T}_\alpha^2\setminus A$ is obvious. To prove continuity on $\mathcal{D}\setminus\lambda_E^{-1}(B)$, first note that it is not hard to see that $u_+$ is jointly continuous on $(x,\lambda)$; this essentially follows from the fact that $u_+(\alpha k+i,x,\lambda,E)$ is a polynomial in $\cos x,\,\sin x$, and $\lambda$. As a result, $R_+(\alpha k-1,x,\lambda,E)$ is also jointly continuous on $(x,\lambda)$. Now, the result follows from $(\ref{eq47'})$ and the fact that $r$ is continuous away from $B$. Finally, since, by assumption on $r$, $Leb(B)=0$, from Proposition $\ref{Leb}$ it follows that $Leb(\lambda_E^{-1}(B))=0$, and thus $K_{E,\alpha}^k(\cdot,\cdot)$ is almost everywhere continuous on $\mathbb{T}_\alpha^2.$
\end{proof}

\begin{lemm}\label{lemma4.6}
We have $$\norm{\tilde T_{E,\alpha}^{k}}_{2,2}\leq 1.$$
\end{lemm}
\begin{proof}
First let us define $$K_1^{k}(x,y)=\frac{r(\lambda(x,y,E))}{\sum_{i=0}^{\alpha-1}u_+^2(\alpha k+i,y,\lambda(x,y,E),E)},$$ and 
$$K_2^{k}(x,y)=\frac{R_+^2(\alpha k-1,y,\lambda(x,y,E),E)r(\lambda(x,y,E))}{\sum_{i=0}^{\alpha-1}f_iu_+^2(\alpha k+i,y,\lambda(x,y,E),E)},$$ whenever $\lambda(x,y,E)$ exists, else we set them both equal to zero. 

Next, we compute
\begin{align*}
\int_{\mathbb{T}_\alpha} K_2^{k}(x,y)dx&=\int_{\mathbb{T}_\alpha}\frac{R_+^2(\alpha k-1,y,\lambda(x,y,E),E)r(\lambda(x,y,E))}{\sum_{i=0}^{\alpha-1}f_iu_+^2(\alpha k+i,y,\lambda(x,y,E),E)}dx.
\end{align*}
Similarly as in Lemma $\ref{lemm4.3.}$ we introduce the following change of variables $J:\lambda\to x.$  From $\phi_+(\alpha k-1,y,\lambda, E)=x$, we get 
\begin{align*}
\frac{\partial x}{\partial \lambda}&=\frac{\partial \phi_{+}(\alpha k-1,y,\lambda, E)}{\partial \lambda}\\
&=-R_{+}^{-2}(\alpha k -1,y,\lambda,E)\sum_{i=0}^{\alpha-1}f_iu_{+}^2(\alpha k+i,y,\lambda,E).
\end{align*} Then, we get
\begin{align*}
\int_{\mathbb{T}_\alpha} K_2^{k}(x,y)dx&=\int_{\mathbb{T}_\alpha}\frac{R_+^2(\alpha k-1,y,\lambda,E)r(\lambda)}{\sum_{i=0}^{\alpha-1}f_iu_+^2(\alpha k+i,y,\lambda,E)}\left|\frac{d x}{d\lambda}\right|d\lambda\\
&=\int_{\mathbb{T}_\alpha}r(\lambda)d\lambda\\
&=1.
\end{align*}
To compute the second integral, we note that from Lemmas $\ref{lemmm3.2}$ and $\ref{lemmm3.3}$ we get $$ \sum_{i=0}^{\alpha-1}f_iu_+^2(\alpha k+i,y,\lambda(x,y,E),E)=R_{-}^{-2}(\alpha k+\alpha-1,x,\lambda(x,y,E),E)\sum_{i=0}^{\alpha-1}f_iu_-^2(\alpha k+i,x,\lambda(x,y,E),E).$$ 
So,
\begin{align*}
\int_{\mathbb{T}_\alpha}K_1^{k}(x,y)dy&=\int_{\mathbb{T}_\alpha}\frac{r(\lambda(x,y,E))}{\sum_{i=0}^{\alpha-1}f_iu_+^2(\alpha k+i,y,\lambda(x,y,E),E)}\\
&=\int_{\mathbb{T}_\alpha}\frac{R_{-}^{2}(\alpha k+\alpha-1,x,\lambda(x,y,E),E)}{\sum_{i=0}^{\alpha-1}f_iu_-^2(\alpha k+i,x,\lambda(x,y,E),E)}r(\lambda(x,y,E))dy\\
&=1,
\end{align*}
where the last step follows after performing the same change of variables as in Lemma $\ref{lemm4.3.}.$

Finally, noting that $K_{E,\alpha}^{k}(x,y)=\sqrt{K_1^{k}(x,y)}\sqrt{K_2^{k}(x,y)}$, the result follows immediately by the Schur Test; specifically, the version that appears in \cite{jw}. 
\end{proof}


\section{The $\norm{\cdot}_{2,2}-$norm of $\tilde T_{E,\alpha}^{k}$}

The goal of this section is to prove the following proposition.

\begin{prop}\label{prop5.1}
 For all $k\geq 0$ we have $$\norm{\tilde T_{E,\alpha}^{k}}_{2,2}<1.$$
\end{prop}

The fact that these operators are $\pi-$periodic, as established in Lemma $\ref{lemma4.5}$, suggests that the operators $\tilde T_{E,\alpha}^{k}$ might be decomposable into a direct sum of integral operators on $L^2(0,\pi).$ This is established in the following lemma. Part $(a)$ of the lemma is common knowledge and the rest is the exact analogue of the continuum version, however, for completeness we provide the statement and its proof with the corresponding modifications and adjustments to the discrete setting. 

\begin{lemm}\label{lemma5.2}
\begin{enumerate}[$(a)$] 
\item Suppose $h$ is continuous on $(\pi n,\pi(n+1))$ for $n=0,1,\dots, 2B(\alpha)-1,$ and for $j\in\{0,1,\dots,2B(\alpha)-1\}$, and $x\in(0,\pi)$, let $$\displaystyle(Uh)_j(x)=\frac{1}{\sqrt{2B(\alpha)}}\sum_{n=0}^{2B(\alpha)-1}e^{-\frac{i\pi j n}{B(\alpha)}}h(x+\pi n).$$ Then $U$ extends to a unitary operator $$U:L^2(\mathbb{T}_{\alpha})\to \bigoplus_{j=0}^{2B(\alpha)-1}L^2(0,\pi).$$
\item We have $$U\tilde T_{E,\alpha}^{k}U^{-1}=\bigoplus_{j=0}^{2B(\alpha)-1} L_{j,E}^k,$$ where $L_{j,E}^k$ is the integral operator on $L^2(0,\pi)$ with kernel $$ L_{E,\alpha, j}^{k}(x,y)=\sum_{n=0}^{2B(\alpha)-1}K_{E,\alpha}^{k}(x,y+n\pi)e^{\frac{i\pi j n}{B(\alpha)}}.$$
\item We have $\norm{\tilde T_{E,\alpha}^{k}}_{2,2}=\norm{L_{0,E}^k}_{2,2}$ where the norms are taken in the respective $L^2-$ spaces. 
\end{enumerate}
\end{lemm}

\begin{proof}
\begin{enumerate}[$(a)$]
\item First, we will show that $U$ is densely defined with dense image. To this end, let $g=(g_m)\in \bigoplus_{m=0}^{2B(\alpha)-1} L^2(0,\pi)$, be any given continuous function. Let us define, $$h(x+\pi n)=\frac{1}{\sqrt{2B(\alpha)}}\sum _{m=0}^{2B(\alpha)-1} e^{\frac{i\pi m n}{B(\alpha)}} g_m(x),$$ for $x\in (0,\pi)$ and $n\in \{0,1,\dots, 2B(\alpha)-1\}.$ We claim that $Uh=g.$ That is, 

\begin{align*}
(Uh)_j(x)&=\frac{1}{\sqrt{2B(\alpha)}}\sum_{n=0}^{2B(\alpha)-1}e^{-\frac{i\pi j n}{B(\alpha)}}h(x+\pi n)\\
&=\frac{1}{\sqrt{2B(\alpha)}}\sum_{n=0}^{2B(\alpha)-1}e^{-\frac{i\pi j n}{B(\alpha)}}\frac{1}{\sqrt{2B(\alpha)}}\sum _{m=0}^{2B(\alpha)-1} e^{\frac{i\pi m n}{B(\alpha)}} g_m(x)\\
&=\frac{1}{2B(\alpha)}\sum_{m=0}^{2B(\alpha)-1}g_m(x)\sum_{n=0}^{2B(\alpha)-1}e^{\frac{i\pi n(m-j)}{B(\alpha)}}\\
&=g_j(x),
\end{align*}
where the last line follows from the fact that the sum of all the $B(\alpha)^{th}$ roots of unity is zero; in particular, for $m\neq j$ we have $\sum_{n=0}^{2B(\alpha)-1}e^{\frac{i\pi n(m-j)}{B(\alpha)}}=0,$ while for $m=j$ we get $\sum_{n=0}^{2B(\alpha)-1}e^{\frac{i\pi n(m-j)}{B(\alpha)}}=2B(\alpha)$.
Next, we show that $U$ is an isometry, and thus it can be extended to a unitary operator from $L^2(\mathbb{T}_\alpha)$ to $\bigoplus_{j=0}^{2B(\alpha)-1} L^2(0,\pi).$ So,

\begin{align*}
\norm{Uh}_2^2&=\sum_{j=0}^{2B(\alpha)-1}\int_0^{\pi}\left|(Uh)_j(x)\right|^2dx\\
&=\sum_{j=0}^{2B(\alpha)-1}\int_0^\pi\left|\frac{1}{\sqrt{2B(\alpha)}}\sum_{n=0}^{2B(\alpha)-1}e^{-\frac{i\pi j n}{B(\alpha)}}h(x+\pi n)\right|^2dx\\
&=\int_0^\pi\sum_{j=0}^{2B(\alpha)-1}\left|\sum_{n=0}^{2B(\alpha)-1}\frac{1}{\sqrt{2B(\alpha)}}e^{-\frac{i\pi j n}{B(\alpha)}}h(x+\pi n)\right|^2dx\\
&=\int_0^\pi\sum_{j=0}^{2B(\alpha)-1}\left|\Big\langle{\frac{1}{\sqrt{2B(\alpha)}}e^{\frac{i\pi j (\cdot)}{B(\alpha)}},h(x+\pi(\cdot))\Big\rangle}_{\ell^2(\{0,\dots,2B(\alpha)-1\})}\right|^2dx\\
&=\int_0^\pi\norm{h(x,\cdot)}_{\ell^2(\{0,\dots,2B(\alpha)-1\})}^2dx\\
&=\int_0^\pi\sum_{j=0}^{2B(\alpha)-1}\left|h(x+\pi j)\right|^2dx\\
&=\int_0^{2\pi B(\alpha)}|h(x)|^2dx\\
&=\norm{h}_2^2.
\end{align*}

Above, we have used the fact that $\Big\{\frac{1}{2B(\alpha)}e^{i\pi j(\cdot)}\Big\}_{j=0}^{2B(\alpha)-1}$ is an orthonormal basis of $\ell^2(\{0,\dots, 2B(\alpha)-1\})$, and Parseval's identity going from line four to five.
\item We have
\begin{align*}
\left(\bigoplus_{l=0}^{2B(\alpha)-1}L_l^kUh\right)_j(y)&=\int_0^\pi L_{E,\alpha, j}^k(y,x)(Uh)_j(x)dx\\
&=\int_0^\pi \sum_{n=0}^{2B(\alpha)-1}K_{E,\alpha}^{k}(y,x+n\pi)e^{\frac{i\pi j n}{B(\alpha)}}\frac{1}{\sqrt{2B(\alpha)}}\sum_{m=0}^{2B(\alpha)-1}e^{-\frac{i\pi j m}{B(\alpha)}}h(x+\pi m)dx\\
&= \frac{1}{\sqrt{2B(\alpha)}}\sum_{n=0}^{2B(\alpha)-1}\sum_{m=0}^{2B(\alpha)-1}e^{-\frac{i\pi j(m- n)}{B(\alpha)}}\int_0^\pi K_{E,\alpha}^{k}(y,x+n\pi)h(x+\pi m)dx.\\
\end{align*}
On the other hand,
\begin{align*}
(U\tilde T_{E,\alpha}^{k}h)_j(y)&=\frac{1}{\sqrt{2B(\alpha)}}\sum_{n=0}^{2B(\alpha)-1}e^{-\frac{i\pi jn}{B(\alpha)}}(\tilde T_{E,\alpha}^{k}h)(y+n\pi)\\
&=\frac{1}{\sqrt{2B(\alpha)}}\sum_{n=0}^{2B(\alpha)-1}e^{-\frac{i\pi jn}{B(\alpha)}}\int_0^{2\pi B(\alpha)}K_{E,\alpha}^{k}(y+n\pi,x)h(x)dx\\
&=\frac{1}{\sqrt{2B(\alpha)}}\sum_{n=0}^{2B(\alpha)-1}e^{-\frac{i\pi jn}{B(\alpha)}}\sum_{m=0}^{2B(\alpha)-1}\int_0^{\pi }K_{E,\alpha}^{k}(y+n\pi,x+m\pi)h(x+m\pi)dx\\
&=\frac{1}{\sqrt{2B(\alpha)}}\sum_{n=0}^{2B(\alpha)-1}\sum_{m=0}^{2B(\alpha)-1}e^{-\frac{i\pi jn}{B(\alpha)}}\int_0^{\pi }K_{E,\alpha}^{k}(y,x+(m-n)\pi)h(x+m\pi)dx\\
&=\frac{1}{\sqrt{2B(\alpha)}}\sum_{\tilde n=0}^{2B(\alpha)-1}\sum_{m=0}^{2B(\alpha)-1}e^{-\frac{i\pi j(m-\tilde n)}{B(\alpha)}}\int_0^{\pi }K_{E,\alpha}^{k}(y,x+\tilde n\pi)h(x+m\pi)dx.\\
\end{align*} The last equality above follows by using the fact that  $K_{E,\alpha}^{k}$ is defined up to $mod\, 2\pi B(\alpha)$. The claim follows by comparing these two expressions.

\item  From above it follows that $\norm{\tilde T_{E,\alpha}^{k}}_{2,2}=\max_{0\leq j\leq 2B(\alpha)-1}\norm{L_j^k}_{2,2}.$ Now, since $K_{E,\alpha}^{k}(x,y)\geq 0$, it follows readily that $\left|L_{E,\alpha, j}^k\right|\leq L_{E,\alpha, 0}^k$ and thus $\norm{L_j^k}_{2,2}\leq \norm{L_0^k}_{2,2},$ for all $j$, thus proving the claim. 
\end{enumerate}
\end{proof}

\begin{proof}[Proof of Proposition $\ref{prop5.1}$] Here, we adapt the argument for the analogous continuum result in \cite{ds}. By Lemmas $\ref{lemma4.5}$ and $\ref{lemma5.2}$ it suffices to show that $\norm{L_0^k}_{2,2}\neq 1.$ We establish this via a proof by contradiction; that is, assume that $\norm{L_0^k}_{2,2}=1.$ Since $L_0^k$ is a compact operator (this follows from Lemma $\ref{lemma4.5}$ and Lemma $\ref{lemma5.2}$ part $(b)$), it implies that $|L_0^k|$ is also compact, hence it follows that $\norm{|L_0^k|}_{2,2}$ is an eigenvalue of $|L_0^k|.$ Then, using the fact that $\norm{L_0^k}=\norm{|L_0^k|}$, and $\norm{|L_0^k|f}_{2}=\norm{L_0^kf}_{2}$ it follows that there exists some $f\neq 0$ such that $\norm{L_0^kf}_2=\norm{f}_2$, where $f$ is chosen to be the eigenvector corresponding to the eigenvalue $\norm{L_0^k}_{2,2}=\norm{|L_0^k|}_{2,2}=1$. Let $\tilde f$ be the $\pi-$periodic extension of $f$ to $\mathbb{T}_\alpha$; that is $\tilde f(x+n\pi)=f(x)$, for all $x\in(0,\pi).$ Then,
\begin{align*}
\left(L_0^kf\right)(y)&=\int_0^\pi \sum_{n=0}^{2B(\alpha)-1}K_{E,\alpha}^{k}(y,x+n\pi) f(x)dx\\
&=\sum_{n=0}^{2B(\alpha)-1}\int_0^\pi K_{E,\alpha}^{k}(y,x+n\pi)\tilde f(x+n\pi)dx\\
&=\int_{\mathbb{T}_\alpha}K_{E,\alpha}^{k}(y,x)\tilde f(x)dx.
\end{align*}
Then,
\begin{align*}
\norm{L_0^kf}_2^2&=\int_0^\pi|\left(L_0^kf\right)(y)|^2dy\\
&=\int_0^\pi\left|\int_{\mathbb{T}_\alpha}K_{E,\alpha}^{k}(y,x)\tilde f(x)dx\right|^2dy\\
&=\int_0^\pi\left|\int_{\mathbb{T}_\alpha}\sqrt{K_1^{k}(y,x)}\sqrt{K_2^{k}(y,x)}\tilde f(x)dx\right|^2dy\\
&\leq \int_0^\pi\left(\int_{\mathbb{T}_\alpha}K_1^{k}(y,x)dx\int_{\mathbb{T}_\alpha}K_2^{k}(y,x)|\tilde f(x)|^2dx\right)dy\\
&=\int_0^\pi\int_{\mathbb{T}_\alpha}K_2^{k}(y,x)|\tilde f(x)|^2dxdy\\
&=\int_0^\pi \sum_{n=0}^{2B(\alpha)-1}\int_{-n\pi}^{-(n-1)\pi}K_2^{k}(y,x)|\tilde  f(x)|^2dx dy\\
&=\int_0^\pi \sum_{n=0}^{2B(\alpha)-1}\int_{0}^{\pi}K_2^{k}(y,x-n\pi)|\tilde  f(x-n\pi)|^2dx dy\\
&=\int_0^\pi \sum_{n=0}^{2B(\alpha)-1}\int_{0}^{\pi}K_2^{k}(y+n\pi,x)|f(x)|^2dx dy\\
&=\int_0^\pi\left(\sum_{n=0}^{2B(\alpha)-1}\int_0^\pi K_2^{k}(y+n\pi,x)dy\right)|f(x)|^2dx \\
&=\int_0^\pi \left(\int_{\mathbb{T}_\alpha} K_2^{k}(y,x)dy\right)|f(x)|^2dx\\
&=\int_0^\pi|f(x)|^2dx\\
&=\norm{f}_2^2\\
&=\norm{L_0^kf}_2,
\end{align*}
where above, among other facts, we have used results appearing in the proof of Lemma $\ref{lemma4.6}$. Since, we have equality in the application of Cauchy-Schwarz inequality, we conclude that for almost every $y\in (0,\pi)$ the functions $\Big|\sqrt{K_1^{k}(y,\cdot)}\Big|^2$ and $\Big|\sqrt{K_2^{k}(y,\cdot)}\tilde f(\cdot)\Big|^2$ are linearly dependent almost everywhere. That is, for $y\in (0,\pi)\setminus N$, with $Leb(N)=0$, there is some constant $C_y>0$ such that $$C_y K_1^{k}(y,\cdot)=K_2^{k}(y,\cdot)\tilde f(\cdot)^2\, a.e.$$ For each $y\in [0,\pi)\setminus N$, set $$M_y\defeq\{x: \lambda(y,x)\in \supp r\}.$$ Then, for almost every $x\in M_y$ we have $C_y=R_+^2(\alpha k-1,x,\lambda(y,x,E), E)\tilde f(x)^2.$ Rewriting this, using an immediate consequence of Lemma $\ref{lemmm3.2}$, we get \begin{equation}\label{eq46} \tilde f(x)^2=C_yR_-^2(\alpha k+\alpha -1,y, \lambda(x,y,E),E).\end{equation}

Our next goal is to show that $\tilde f^2$ is real analytic on $\mathbb{R}$. To this end, let $[a,b]$ be a non-trivial interval contained in $\supp r$. The existence of such an interval is guaranteed by our assumption on $r$. Let $c_y$ and $d_y$ be the unique phases determined by $\lambda(y,c_y,E)=a$ and $\lambda(y,d_y,E)=b$; that is, $\phi_-(\alpha k+\alpha-1,y,a,E)=c_y$ and $\phi_-(\alpha k+\alpha-1,y,b,E)=d_y$. Then, from Lemma $\ref{lemm1}$ it follows that $c_y<d_y$. By construction, $\phi_-(\alpha k+\alpha-1,y,\cdot, E)$ is continuous when viewed as a function from $\supp r$ to $\mathbb{T}_\alpha$. Hence, by the intermediate value theorem, it follows that $[c_y,d_y]\subset M_y.$ Moreover, again by construction, it follows that $c_y$ and $d_y$ are continuous functions of $y$, and by Lemma $\ref{lemm3.2'}$, that they are strictly increasing. From $\phi_-(\alpha k+\alpha-1,y+\pi,\cdot,E)=\phi_-(\alpha k+\alpha-1,y,\cdot,E)+\pi$, we have $[c_{y+\pi},d_{y+\pi}]=[c_y+\pi,d_y+\pi].$ Set $$I\defeq \bigcup_{y\in[0,\pi)\setminus N}(c_y,d_y).$$ From the discussion above it follows that $I$ is an open interval of length greater than $\pi$.\\

Next, since the solution $u_-$ of the difference equation $(\ref{eq024})$ is simply a polynomial in $\lambda$, it follows that in particular it is analytic in $\lambda$. Consequently, $R_-^2(\alpha k+\alpha -1, y,\lambda, E)$, as the sum of the squares of two analytic functions, is also analytic in $\lambda$. From our discussion above and equation $(\ref{eq45})$ it follows that in turn $x(\lambda)$ is analytic in $\lambda$ and by Lemma $\ref{lemm1}$ that $x'(\lambda)>0$, thus implying that the inverse function of $x(\lambda)$, namely $\lambda(y,x, E)$ is analytic in $x$. 

The discussion above and the expression in $(\ref{eq46})$, shows that, in particular, for each fixed $y\in (0,\pi)\setminus N$ the function $\tilde f^2$ is real analytic  on $(c_y,d_y)$, and as a result on the entire $I$ as well. Moreover, since $\tilde f$ is $\pi-$periodic, we conclude that $\tilde f^2$ is analytic on the entire real line. Now, we fix again $y\in [0,\pi)\setminus N$, and by analytic continuation conclude that in fact $(\ref{eq46})$ holds for all $x$ for which $\lambda(x,y,E)$ exists. Next, since all  real analytic periodic functions are bounded on $\mathbb{R}$, and since $\tilde f^2$ is a $\pi-$periodic analytic function on $\mathbb{R}$, from $(\ref{eq46})$ we conclude that $R_-^2(\alpha k+\alpha-1,y,\lambda(x,y,E),E)$ is bounded in $x$. 
 But since $\lambda(x,y,E)$ takes on arbitrary values as $x$ varies, we get \begin{equation}\sup_{\lambda \in\mathbb{R} }R_-(\alpha k+\alpha-1,y,\lambda,E)<\infty,\end{equation} contradicting Lemma $\ref{lemma5.3}$ below. This contradiction shows that $\norm{\tilde T_{E,\alpha}^{k}}_{2,2}=\norm{L_0^k}_{2,2}<1$, and thus proving the claim!

\end{proof}

\begin{lemm}\label{lemma5.3}
It holds that $$\lim_{\lambda\to\infty} R_{-}(\alpha k+\alpha-1, y, \lambda,E)=\infty.$$
\end{lemm}

\begin{proof}
We have

\begin{align*}
R_-(\alpha k+\alpha-1, y,\lambda, E)&=\prod_{i=\alpha-1}^0\Big|\frac{\cos\phi_-(\alpha k+i-1, y,\lambda, E)}{\sin\phi_-(\alpha k+i,y,\lambda, E)}\Big|\\
&=\big|\cos\phi_-(\alpha k-1,y,\lambda, E)\cot \phi_-(\alpha k,y,\lambda, E)\dots\cot\phi_-(\alpha k+\alpha -2,y,\lambda, E)\big|\\
&\cdot\Big|\frac{1}{\sin \phi_-(\alpha k+\alpha -1,y,\lambda, E)}\Big|\\
&=\left|\cos y(\lambda)\Big(E-f_0\lambda-\tan y(\lambda)\Big)\left(E-f_1\lambda-\frac{1}{E-f_0\lambda-\tan y(\lambda)}\right)\right|\dots\\
&\dots\left|\left(E-f_{\alpha-2}\lambda-\frac{1}{E-f_{\alpha-3}\lambda-\frac{1}{\dots+\frac{1}{E-f_0\lambda-\tan y(\lambda)}}}\right)\right|\cdot\\
&\sqrt{1+\left(E-f_{\alpha-1}\lambda-\frac{1}{E-f_{\alpha-2}\lambda-\frac{1}{E-f_{\alpha-3}\lambda-\frac{1}{\dots+\frac{1}{E-f_0\lambda-\tan y(\lambda)}}}}\right)^2}.
\end{align*}
The first equality follows from $(\ref{eq030})$, and to go from the second line to the third we have used the analogue of equation $(\ref{eq3})$ for $\phi_-$ iteratively, and the trig identity $\left|\left(\sin A\right)^{-1}\right|=\sqrt{1+\left(\cot A\right)^2}.$ We will consider in more detail the case when $\alpha=3$, then the same argument can be extended to the general case. Before we do so, let us make a crucial observation.  Since the coupling constant $\lambda$ defines the $\alpha-$block starting at site $\alpha k$ and ending at site $\alpha k+\alpha-1$, and since the Pr\"ufer angle $y$ (i.e. in previous  notation this is the angle we denote by $\theta_{k-1}$, see the paragraph where equation $(\ref{eq024})$ is discussed!) is determined by the solutions at sites $\alpha k-1$ and $\alpha k$, when solving from left to right, as such it is fixed and thus independent of the coupling constant $\lambda$. For $\alpha=3$ we have

\begin{align*}
R_-(3k+2, y,\lambda, E)&=\cos y(\lambda)(E-f_0\lambda -\tan y(\lambda))\left(E-f_1\lambda-\frac{1}{E-f_0\lambda-\tan y(\lambda)}\right)\\
&\cdot\sqrt{1+\left(E-f_2\lambda-\frac{1}{E-f_1\lambda-\frac{1}{E-f_0\lambda-\tan y(\lambda)}}\right)^2}\\
&=\cos y(\lambda)(E-f_0\lambda -\tan y(\lambda))\sqrtx{\left(E-f_1\lambda-\frac{1}{E-f_0\lambda-\tan y(\lambda)}\right)^2}{+\left((E-f_2\lambda)\left(E-f_1\lambda-\frac{1}{E-f_0\lambda-\tan y(\lambda)}\right)-1\right)^2}\\
&=\cos y(\lambda)\sqrtx{\Big[(E-f_1\lambda)\big(E-f_0\lambda -\tan y(\lambda)\big)-1\Big]^2}{+\Big[(E-f_2\lambda)\Big((E-f_1\lambda)\big(E-f_0\lambda -\tan y(\lambda)\big)-1\Big)-\big(E-f_0\lambda -\tan y(\lambda)\big)\Big]^2}\\
&=\sqrtx{\Big[(E-f_1\lambda)\big(E-f_0\lambda)\cos y(\lambda)-(E-f_1\lambda)\sin y(\lambda)-\cos y(\lambda)\Big]^2}{+\Big[(E-f_2\lambda)\Big((E-f_1\lambda)(E-f_0\lambda)-(E-f_0\lambda)-1\Big)\cos y(\lambda)}\\
&\overline{-\Big((E-f_2\lambda)(E-f_1\lambda)+1\Big)\sin y(\lambda)\Big]^2}.\\
\end{align*}
Now, since $y(\lambda)$ is fixed and independent of $\lambda$, we note that, regardless of its value, the expression under the square root above will either be a polynomial of degree three or of degree two in $\lambda$ (i.e. depending on whether $\sin y(\lambda)$ or $\cos y(\lambda)$ is zero). Thus, it is clear that $R_-(3k+2,y,\lambda, E)\to \infty$ as $\lambda\to\infty$, proving the claim!
\end{proof}


\section{ The dependence of the $\norm{\cdot}_{2,2}-$norm of $\tilde T_{E,\alpha}^{k}$ on $k$ and the energy $E$}\label{sec7}

The goal of this section is to show that $\tilde T_{E,\alpha}^{k}$ depends continuously on $E$. More precisely, we prove the following theorem.

\begin{thm}\label{thm6.1}
The real valued map $E\mapsto \norm{\tilde T_{E,\alpha}^{k}}_{2,2}$ is continuous on $\Sigma_0$.
\end{thm}

Before we begin the proof of the above theorem, let us first observe that the norm of $\tilde T_{E,\alpha}^{k}$ is independent of $k$. The way to see this, is by noting that all the quantities involved in the definition of $\tilde T_{E,\alpha}^{k}$, specifically the Pr\'ufer amplitude $R_+$ and the solution $u_+$ of the difference equation $(\ref{eq024})$, are defined locally on each $\alpha-$block, and the fact that the random variables are identically distributed.\\

Next, we begin the proof of Theorem $\ref{thm6.1}$, by first proving two preparatory lemmas by adapting the analogous arguments for the continuum case. 

\begin{lemm}\label{lemma6.2}
Suppose $E_n\to E$. Then, $$\lambda(y,x,E_n)\to \lambda(y,x,E),$$ whenever $\lambda(y,x,E)$ exists.
\end{lemm}
\begin{proof}
Suppose $\lambda(y,x,E)$ exists, that is, $\lambda_0=\lambda(y,x,E)$ for some $\lambda_0\in\mathbb{R}.$ By definition, this means that $\phi_+(\alpha k-1,x,\lambda_0,E)=y.$ Now, since $\phi_+$ is strictly increasing in the coupling constant $\lambda$, it follows that given any $\epsilon>0$ we have $$\phi_+(\alpha k-1,x,\lambda_0-\epsilon,E)<y<\phi_+(\alpha k-1,x,\lambda_0+\epsilon, E).$$ Then, since by construction, $\phi_+$ depends continuously on the energy $E$, it follows that there exists some $n_0$ such that for all $n\geq n_0$ we have 
$$\phi_+(\alpha k-1,x,\lambda_0-\epsilon, E_n)<y<\phi_+(\alpha k-1, x,\lambda_0+\epsilon, E_n).$$ So, by the intermediate value theorem, there is some $\lambda_n\in(\lambda_0-\epsilon,\lambda_0+\epsilon)$ such that $\phi_+(\alpha k-1,x,\lambda_n, E_n)=y$; that is, $\lambda(y,x,E_n)=\lambda_n$. In conclusion, we have shown that, given any $\epsilon>0$, there is some $n_0$ such that for all $n\geq n_0$ we have $\lambda(y,x,E_n)\in(\lambda(y,x,E)-\epsilon,\lambda(y,x,E)+\epsilon)$, which shows that $\lambda(y,x,E_n)\to\lambda(y,x,E)$ as $n\to\infty$, thus proving the claim.
\end{proof}

\begin{lemm}\label{lemma6.3}
Suppose $E_n\to E.$ Then, we have $$\mathbb{R}^2\setminus A(E)\subset\liminf_{n\to\infty}\mathbb{R}^2\setminus A(E_n),$$ where $A(E)=\lambda(\cdot,\cdot,E)^{-1}\left([-M,M]\right)$ and $\supp r\subset[-M,M].$
\end{lemm}

\begin{proof}
Let $(y,x)\in \mathbb{R}^2\setminus A(E);$ that is, $\lambda(y,x,E)$ either does not exist, or it exists but lies outside of $[-M,M].$ Suppose there is a subsequence $n_j\to \infty$ such that $\lambda(y,x,E_{n_j})$ exists and belongs to some interval $[a,b]\subset[-M,M]$ for every $j$. Then, $\phi_-(\alpha k+\alpha-1,y, \lambda(y,x,E_{n_j}), E_{n_j})=x.$ By monotonicity, we have $$\phi_-(\alpha k+\alpha -1, y, a, E_{n_j})\leq x\leq \phi_-(\alpha k+\alpha-1,y,b,E_{n_j}).$$ Since $\phi_-$ is continuous as a function of the energy $E$, taking $j\to\infty$, we find that $$\phi_-(\alpha k+\alpha-1,y,a, E)\leq x\leq \phi_-(\alpha k+\alpha-1,y,b,E).$$ By the intermediate value theorem, this  would imply that there exist some $\lambda_0\in[-M,M]$ such that $\phi_-(\alpha k+\alpha-1,y,\lambda_0,E)=x$, which is the same as saying that $\lambda(y,x,E)$ exists and is equal to $\lambda_0$, thus leading us to a contradiction.
\end{proof}

\begin{lemm}\label{lemma6.4}
Suppose $E_n\to E$. Then $$\lim_{n\to\infty} K_{E_n,\alpha}^{k}(y,x)=K_{E,\alpha}^{k}(y,x)$$ for almost every $(y,x)\in \mathbb{R}^2.$
\end{lemm}

\begin{proof}
The case where $(y,x)\in \mathbb{R}^2\setminus A(E)$ follows trivially from Lemma $\ref{lemma6.3}.$
Next, we consider the case where $\lambda(y,x,E)$ exists. First, let $B$ denote the set of discontinuities of $r$, and let $K\defeq \lambda_E^{-1}(B)\bigcup \lambda_{E_n}^{-1}(B)$. Since, $Leb(B)=0$, then by Proposition $\ref{Leb}$ it follows that $Leb(K)=0.$ Then by Lemma $\ref{lemma6.2}$ we know that $\lambda(y,x,E_n)\to \lambda(y,x,E)$ as $n\to\infty$. So, using the fact that $u_+$ and $R_+$ depend continuously on $E$, and that $r$ is continuous away from $B$, for every $(x,y)\in\mathbb{R}^2\setminus K$, we get that
$$K_{E_n,\alpha}^{k}(y,x)= \frac{R_{+}(\alpha k -1,x,\lambda(y,x,E_n),E_n)}{\sum_{i=0}^{\alpha-1}f_iu_{+}^2(\alpha k+i, x,\lambda(x,y,E_n),E_n)}r(\lambda(y,x,E_n))\to K_{E,\alpha}^{k}(y,x),\, \text{as }\, n\to \infty,$$ concluding the proof.
\end{proof}

Now, we are ready to give the proof of Theorem $\ref{thm6.1}$.

\begin{proof}[Proof of Theorem $\ref{thm6.1}$]
 Let $E$ be any point in $\Sigma_0$, and let ${E_n}\in \Sigma_0$ be an arbitrary sequence converging to $E$; that is, $|E_n-E|\to0$ as $n\to\infty.$ We will show that $$\lim_{n\to\infty}\norm{\tilde T_{E_n,\alpha}^{k}}_{2,2}=\norm{\tilde T_{E,\alpha}^{k}}_{2,2}.$$ From Lemma $\ref{lemma5.2}$, since $\norm{\tilde T_{E,\alpha}^{k}}_{2,2}=\norm{L_{0,E}^k}_{2,2}$, it suffices to show that $\norm{L_{0,E_n}^k-L_{0,E}^k}_{2,2}\to 0$, as $n\to\infty$, where $$\big(L_{0,(\cdot)}^{k}f\big)(x)=\int_0^\pi L_{(\cdot), \alpha, 0}^k(x,y)f(y)dy,\, \text{ and }\, L_{(\cdot),\alpha,0}^k=\sum_{n=0}^{2B(\alpha)-1}K_{(\cdot), \alpha}^{k}(x,y+n\pi).$$
From the definition of $K_{E,\alpha}^{k}$ in Definition $\ref{defi5.4}$, and Lemma $\ref{lemm3.10}$, it follows that $L_{(\cdot), \alpha, 0}^k(x,y)$ is uniformly bounded on  on $\{E\}\cup\{E_n\}$. On the other hand, Lemma $\ref{lemma6.4}$ shows that for almost every $(x,y)\in\mathbb{R}^2$ we have $L_{E_n,\alpha, 0}^k(\cdot,\cdot)\to L_{E,\alpha,0}^k(\cdot,\cdot)$, as $n\to \infty$. Thus, using dominated convergence theorem, we get
\begin{align*}
\norm{L_{0,E_n}^k-L_{0,E}^k}_{2,2}&\leq\norm{L_{0,E_n}^k-L_{0,E}^k}_{HS}\\
&=\int_0^\pi\int_0^\pi\left|L_{E_n,\alpha,0}^k(x,y)-L_{E,\alpha,0}^k(x,y)\right|^2dydx\to 0,\hspace{.5cm}\text{as}\hspace{.5cm}n\to\infty,\\
\end{align*}
concluding the proof!
\end{proof}

\begin{cor}\label{cor6.5} There exists some constant $0<q<1$, such that $$\sup_{E\in\Sigma_0}\norm{\tilde T_{E,\alpha}^{k}}_{2,2}\leq q<1.$$
\end{cor}
\begin{proof}
Since $\Sigma_0$ is compact, from Theorem $\ref{thm6.1}$ it follows that $\norm{\tilde T_{E,\alpha}^{k}}_{2,2}$ attains its maximum in $\Sigma_0.$ That is, there is some $E_0\in\Sigma_0$, such that $\norm{\tilde T_{E,\alpha}^{k}}_{2,2}\leq \norm{\tilde T_{E_0,\alpha}^{k}}_{2,2}$, for all $E\in\Sigma_0$, hence we also have $$\sup_{E\in\Sigma_0}\norm{\tilde T_{E,\alpha}^{k}}_{2,2}\leq \norm{\tilde T_{E_0,\alpha}^{k}}_{2,2}.$$ From Proposition $\ref{prop5.1}$ it follows that there is some $0<q<1$ such that $\norm{\tilde T_{E_0,\alpha}^{k}}_{2,2}\leq q$, proving the claim.
\end{proof}


\section{Proof of Theorem $\ref{thm2.1}$} 

Now, we are in a position to give the proof of Theorem $\ref{thm2.1}$, which will merely be a collection of the facts proven above. As we have already remarked, it is sufficient to show that $\rho_L(m,0)\leq Ce^{-\gamma \left|\lfloor{\frac{m}{\alpha}\rfloor}\right|}.$ First, with no loss of generality, suppose that $m=\alpha k_0+j$, for some fixed $k_0$ and some $j\in\{0,1,\dots,\alpha-1\}$, then we have

\[\hspace{-10cm}\int_\Omega\left(\sup_{t\in\mathbb{R}}\left|\langle{\delta_m,e^{-itH_\omega}\delta_0\rangle}\right|\right)d\mu(\omega)\]
\begin{align*}&=a(m,0)\\
&\leq \liminf_{L\to\infty} a_L(m,0)\\
&\leq\liminf_{L\to\infty} \rho_L(m,0)\\
&\leq \liminf_{L\to\infty}\int_{\Sigma_0}\rho_L(m,0,E)dE\\
&\leq \tilde C \liminf_{L\to\infty}\int_{\Sigma_0}\sum_{N=0}^{2B(\alpha)-1}\Big\langle{\tilde T_{E,\alpha}^{1}\dots \tilde T_{E,\alpha}^{k_0}T_{E,\alpha}^{k_o+1}\dots T_{E,\alpha}^{L-2}\psi_{E,\alpha,N}^L,T_{E,\alpha}^{0}T_{E,\alpha}^{-1}\dots T_{E,\alpha}^{-L+1}\psi_{E,\alpha}^{-L}\Big\rangle}_{L^2(\mathbb{T},d\theta_0)}d E.\\
&\leq \tilde C\liminf_{L\to\infty}\sum_{N=0}^{2B(\alpha)-1}\int_{\Sigma_0}\norm{\tilde T_{E,\alpha}^{1}\dots \tilde T_{E,\alpha}^{k_0}T_{E,\alpha}^{k_o+1}\dots T_{E,\alpha}^{L-2}\psi_{E,\alpha}^L}_2\norm{T_{E,\alpha}^{0}T_{E,\alpha}^{-1}\dots T_{E,\alpha}^{-L+1}\psi_{E,\alpha}^{-L}}_2 d E\\
&\leq \tilde C\liminf_{L\to\infty}\sum_{N=0}^{2B(\alpha)-1}\int_{\Sigma_0}\prod_{i=1}^{k_0}\norm{\tilde T_{E,\alpha}^i}_{2,2}\norm{T_{E,\alpha}^{k_0+1}}_{1,2}\prod_{i=k_0+2}^{L-2}\norm{T_{E,\alpha}^i}_{1,1}\norm{\psi_{E,\alpha}^L}_1\norm{T_{E,\alpha}^{0}}_{1,2}\prod_{i=-1}^{-L+1}\norm{T_{E,\alpha}^i}_{1,1}\norm{\psi_{E,\alpha}^{-L}}_1 d E\\
&\leq \tilde C\liminf_{L\to\infty}\sum_{N=0}^{2B(\alpha)-1}\int_{\Sigma_0} q^{k_0}\cdot C_0\cdot C_1\cdot C_0\cdot C_2 d E\\
&=\bar C 2B(\alpha) Leb(\Sigma_0)q^{k_0}\\
&=Ce^{-\gamma \lfloor{\frac{m}{\alpha}\rfloor}},
\end{align*} where $\gamma=\log\left(q^{-1}\right)>0$, since $q<1.$

The first equality follows by definition; in the second and third lines we have used Lemmas $\ref{lemma02.3}$ and $\ref{lemma02.4}$; going from line three to line four we have used Lemma \ref{lem04.7} ; from line four to five we have used $\ref{lemma4.9}$; from line five to six the Cauchy-Schwarz inequality; from line five to six standard results for operator norms; finally, in the last inequality we have used Lemmas $\ref{lemma05.2}$, $\ref{lemm4.3.}$, and Corollary $\ref{cor6.5}$.


\section{The operator $T_1$ is a strict contraction}

The purpose of this section is to deduce that the single operator $T_1$, defined in the original Kunz-Souillard setting, is a strict contraction. We remind the reader that in the original Kunz-Souillard work, it is originally only known for the second iterate of $T_1$ to be a contraction, hence, from this point of view, the result we present here is an improvement of the original work. 
We achieve this via appealing to the norm estimates in the previous section, and by showing that, in the special case for $\alpha=1$, and $f_0=1$,  the integral operators $\tilde T_{E,\alpha}^{k}\equiv \tilde T_E^{k}$, are unitarily equivalent to the $T_1$ operator defined in \cite{dks}, \cite{ks}.\\

For each $k\geq 0$ we consider the following operator on $L^2(\mathbb{T}_1),$ where $\mathbb{T}_1=\mathbb{R}/(2\pi B(1)\mathbb{Z})$, \begin{equation}\label{eq1} \left(\tilde T_E^{k}f\right)(x)=\int_{\mathbb{T}_1}\frac{R_{+}(k-1,y,\lambda(x,y,E),E)}{u_{+}^2(k,y,\lambda(x,y,E),E)}r(\lambda(x,y,E))f(y)dy.\end{equation}

Before we rewrite the above integral, let us note that it follows immediately from the proof of Lemma $\ref{lemma03.1}$ that in the case $\alpha=1$ we can actually take $B(1)$ to be $1$. Next, let us rewrite this integral operator using the definition of $R_+$ and $u_+$. That is, since $u_{+}(k,y,\lambda,E)=\sin y$ and $$R_+(k,y,\lambda, E)\sin\phi_+(k,y,\lambda,E)=R_+(k-1,y,\lambda,E)\cos \phi_+(k-1,y,\lambda,E)$$ with $R_+(k,y,\lambda,E)=1$, $\phi_+(k,y,\lambda,E)=y$ and $\phi_+(k-1,y,\lambda, E)=x$ we get $$\frac{ R_+(k-1,y,\lambda,E)}{u^2_+(k,y,\lambda,E)}=\frac{1}{\cos x\sin y}.$$ 
Next, we claim that $\lambda(x,y,E)=E-\tan x-\cot y$. We show that, given any $x,y\in \mathbb{T}$, there is a coupling constant $\lambda=E-\tan x-\cot y$, such that $\phi_+(k-1,y,\lambda, E)=x.$

From $$\cot \phi_+(k,y,\lambda, E)+\tan \phi_+(k-1,y,\lambda, E)=E-\lambda$$ we get $\phi_+(k-1,y,\lambda,E)=\tan^{-1}(E-\lambda-\cot y).$
So, given any $x,y\in\mathbb{T}_1$, taking $\lambda$ to be $E-\tan x-\cot y$, we clearly get $\phi_+(k-1,y,\lambda,E)=x.$ Whenever such a coupling constant exists, we set $\lambda(x,y,E)=\lambda=E-\tan x-\cot y$, as claimed. 

Now, we rewrite $(\ref{eq1})$: 
\begin{equation}
\left(\tilde T_E^{k}f\right)(x)=\frac{1}{\cos x}\int_{\mathbb{T}_1}r\left(E-\tan x-\frac{1}{\tan y}\right)\frac{1}{\sin y}f(y)dy.\end{equation}

Next, let us define an operator $T_{1,E}$ on $L^2(\mathbb{R})$, as follows
\begin{equation}
\left(T_{1,E}f\right)(u)=\int_{\mathbb{R}}r(E-u-v^{-1})|v|^{-1}f(v)dv.\end{equation}

We will show that these two operators are conjugates of one another!\\

Define $\mathcal{U}_1:L^{2}(\mathbb{R})\to L^2(\mathbb{T}_1)$ by $$\left(\mathcal{U}_1f\right)(x)=\sec xf(\tan x)$$ and $\mathcal{U}_2:L^2(\mathbb{T}_1)\to L^2(\mathbb{R})$ by $$\left(\mathcal{U}_2f\right)(x)=\frac{1}{\sqrt{1+x^2}}f(\tan^{-1}x),$$ and note that they are both unitary operators and inverses of one another. 

\begin{prop}\label{prop8.1}
With the same notation as above, we have $\tilde T_E^{k}=\mathcal{U}_1 T_{1,E}\mathcal{U}_2.$
\end{prop}

\begin{proof}
\begin{align*}
\left(\mathcal{U}_1T_{1,E}\mathcal{U}_2f\right)(x)&=\sec x\left(S_E\mathcal{U}_2f\right)(\tan x)\\
&=\sec x\int_{\mathbb{R}} r\left(E-\tan x-v^{-1}\right)|v|^{-1}\left(\mathcal{U}_2f\right)(v)dv\\
&=\sec x\int_{\mathbb{R}}r\left(E-\tan x-v^{-1}\right)|v|^{-1}\frac{1}{\sqrt{1+v^2}}f\left(\tan^{-1}v\right)dv.\\
&=\frac{1}{\cos x}\int_{\mathbb{T}_1}r\left(E-\tan x-\frac{1}{\tan y}\right)\frac{1}{\tan y}\frac{1}{\sec y}f(y)\sec^2 y\,dy\\
&=\frac{1}{\cos x}\int_{\mathbb{T}_1}r\left(E-\tan x-\frac{1}{\tan y}\right)\frac{1}{\sin y}f(y)dy\\
&=\left(\tilde T_E^{k}f\right)(x).
\end{align*}
Above, in the third line, we have performed the following change of variables $y=\tan^{-1}v$. Then, $dv=\sec^2y\,dy$.

\end{proof}

\begin{thm}There is some $0<q<1$ such that $$\sup_{E\in\Sigma_0}\norm{T_{1,E}}_{2,2}\leq q<1.$$

\end{thm}
\begin{proof}
In Proposition $\ref{prop8.1}$ we established that $T_{1,E}$ and $\tilde T_E^{k}$ are unitary equivalent operators. 
So, it follows that $$\norm{T_{1,E}}_{2,2}= \norm{\tilde T_{E}^{k}}_{2,2}.$$
Finally, the result follows by taking the supremum over $\Sigma_0$ of both sides and using Corollary $\ref{cor6.5}$ with $\alpha=1.$
\end{proof}


\section{Positivity of Lyapunov exponents}

\subsection{Introduction and F\"urstenberg's Theorem}
In this section,we prove positivity of the Lyapunov exponents at all energies, for the {\it generalized Anderson model}. We do so by appealing to F\"urstenberg's theorem. Let us first define Lyapunov exponents in this context. Let $d\tilde \mu$ be a probability measure on $SL(2,\mathbb{R})$ which satisfies \begin{equation}\label{eq1sec10}\int \log\norm{M}d\tilde\mu(M)<\infty.\end{equation} Let, $T_1,T_2,\dots$ be $i.i.d$ matrices each with distribution $\mu$. Then, we are interested in the Lyapunov exponent $L\geq 0$, given by $$L=\lim_{n\to\infty}\frac{1}{n}\log\norm{M_n}, \hspace{1cm} \tilde \mu^{\mathbb{Z}_+}- a.s$$ where $M_n=T_n\cdots T_1.$

\begin{thm}[F\"urstenberg's Theorem]\label{FT} Let $\tilde \mu$ be a probability measure on $SL(2,\mathbb{R})$ which satisfies $(\ref{eq1sec10})$. Denote by $G_{\tilde \mu}$ the smallest closed subgroup of $SL(2, \mathbb{R})$ which contains $\supp\tilde \mu.$ 

Assume
\begin{enumerate}[$(i)$]
\item $G_{\tilde \mu}$ is not compact; 

and one of the following: 
\item There is no finite non-empty set $L\subset \mathbb{P}^1$ such that $M(L)=L$ for all $M\in G_{\tilde \mu}$.
\item There is no set $L\subset \mathbb{P}^1$ of cardinality $1$ or $2$ such that $M(L)=L$ for all $M\in G_{\tilde \mu}$.
\end{enumerate}
Then, $L>0.$
\end{thm}


\subsection{F\"urstenberg at all energies}
Now, we state and prove the main theorems of this section.

\begin{thm}\label{posexp}
Suppose that $\#(\supp \nu)\geq 5$. Then, for  the discrete {\it generalized Anderson model}, with $\alpha=2$, we have $L(E)>0$ for all $E\in \mathbb{R}$. 
\end{thm}

\begin{proof}
For every $E\in\mathbb{R}$, the measure $\nu$, as defined in the first section, induces a measure $\tilde \nu$ in $SL(2,\mathbb{R})$ via the map $$\lambda\mapsto\left(\begin{array}{cc}X_1(\lambda)&-1\\
1&0\end{array}\right)\left(\begin{array}{cc}X_0(\lambda)&-1\\
1&0\end{array}\right)=\left(\begin{array}{cc} X_0(\lambda)X_1(\lambda)-1&-X_1(\lambda)\\
X_0(\lambda)&-1\end{array}\right)$$ where $X_0(\lambda)=E-f_0\lambda$ and $X_1(\lambda)=E-f_1\lambda$, $f_0,f_1>0$ and $\lambda\in \supp \nu$. 
So, the random $i.i.d$ matrices for us are the two step transfer matrices. Next, since $\supp \nu$ is uncountable, so is $\supp \tilde \nu$. Let $G_{\tilde \nu}$ be the smallest closed subgroup of $SL(2,\mathbb{R})$ which contains $\supp \tilde \nu.$ 
We will show positivity of the Lyapunov exponent by establishing conditions $(i)$ and $(iii)$ in Theorem $\ref{FT}$. 

To establish $(i)$, by means of contradiction, we suppose that $G_{\tilde \nu}$ is compact. Since $SO(2)$ is a maximal compact subgroup of $SL(2, R)$ we know that every other maximal compact subgroup will be a conjugate of $SO(2)$. Hence, $G_{\tilde \nu}$ will belong to one of these conjugate classes of $SO(2)$. In particular, since $SO(2)$ is abelian, so will every conjugate class of it, and hence $G_{\tilde \nu}$ has to be abelian as well. Below, we will show that $G_{\tilde \nu}$ is not abelian, by means of producing two elements that do not commute, and hence contradicting our above assumption. To this end, let $a,b\in \supp \nu$ such that $a\neq b$, and thus let $M_a, M_b$ be two distinct elements of $G_{\tilde \nu}$, different from $\pm I_2$; that is,

$$M_a=\left(\begin{array}{cc} X_0(a)X_1(a)-1&-X_1(a)\\
X_0(a)&-1\end{array}\right)\, \text{ and }\, M_b=\left(\begin{array}{cc} Y_0(b)Y_1(b)-1&-Y_1(b)\\
Y_0(b)&-1\end{array}\right),$$ where $Y_0(b)=E-f_0b,\, Y_1(b)=E-f_1b$. For ease of notation, from now on, we will suppress the $a$ and $b$ dependence. 

Suppose that $M_aM_b=M_bM_a$, that is,

$$\footnotesize 
\arraycolsep=3pt
\medmuskip = 1mu\left(\begin{array}{cc}(X_0X_1-1)(Y_0Y_1-1)-X_0Y_1& -Y_0(X_0X_1-1)+X_0\\ X_1(Y_0Y_1-1)-Y_1& -Y_0X_1+1\end{array}\right)=\left(\begin{array}{cc}(X_0X_1-1)(Y_0Y_1-1)-Y_0X_1& -X_0(Y_0Y_1-1)+Y_0\\ Y_1(X_0X_1-1)-X_1& -X_0Y_1+1\end{array}\right).$$ So, in particular, we should have

\begin{eqnarray}\label{sec10eq56} \nonumber
Y_0X_0(X_1-Y_1)&=&0\\\nonumber
X_1Y_1(X_0-Y_0)&=&0\\
Y_0X_1-X_0Y_1&=&0\\ \nonumber
\nonumber\end{eqnarray}

Since, by hypothesis, $a\neq b$, it follows that $X_1\neq Y_1$ and $X_0\neq Y_0.$ So, looking at the first two equations, it follows that the only option is that $X_0Y_0=0$ and $X_1Y_1=0.$ Suppose $X_0=0$. Then, from the third equation it would follow that $Y_0X_1=0$. Since, $a\neq b$, in this case, we must have $Y_0\neq 0,$ and thus $X_1=0,$ which is impossible. One argues in a similar way for other cases. Therefore, we conclude that $G_{\tilde \nu}$ is non-compact.

Next, let $L$ be a subset of $\mathbb{P}^1$, where $L=\{v\}$ or $L=\{v,w\}$. For each fixed $E\in \mathbb{R}$, as above, let $G_{\tilde\nu}(E)$ be the smallest closed subgroup of $SL(2,\mathbb{R})$ containing $\supp \tilde \nu.$ To establish $(iii)$, we will break the argument into cases. To this end, we begin by computing the fixed points of the matrix 

$$M_\lambda=\left(\begin{array}{cc} X_0(\lambda)X_1(\lambda)-1&-X_1(\lambda)\\X_0(\lambda)&-1\end{array}\right),$$ where $M_\lambda\neq \pm I_2.$

We compute the fixed points for each of the following three cases:

\begin{enumerate}[$(a)$]
\item $|Tr[M_\lambda]|<2$
\item $|Tr[M_\lambda]|=2$
\item $|Tr[M_\lambda]|>2$
\end{enumerate}
where $Tr[M_\lambda]=X_0(\lambda)X_1(\lambda)-2$ denotes the trace of the matrix $M_\lambda.$ The above three cases correspond to elliptic, parabolic, and hyperbolic systems, respectively. 

If we are in the first case, then it is known that $M_\lambda$ is conjugate to a rotation. Hence, it will have no fixed points in $\mathbb{P}^1$, unless it is a rotation by a multiple of $\pi$, which is equivalent to $M_\lambda$ being equal to $\pm I_2$. Hence, in this case, if $M_\lambda\neq \pm I_2$ then it will not fix any points in $\mathbb{P}^1.$

Next, suppose that $|Tr[M_\lambda]|=2,$ that is, $|X_0(\lambda)X_1(\lambda)-2|=2.$ From here, there are two possibilities:

\begin{enumerate}[$(i)$]
  \item $X_0(\lambda)X_1(\lambda)=0$
  \item $X_0(\lambda)X_1(\lambda)=4$.
 \end{enumerate}
Consider the case where $X_0X_1=0$. That is, our matrices will have the form 
$$M_\lambda=\left(\begin{array}{cc}-1&-X_1(\lambda)\\X_0(\lambda)&-1\end{array}\right).$$

Now, we consider the following two sub-cases: $X_1(\lambda)=0.$ That is, we want to find fixed points of the matrix 
$$M_\lambda=\left(\begin{array}{cc} -1&0\\X_0(\lambda)&-1\end{array}\right).$$ Suppose that $M_\lambda v=v$, for some $v\in \mathbb{P}^1.$ A simple direct calculation shows that $v$ would have to be $[0:1]$, where with $[a:b]$ we denote the equivalence class of the vector $(a,b)^T$.

Next, suppose $X_0=0$. So, we compute the fixed points of the matrix 

$$M_\lambda=\left(\begin{array}{cc} -1&-X_1(\lambda)\\0&-1\end{array}\right).$$ Similarly, as above, one finds that now the fixed point of $M_\lambda$ is $[1:0]$.

Next, we consider case $(ii)$. So, we are interested in computing the fixed point of the matrix 

\begin{equation}\label{eq1058}M_\lambda=\left(\begin{array}{cc}3&-X_1(\lambda)\\X_0(\lambda)&-1\end{array}\right).\end{equation} Let $v\in \mathbb{P}^1$, such that $M_\lambda v=v,$ that is, 

$$\left(\begin{array}{cc}3&-X_1(\lambda)\\X_0(\lambda)&-1\end{array}\right)\left(\begin{array}{c}v_1\\v_2\end{array}\right)=k\left(\begin{array}{c}v_1\\v_2\end{array}\right).$$ This is equivalent to 
\begin{equation}\label{sec10eq57}\frac{3z-X_1(\lambda)}{X_0(\lambda)z-1}=z,\end{equation} where $z=\frac{v_1}{v_2}.$ That is, it is equivalent to finding the fixed points of a M\"obius transformation. Rewriting equation $\ref{sec10eq57})$, we find that the fixed point of $M_\lambda$, in this case, has to be a root of the quadratic equation: $$X_0z^2-4z+X_1=0.$$ Using the fact that $X_0X_1=4$, we can rewrite this equation as $$\left(z-\frac{2}{X_0}\right)^2=0,$$ from where we find that $z=\frac{2}{X_0}$ is the fixed point of this Mobius transformation. In particular, the matrix in $(\ref{eq1058})$ fixes the point $\left[\frac{2}{X_0}: 1\right]$.

Finally, we compute the fixed points of the matrix $M_\lambda$, in the case when $|Tr[M_\lambda]|>2.$ As before, let $v\in\mathbb{P}^1$ be such that $M_\lambda v=v$; that is, 

$$\left(\begin{array}{cc} X_0(\lambda)X_1(\lambda)-1&-X_1(\lambda)\\X_0(\lambda)&-1\end{array}\right)\left(\begin{array}{c}v_1\\v_2\end{array}\right)=k\left(\begin{array}{c}v_1\\v_2\end{array}\right).$$ This is equivalent to $$\frac{(X_0X_1-1)v_1-X_1v_2}{X_0v_1-v_2}=\frac{v_1}{v_2}.$$ Letting $z=\frac{v_1}{v_2}$, we can rewrite this as 
$$\frac{(X_0X_1-1)z-X_1}{X_0z-1}=z.$$ From above, we find that the fixed points are roots of the following quadratic equation 

\begin{equation}\label{eq1059}
z^2-X_1(\lambda)z+\frac{X_1(\lambda)}{X_0(\lambda)}=0.
\end{equation}
In particular, we find that, in this case, the matrix $M_\lambda$ fixes the points $[z_1(\lambda):1]$ and $[z_2(\lambda):1]$ where $z_1(\lambda)$ and $z_2(\lambda)$ are roots of the quadratic equation in $(\ref{eq1059})$.

Next, let $L=\{v\}\subset\mathbb{P}^1$. We show that $L$ cannot be fixed by all elements of $G_{\tilde \nu}(E).$ Suppose that there is some $a\in\supp \nu$, such that, $M_av=v$. From our discussion above, we know that $|Tr[M_a]|=2$ or $|Tr[M_a]|>2$. Suppose that $|Tr[M_a]|=2.$ Then, we are in case $(b)$ above, so we know that $v$ would have to be one of $[0:1], [1:0]$, or $\left[\frac{2}{X_0(a)}:1\right]$. It easily follows from the way these fixed points arise, that if $b\neq a$ and $|Tr[M_b]|=2$, then, $M_b$ cannot fix the same point as $M_a$. Next, suppose that $|Tr[M_b]|>2,$ and suppose further that $b$ is not a solution of 
\begin{equation}\label{quadeq1}\frac{X_1(a)}{X_0(a)}+\frac{X_1(\lambda)}{X_0(\lambda)}-2\frac{X_1(\lambda)}{X_0(a)}=0,\end{equation} 
where above, we are thinking of it as a quadratic equation in $\lambda$. First, observe that $a$ satisfies the quadratic equation in $(\ref{quadeq1})$. So, actually, there is at most one point in $\supp \nu$, different from $a$, that satisfies $(\ref{quadeq1})$. Thus, by our assumption on $\#(\supp \nu)$ (i.e.$\#(\supp \nu)\geq 5$), there is some $b\in\supp \nu$ that is not a root of $(\ref{quadeq1})$. Now, since the fixed points of $M_b$ are $[z_1(b):1]$ and $[z_2(b):1]$ where $z_1(b),z_2(b)$ are roots of the quadratic equation in $(\ref{eq1059})$, for $\lambda=b$, and since zero is not a root of this equation, it is clear that $v$ could only possibly be $\left[\frac{2}{X_0(a)}:1\right]$. But, this in turn would imply that $\frac{2}{X_0(a)}$ is a root of $(\ref{eq1059})$, for $\lambda=b,$ which is not possible by our choice of $b$. Thus, we conclude that $M_b$ cannot fix $\left[\frac{2}{X_0(a)}:1\right]$, which is what we wanted to show.

Next, suppose that $|Tr[M_a]|>2,$ and that $M_a(\{v\})=\{v\}.$ From part $(c)$ above we know that $v$ is $[z_1(a):1]$ where $z_1(a)$ is one of the roots of $(\ref{eq1059})$, for $\lambda=a$. Now, let $b\in\supp\mu$, such that $b\neq a$, and that $b$ is not a solution of \begin{equation}\label{eq1061} z_1(a)X_0(\lambda)X_1(\lambda)-z_1(a)^2X_0(\lambda)-X_1(\lambda)=0,\end{equation}  and is not a solution of $$\frac{X_1(a)}{X_0(a)}+\frac{X_1(\lambda)}{X_0(\lambda)}-2\frac{X_1(a)}{X_0(\lambda)}=0$$when viewed as an equation in $\lambda$. Again, since there are at most three points in $\supp \mu$ that could be roots of the above equations, by our assumption on the cardinality of $\supp \nu$, we note that there are points in $\supp \nu$ that satisfy this condition. Suppose that $M_bv=v.$ As before, there are two subcases we need to consider. If $|Tr[M_b]|>2$, then $v$ would have to be $[z_1(b):1]$ where $z_1(b)$ is a root of $(\ref{eq1059})$ with $\lambda=b.$ So, as a result, we would have $z_1(a)=z_1(b)$.  Since $M_a$ and $M_b$ will each have two fixed points, if it is the case that we also have $z_2(a)=z_2(b)$, then this would imply that $X_1(a)=z_1(a)+z_2(a)=z_1(b)+z_2(b)=X_1(b)$, which would in turn imply that $a=b$, a contradiction. So, the only possibility is that $z_2(a)\neq z_2(b).$  Then, using
\begin{eqnarray*}
z_1(a)+z_2(b)&=&X_1(b)\\
z_1(a)z_2(b)&=&\frac{X_1(b)}{X_0(a)}
\end{eqnarray*}
and substituting for $z_2(b)$ from the first equation into the second, we find that $b$ is a root of $(\ref{eq1061})$, contradicting our choice of $b$.  Now, if $|Tr[M_b]|=2$, then one argues similarly as in one of the cases above.

Now, let $L=\{v,w\}$ be a subset of $\mathbb{P}^1$. As above, we will show that $L$ cannot be fixed by all elements of $G_{\tilde \mu}.$ Suppose there is some $a\in \supp \nu$ such that $M_a(L)=L$. Observe that the only possibility is for $M_a$ to be hyperbolic. This is the case since elliptic matrices have no fixed points, and parabolic matrices have only one fixed point and don't have any periodic points, in particular, they have no points of period two. Then, as before, we know that $v$ must be $[z_1(a):1]$ or $[z_2(a):1]$, where $z_1(a), z_2(a)$ are roots of $(\ref{eq1059})$ with $\lambda=a.$ Now, suppose there is some $b\neq a$ such that $M_bv=v$. For the same reason as before, $M_b$ must be hyperbolic. Then, $v$ must also be one of the following $[z_1(b):1]$ or $[z_2(b):1]$, where $z_1(b),z_2(b)$ are roots of $(\ref{eq1059})$ with $\lambda=b.$ But then we would have $X_1(a)=z_1(a)+z_2(a)=z_1(b)+z_2(b)=X_2(b),$ which would imply that $a=b$, contradicting our choice of $b.$ In this way, we have covered all possible cases, and thus shown that $G_{\tilde \nu}$ cannot fix any subset of $\mathbb{P}^1$ of cardinality one or two. Finally, the result follows by conditions $(i)$ and $(iii)$ in F\"urstenberg's Theorem. 
\end{proof}
\noindent
Below, following the line of arguments in Theorem $\ref{posexp}$, we show that even for the general case, taking the support of $\nu$ to be large enough, we have uniform positivity of the Lyapunov exponents for the {\it generalized Anderson model}.

\begin{thm}
Fix any $\alpha\in\mathbb{Z}_+.$ If $\#(\supp \nu) \geq 35\alpha$, where $\nu$ is as above, then, for the discrete {\it generalized Anderson model}, we have $L(E)>0$ for all $E\in \mathbb{R}$.
\end{thm}

\begin{proof} The random $i.i.d.$ matrices in this case are
$$M_\lambda=\left(\begin{array}{cc} P_\alpha(\lambda)&R_{\alpha-1}(\lambda)\\Q_{\alpha-1}(\lambda)&S_{\alpha-2}(\lambda)\end{array}\right),$$  where $P_m, Q_m, R_m, S_m$ are polynomials of degree $m$, and we only consider $M_\lambda\neq \pm I_2.$ The first part of the proof, including non-commutativity of $G_{\tilde \nu}$, is identical as above, so we skip this part here. However, note that the non-commutativity of $G_{\tilde \nu}$ to be guaranteed it is sufficient to have $\#(\supp \nu)\geq \alpha+1,$ which we do by hypothesis.
As before, we compute the fixed points for each of the following three cases:

\begin{enumerate}[$(a)$]
\item $|Tr[M_\lambda]|<2$
\item $|Tr[M_\lambda]|=2$
\item $|Tr[M_\lambda]|>2$
\end{enumerate}
where $Tr[M_\lambda]=P_\alpha(\lambda)+S_{\alpha-2}(\lambda)$ denotes the trace of the matrix $M_\lambda.$ 
As argued before, we can disregard the first case.

If $|Tr[M_\lambda]|=2$, then $M_\lambda$ has only one fixed point, $[H(\lambda):1]$, where, depending on which sub-case we are in, $H$ can be one of the following, $$H(\lambda)\in \left\{\frac{P_\alpha(\lambda)+1}{Q_{\alpha-1}(\lambda)}, \frac{P_\alpha(\lambda)-1}{Q_{\alpha-1}(\lambda)}, \frac{1}{2}\frac{R_{\alpha-1}(\alpha)}{P_\alpha(\lambda)-1}, \frac{1}{2}\frac{R_{\alpha-1}(\lambda)}{P_\alpha(\lambda)+1}\right\}.$$

If $|Tr[M_\lambda]|>2$, then the fixed points of $M_\lambda$ are $[z_1(\lambda):1]$ and $[z_2(\lambda):1]$, where $z_1, z_2$ are roots of the quadratic equation \begin{equation} \label{eq1067}Q_{\alpha-1}(\lambda)z^2-(P_\alpha(\lambda)-S_{\alpha-2}(\lambda))z-R_{\alpha-1}(\lambda)=0,\end{equation} in case $Q_{\alpha-1}(\lambda)\neq 0$, otherwise, the only fixed double point will be $\left[\frac{R_{\alpha-1}(\lambda)}{P_\alpha(\lambda)-S_{\alpha-2}(\lambda)}:1\right].$

Let $L=\{v\}\subset \mathbb{P}^1$. Let $a\in\supp \mu$ such that $M_av=v.$ We proceed by considering multiple cases. 

Case 1: If $|Tr[M_a]|=2$, then $v=[H(a):1]$, where $H$ is as above. Now, pick $b\neq a$ such that $H(b)\neq H(a)$, and such that $$Q_{\alpha-1}(b)H(a)^2-(P_\alpha(b)-S_{\alpha-2}(b))H(a)-R_{\alpha-1}(b)\neq 0,$$ for all $H$ as above. By our choice of $b$, it follows that $M_b$ cannot fix $v$. Since there are four choices for $H$ and the highest degree is at most $\alpha$ in $(\ref{eqqq01})$ and there are at most $16$ different ways that $H(b)$ can be equal to $H(a)$, we conclude that there are at most $16\alpha+4\alpha$ possible choices for $b\in \supp \nu$ such that $H(b)=H(a)$ or \begin{equation}\label{eqqq01}Q_{\alpha-1}(b)H(a)^2-(P_\alpha(b)-S_{\alpha-2}(b))H(a)-R_{\alpha-1}(b)= 0.\end{equation} Since, by assumption, $\#(\supp \nu)\geq 35\alpha>20\alpha+1$, we can always find some $b\in\supp\nu$, such that it satisfies the above conditions. 

Case 2: If $|Tr[M_a]|>2$, then $v$ is $[z_1(a):1]$, or $[z_2(a):1]$ where $z_1, z_2$ are roots for $(\ref{eq1067})$, when $Q_{\alpha-1}(a)\neq 0,$ or $[K(a):1]$, where $K(a)=\frac{R_{\alpha-1}(a)}{P_\alpha(a)-S_{\alpha-2}(a)},$ otherwise. Pick $b\neq a$ such that \begin{equation}\label{eq1068}Q_{\alpha-1}(a)H(b)^2-(Tr[M_a]-2S_{\alpha-2}(a))H(b)-R_{\alpha-1}(b)\neq 0,\end{equation} for all $H$ as above. Then, from this choice of $b$, it follows that if $M_b$ fixes $v$, it cannot have trace $\pm 2.$ Because the highest degree in $(\ref{eq1068})$ is $2\alpha$, and there are four choices for $H(b)$ we note that there are at most $8\alpha$ points $b\in \supp \nu$ for which the expression in $(\ref{eq1068})$ would possibly not hold. Since $\#(\supp \nu)\geq 35\alpha >20\alpha+8\alpha+1$, we can always find some $b\in \supp \nu$ for which $(\ref{eq1068})$ holds. 
Next, suppose $|Tr[M_b]|>2$, and that $M_bv=v$. As before, $v$ is $[z_1(b):1]$, or $[z_2(b):1]$ where $z_1, z_2$ are roots for $(\ref{eq1067})$, when $Q_{\alpha-1}(b)\neq 0,$ or $[K(b):1]$, where $K(b)=\frac{R_{\alpha-1}(b)}{P_\alpha(b)-S_{\alpha-2}(b)},$ otherwise. If we further pick $b$ such that
 $$Q_{\alpha-1}(a)K(b)^2-(Tr[M_a]-2S_{\alpha-2}(a))K(b)-R_{\alpha-1}(b)\neq 0,$$ and
$$Q_{\alpha-1}(b)(Tr[M_a]-2S_{\alpha-2}(a))\neq Q_{\alpha-1}(a)(Tr[M_b]-2S_{\alpha-2}(b)),$$ and $$z_1(a)(Tr[M_b]-2S_{\alpha-2}(b))-z_1^2(a)Q_{\alpha-1}(b)+R_{\alpha-1}(b)\neq 0,$$ then it actually follows that for this choice of $b$, even in this case $M_b$ cannot fix $v$. Because the highest degree in the first equation is at most $3\alpha$, and in the last two $\alpha$, we note that there are at most $3\alpha+\alpha+\alpha$ points $b\in \supp \nu$, for which the above three conditions could possibly fail. However, since $\#(\supp \nu)\geq 35\alpha\geq20\alpha+8\alpha+5\alpha+1$, we can always find some $b\in\supp \nu$ for which the above three conditions hold. Hence, $L=\{v\}$ cannot be fixed by all of $G_{\tilde \nu}$. 

Finally, let $L=\{v,w\}\subset \mathbb{P}^1$, and suppose that $M_a(L)=L$ for some $a\in\supp \nu$. The only possibility is that $|Tr[M_a]|>2.$ Pick some $b\in\supp \nu$ such that $b\neq a$ and $$Q_{\alpha-1}(b)(Tr[M_a]-2S_{\alpha-2}(a))\neq Q_{\alpha-1}(a)(Tr[M_b]-2S_{\alpha-2}(b)),$$ then since the only possibility is that $Tr[M_b]|>2$ as well, we find that, for this choice of $b$, $M_b$ cannot fix the same set $L$ as $M_a.$ So, we have verified conditions $(i)$ and $(iii)$ in F\"urstenberg's Theorem, and thus proved the claim!
\end{proof}

\begin{rmk} We wish to point out that the lower bound on $\#(\supp \nu)$ in the above theorem is not necessarily optimal.\end{rmk}

\subsection{Existence of exceptional energies -- Examples}

Next we will demonstrate that if one begins by fixing the size of the support for the distribution $\nu$, of the random $i.i.d.$, then it is not possible to obtain uniform positivity of the Lyapunov exponents for all $\alpha$'s. In other words, if you fix the cardinality of the support of $\nu$, then we can always find some $\alpha$ for which the set of exceptional energies will be nonempty. We emphasize that in all of these examples we are primarily concerned with the existence of the exceptional energies and not necessarily with the size of the set of exceptional energies. Let us begin with some warm-up examples. In all of these examples we will take $f_i\equiv 1.$

\begin{exm} Suppose that $\supp \nu=\{0,1\}$. We show that for $\alpha=2m$ or $\alpha=3m,$ where $m\in\mathbb{Z}_+$, the energies $E=0,1$ will be exceptional energies; that is, the Lyapunov exponent will vanish at $0$ and $1$. Note that the random $i.i.d.$ matrices will be the $\alpha$-step transfer matrices 
$$A_0\defeq M_0^\alpha=\left(\begin{array}{cc}E&-1\\1&0\end{array}\right)^\alpha\, \hspace{5mm}\text{and}\,\hspace{5mm} A_1\defeq M_1^\alpha=\left(\begin{array}{cc}E-1&-1\\1&0\end{array}\right)^\alpha$$

For $E=0$, we begin with the observation that $$M_0^2=\left(\begin{array}{cc}0&-1\\1&0\end{array}\right)^2=-I_2,\hspace{5mm} M_1^3=\left(\begin{array}{cc}-1&-1\\1&0\end{array}\right)^3=I_2.$$ Similarly, for $E=1$
$$M_0^3=\left(\begin{array}{cc}1&-1\\1&0\end{array}\right)^3=-I_2,\hspace{5mm} M_1^2=\left(\begin{array}{cc}0&-1\\1&0\end{array}\right)^2=-I_2.$$

So, if $\alpha=2m$ or $3m$, then in each of these cases either the norm of $A_0$ or $A_1$ will be one, so we may disregard it from the random product when computing the Lyapunov exponent; that is, $L(0)$ and $L(1)$ will be computed purely in terms of powers of either $A_0$ or $A_1$, but not both simultaneously. Since, in each of these cases, the trace of both $M_0, M_1$ is strictly less than $2$, we know they will be conjugate to a rotation, so the norms of their powers will remain bounded. Hence, the Lyapunov exponents for $E=0,1$ will vanish. 
\end{exm}

One natural question to ask is if there are exceptional energies as one starts to increase the support of the distribution $\nu$, and if so, for what block size $\alpha$?

The next example answers this question in the special case $\alpha=3,$ and the next proposition gives a more general answer. 

\begin{exm} Let $\supp \nu=\{-1,0,1\}$. We show that for $\alpha=3m$, the energy $E=0$ is an exceptional energy, that is, $L(E)=0.$ The $i.i.d$ matrices for this case will be 

$$A_{-1}\defeq M_{-1}^\alpha=\left(\begin{array}{cc}E+1&-1\\1&0\end{array}\right)^\alpha; \hspace{2mm}A_0\defeq M_0^\alpha=\left(\begin{array}{cc}E&-1\\1&0\end{array}\right)^\alpha\, \hspace{2mm}\text{and}\,\hspace{2mm} A_1\defeq M_1^\alpha=\left(\begin{array}{cc}E-1&-1\\1&0\end{array}\right)^\alpha$$

For $E=0$, we observe that $$M_{-1}^3=\left(\begin{array}{cc}1&-1\\1&0\end{array}\right)^3=-I_2,\hspace{5mm} M_1^3=\left(\begin{array}{cc}-1&-1\\1&0\end{array}\right)^3=I_2.$$ 

Hence, for $\alpha=3m$, we can ignore $A_{-1}$ and $A_1$ in the computation of $L(0)$, as their norms are one and, more importantly, they commute with the other matrices. In other words, we only deal with powers of $A_0$. As before, since $M_0$ has trace strictly less than $2$, the powers of $A_0$ will remain uniformly bounded, resulting in zero Lyapunov exponent. 
\end{exm}

Finally, we state and prove the following proposition.

\begin{prop}\label{prop10.6}
For any given $N\in\mathbb{Z}_+$ there exists a distribution $\nu$ with $\#(\supp \nu)=\frac{N}{2}-1$ or $ \big\lfloor{\frac{N}{2}\big\rfloor}$, if $N$ is even or odd, respectively, and some energy $E$, such that $L(E)=0$ for $\alpha=N$. 
\end{prop}

\begin{proof} Let $N\geq 2$ be a given integer. Let $\nu$ be a measure with support $\big\{-2\cos \left(\frac{2\pi j}{N}\right): j=1,2,\dots, \frac{N}{2}-1\big\},$ or $\big\{-2\cos \left(\frac{2\pi j}{N}\right): j=1,2,\dots, \big\lfloor{\frac{N}{2}\big\rfloor}\big\}$, if $N$ is even or odd, respectively. For $\alpha=N$, the random $i.i.d$ matrices will be the $N-$step transfer matrices. That is, $$ A_j(E)=\left(\begin{array}{cc}E+2\cos \left(\frac{2\pi j}{N}\right)&-1\\1&0\end{array}\right)^N.$$ We claim that $$A_j(0)=\left(\begin{array}{cc}2\cos \left(\frac{2\pi j}{N}\right)&-1\\1&0\end{array}\right)^N=I_2$$ for all $j$, and hence $L(0)=0$, since the norm of the random products of $A_j$'s will be constantly one. To see that $A_j(0)=I_2$, we begin by noting that $|Tr[M_j(0)]|<2$, where $$M_j(0)=\left(\begin{array}{cc}2\cos \left(\frac{2\pi j}{N}\right)&-1\\1&0\end{array}\right).$$ So, for each $j$, the matrix $M_j(0)$ is conjugate to a rotation matrix $$R_{\theta_j}=\left(\begin{array}{cc}\cos \theta_j&-\sin \theta_j\\ \sin\theta_j&\cos\theta_j\end{array}\right),$$ that is, $M_j(0)=P_j^{-1}R_{\theta_j} P_j,$ for some invertible matrix $P_j$. Since, the trace of a matrix is invariant under conjugation, we get $$2\cos \theta_j=Tr[R_{\theta_j}]=Tr[M_j(0)]=2\cos \left(\frac{2\pi j}{N}\right). $$So, in particular, $\theta_j=\frac{2\pi j}{N}+2\pi m$. Then, $A_j(0)=M_j(0)^N=(P_j^{-1}R_{\theta_j} P_j)^N=P_j^{-1}R_{N{\theta_j}}P_j=I_2$, as claimed.  

\end{proof}
\begin{rmk}Note that above if one translates the support of $\nu$ by any real number $\beta$ then in an identical way it follows that $E=\beta$ is an exceptional energy. \end{rmk}
\begin{rmk}
The above proposition is essentially saying that for any given size block $\alpha$, there will exist some distribution $\nu$ for which the Lyapunov exponent of the discrete {\it generalized Anderson model} will vanish for at least one energy. 
\end{rmk}

We get the following immediate corollary. 
\begin{cor}
For any given $N\in\mathbb{Z}_+$ there exists a distribution $\nu$ with $\#(\supp \nu)=\frac{N}{2}$ or $ \big\lfloor{\frac{N}{2}\big\rfloor}+1$, if $N$ is even or odd, respectively, and some energy $E$, such that $L(E)=0$ for any $\alpha=mN$, where $m\in\mathbb{Z}_+$. 
\end{cor}
\begin{proof}
 Let $N\geq 2$ be a given integer. Let $\nu$ be a measure with support $\big\{-2\cos \left(\frac{2\pi j}{N}\right): j=1,2,\dots, \frac{N}{2}-1\big\}\bigcup\{0\},$ or $\big\{-2\cos \left(\frac{2\pi j}{N}\right): j=1,2,\dots, \big\lfloor{\frac{N}{2}\big\rfloor}\big\}\bigcup\{0\}$, if $N$ is even or odd, respectively. From the proof of Proposition $\ref{prop10.6}$ it follows that $M_j(0)^{mN}=I_2$, so to compute $L(0)$ we only need to consider powers of the matrix 
$$\left(\begin{array}{cc}0&-1\\1&0\end{array}\right)^{mN}.$$ Since the norm of any power of this matrix is obviously bounded, it follows that $L(0)=0.$
\end{proof}


\section{Appendix}

The purpose of this section is to outline the main changes in the arguments presented above in the general case where the background potential $V_0$ is no longer identically zero. Since no substantial obstructions arise in this more general case, but rather only notational ones, we have decided to simply give a brief discussion here to convince the reader that all the arguments presented above go through with no significant changes. 

The first change is the form of $H_\omega^L$ the restriction of $H_\omega$ to $\ell^2\left([-\alpha L,\alpha L-1]\cap\mathbb{Z}\right).$ That is, it now takes the form

\begin{equation*}\footnotesize
\arraycolsep=3pt
\medmuskip = 1mu
  H_{\omega}^L=\left(
  \begin{array}{cccccc} 
    V_{\omega}(-\alpha L)+V_0(-\alpha L)& 1& &&&  0 \\ 
    1& V_{\omega}(-\alpha L+1)+V_0(-\alpha L+1)&&& &  0\\ 
 0&1&&&&\\
    \vdots&  &&\ddots&1  & \vdots \\ 
    &&&&V_{\omega}(\alpha L-2)+V_0(\alpha L-2)&1\\
    0& &&\ldots&1& V_{\omega}(\alpha L-1)+V_0(\alpha L-1) \\ 
  \end{array}
  \right).
\end{equation*}

All the results in Section $\ref{sec2.3.}$ hold with no changes in the arguments. 
The next important difference will be when we introduce the Pr\"ufer phase and amplitude, namely, now they will also depend on the background potential. Specifically, as before, let $u_{-L}(\cdot, \omega, E, V_0)$ be the solution to \begin{equation} u(n+1)+u(n-1)+\left(V_\omega(n)+V_0(n)\right)u(n)=Eu(n),\end{equation} with the same set-up and boundary conditions as before. Then, the Pr\"ufer phase $\phi_{-L}(\cdot, \omega, E, V_0)$ and amplitude $R_{-L}(\cdot, \omega, E, V_0)$ will be defined in the same way as before, with the only difference being that now they also depend on $V_0$. Even with the background potential present, one obtains all the results and the same expressions as in Section $\ref{dpv}$, with the only change being in expressions $(\ref{eq4})$ where we have to add the background potential $V_0(n)$; that is;
\begin{equation}
\cot \phi_{-L}\left(n,\{\omega_j\}_{j=-L}^{\lfloor{\frac{n}{\alpha}\rfloor}},E,\{V_0(j)\}_{j=-\alpha L}^{n}\right)+\tan\phi_{-L}\left(n,\{\omega_j\}_{j=-L}^{\lfloor{\frac{n-1}{\alpha}\rfloor}}, E,\{V_0(j)\}_{j=-\alpha L}^{n-1}\right)=E-\omega_{\lfloor{\frac{n}{\alpha}\rfloor}}-V_0(n).
\end{equation}

The next key modification one needs to do is in equation $(\ref{eq024})$. Namely, for $E,\lambda, \theta_{k-1},\theta_k\in\mathbb{R}$ let $u_{-}(\cdot, \theta_{k-1},\lambda, E, V_0)$ be the unique solution of the difference equation \begin{equation} u(n+1)+u(n-1)+\left(\lambda f(n-\alpha k)+V_0(n)\right) u(n)=Eu(n),\end{equation} with the same initial conditions and set up as before. In the same way as before we have the Pr\"ufer phases and amplitudes: $\phi_-(\cdot, \theta_{k-1}, \lambda, E, V_0), \, R_{-}\left(\cdot, \theta_{k-1}, \lambda, E, V_0\right),\, \phi_{+}\left(\cdot, \theta_k, \lambda, E, V_0\right),$ $ R_+(\cdot, \theta_k, \lambda, E, V_0).$ As a result of this change, now the function $\lambda_{E, V_0}:\mathbb{T}^2_\alpha\to[-M,M]$ will be given by $\lambda_{E, V_0}(x,y)\defeq\lambda(x,y,E,V_0),$ where $\lambda(x,y,E,V_0)$ is defined in exactly the same way as before. Specifically, given any $x,y\in\mathbb{T}_\alpha$, if there is some coupling constant $\lambda\in[-M,M]$ such that $\phi_{-}(\alpha k+\alpha-1, y, \lambda, E, V_0(\alpha k),\dots, V_0(\alpha k+\alpha-1))=x$ or $\phi_{+}(\alpha k-1, x,\lambda, E, V_0(\alpha k),\dots, V_0(\alpha k+\alpha-1)=y$ we set $\lambda(y,x,E, V_0(\alpha k),\dots, V_0(\alpha k+\alpha-1))=\lambda.$
Then, the arguments all the way to Section $\ref{CJ}$ go through in exactly the same way as before. The family of integral operators that we define in Section $\ref{IOF}$ need to be modified slightly when the background potential $V_0$ is present. The change is indeed minor, all one needs to do is essentially replace the corresponding $\lambda(x,y,E)$ with $\lambda(y,x,E, V_0(\alpha k),\dots, V_0(\alpha k+\alpha-1)),$ in Definitions $\ref{defi3.7}$ and $\ref{defi3.8.}$. Then, there are essentially no changes in the arguments one needs to make until Section $\ref{sec7}$. There, the main change one needs to do is in Theorem $\ref{thm6.1}$. Specifically, now one needs to prove that the real valued map $\displaystyle (\bar x,E)\mapsto \norm{\tilde T_{\bar x,E,\alpha}^{k}}_{2,2}$ is continuous on $\big[-\|V_0\|_\infty,\|V_0\|_\infty\big]^\alpha\times\Sigma_0$. This is also done in an identical way as before. Below wel provide the modified statement of Lemma $\ref{lemma6.2}$, and the reader should convince himself that the rest of the results in Section $\ref{sec7}$ follow in an almost identical way as before. 

\begin{prop}Suppose $\bar z_n\defeq(E_n,x_1^n,\dots, x_\alpha^n)\to(E,x_1,\dots, x_\alpha)\defeq \bar z$ as $n\to \infty.$ Then $$\lambda(y,x,\bar z_n)\to \lambda(y,x,\bar z),$$ whenever $\lambda(y,x,\bar z)$ exists. 
\end{prop}

The proof o this lemma and the other results in Section $\ref{sec7}$ follow in an identical manner where one replaces $E_n$ by $\bar z_n$, and $E$ by $\bar z$.

\vskip 20pt
\noindent {\sc Acknowledgments:}  I would like to thank my PhD advisor, David Damanik, for introducing me to this problem and also for his guidance and many invaluable discussions. I am also grateful to Zhenghe Zhang for helpful discussions regarding the proof of Theorem $\ref{posexp}$. I would also like to thank Jorge Acosta, Tom VandenBoom, and Vitalii Gerbutz, for useful discussions regarding Proposition $\ref{prop10.6}$.

\vskip12pt

\noindent\textsc{Department of Mathematics, Rice University, Houston, TX 77005, U.S.A}\\
{\it E-mail address}: \ttfamily{valmir.bucaj@rice.edu}

\end{document}